\documentclass[11pt, twoside]{amsart}
\usepackage{amsmath}
\usepackage{amssymb}
\usepackage{amsfonts}
\usepackage[latin1]{inputenc}
\numberwithin{equation}{section}

\newtheorem{theorem}{Theorem}[section]
\newtheorem{lemma}[theorem]{Lemma}
\newtheorem{definition}[theorem]{Definition}

\newtheorem{proposition}[theorem]{Proposition}
\newtheorem{corollary}[theorem]{Corollary}
\newtheorem{fact}[theorem]{Fact}
\newtheorem{claim}[theorem]{Claim}

\newtheorem{subclaim}[theorem]{Subclaim}

\theoremstyle{remark}  
\newtheorem*{remark}{Remark}

\newtheorem*{explanation}{Explanation}
\newtheorem*{notation}{Notation}

\newcommand{\nc}{\newcommand}
\nc{\wittt}{}
\DeclareMathOperator{\tcf}{tcf} 
\nc{\prodt}{{\textstyle \prod}} 
\nc{\nco}{\DeclareMathOperator}
\nco{\dom}{dom}
\nco{\card}{card}
\nco{\lh}{lh}
\nco{\lgg}{lg}
\nco{\rge}{range}
\nco{\cf}{cf}
\nco{\nex}{next}
\nc{\uhr}{\restriction}
\nco{\supt}{supt}
\nco{\supp}{supp}
\nco{\Lim}{lim}
\nco{\Leb}{leb}
\nco{\modd}{mod}
\nco{\invariant}{inv}

\nc{\potom}{\ensuremath{{\mathcal P}(\omega)}}
\nc{\potinf}{\ensuremath{[\omega]^\omega}}
\nc{\potfin}{\ensuremath{[\omega]^{<\omega}}}
\nc{\inn}{\ensuremath{{\omega^{\uparrow \omega}}}}
\nc{\hoch}{^{<\omega}}
\nc{\hocho}{^{\omega}}
\nc{\tree}[1]{{[} #1 {]}_0}
\nc{\tre}[2]{ {#1}_{#2}}
   
\nc{\prooff}[1]{{\bf Proof} of #1:}
\nc{\proofend}{\makebox{} \hfill $\square$ \\}
\nc{\proofendof}[1]{\makebox{} \hfill $\square_{\rm #1}$ \\}
\nc{\beq}{\begin{eqnarray*}}
\nc{\eeq}{\end{eqnarray*}}
\nc{\bde}{\begin{list}}
\nc{\ede}{\end{list}}

\newenvironment{myrules}
{\begin{list}{}
{
 \setlength{\leftmargin}{0.3in}
 \setlength{\labelwidth}{0.8cm}
 \setlength{\labelsep}{0.1in} 
 \setlength{\parsep}{0.5ex plus 0.2ex minus 0.1 ex}
 \setlength{\itemsep}{0.3ex plus 0.2 ex minus 0ex}
}}{\end{list}}

\newenvironment{myrules1}
{\begin{list}{}
{
 \setlength{\leftmargin}{0.3in} 
 \setlength{\labelwidth}{0.8cm}
 \setlength{\labelsep}{0.1in}
 \setlength{\parsep}{0.5ex plus 0.2ex minus 0.1 ex}
 \setlength{\itemsep}{0.5ex plus 0.2 ex minus 0ex}
}}{\end{list}}

\newcount\skewfactor

\def\mathunderaccent#1#2 {\let\theaccent#1\skewfactor#2
\mathpalette\putaccentunder}
\def\putaccentunder#1#2{\oalign{$#1#2$\crcr\hidewidth
\vbox to.2ex{\hbox{$#1\skew\skewfactor\theaccent{}$}\vss}\hidewidth}}
\def\name{\mathunderaccent\tilde-3 }
\def\Name{\mathunderaccent\widetilde-3 }

\nc{\al}{$\alpha$\  }
\nc{\om}{\omega}
\nc{\omm}{\ensuremath{\omega_1}}
\nc{\ep}{\varepsilon}
\nc{\tk}{\tilde{K}}
\nc{\concat}{\,\hat{} \,}   
\nc{\force}{\Vdash}
\nc{\fb}{f_{\bar{M}}}
\nc{\such}{\, | \,}   

\nc{\meager}{\ensuremath{{\mathcal M}}}
\nc{\lebesgue}{\ensuremath{{\mathcal N}}}
\nc{\nulll}{\ensuremath{{\mathcal N}}}
\nc{\ksigma}{\ensuremath{{\bf K}_\sigma}}
\nc{\ideal}{\ensuremath{{\mathcal I}}}
\nc{\ga}{\ensuremath{\frak a}}
\nc{\AAA}{{\mathcal A}}   
\nc{\gc}{\ensuremath{\frak c}}
\nc{\gs}{\ensuremath{\frak s}}
\nc{\gh}{\ensuremath{\frak h}}
\nc{\gd}{\ensuremath{\frak d}}
\nc{\gb}{\ensuremath{\frak b}}
\nc{\gro}{\ensuremath{\frak g}}
\nc{\gu}{\ensuremath{\frak u}} 
\nc{\gr}{\ensuremath{\frak r}} 
\nc{\fff}{\ensuremath{\frak f}}
\nc{\gm}{\ensuremath{\mathfrak{mcf}}}
\nc{\gge}{\ensuremath{\mathfrak e}}
\nc{\cfupro}{\ensuremath{\cf(\upro)}}
\nc{\cfvpro}{\ensuremath{\cf(\vpro)}}
\nc{\gp}{\ensuremath{\frak p}}
\nc{\gk}{\ensuremath{\frak k}}

\nc{\add}[1]{\mbox{\ensuremath{{\rm add}(#1)}}}
\nc{\cov}[1]{\mbox{\ensuremath{{\rm cov}(#1)}}}
\nc{\unif}[1]{\mbox{\ensuremath{{\rm unif}(#1)}}}
\nc{\cof}[1]{{\mbox{\ensuremath{\rm cof}(#1)}}}

\nc{\addd}[2]{\mbox{\ensuremath{{\rm add}^{#1}(#2)}}}   
\nc{\covv}[2]{\mbox{\ensuremath{{\rm cov}^{#1}(#2)}}}   
\nc{\uniff}[2]{\mbox{\ensuremath{{\rm unif}^{#1}(#2)}}} 
\nc{\coff}[2]{{\mbox{\ensuremath{\rm cof}^{#1}(#2)}}}

\nc{\cd}{Cicho\'n's Diagram}

\nc{\MA}{\mbox{\rm MA}}
\nc{\GCH}{\mbox{\rm GCH}}
\nc{\CH}{\mbox{\rm CH}}
\nc{\zfc}{\mbox{\rm ZFC}}
\nc{\sch}{\mbox{\rm SCH}} 
\nc{\ZF}{\mbox{\rm ZF}}
\nc{\NCF}{\mbox{\rm NCF}} 
\nc{\FD}{\mbox{\rm FD}}   

\nc{\av}[1]{{\rm Av}_{#1}}
\nc{\eps}{\varepsilon}
\nc{\n}{{\bf n}}                 
\nc{\m}{{\bf m}}

\nc{\marginparr}[1]{}

\newcommand{\nothing}[1]{}
\nc{\betta}{\beta}

\begin{document}

\title[Changing Cardinal Characteristics]
{Changing Cardinal Characteristics Without Changing
$\omega$-Sequences or Cofinalities}

\author{Heike Mildenberger and Saharon Shelah}

\thanks{The first author was partially supported by a Lise 
Meitner Fellowship of the State of North Rhine Westphalia} 

\thanks{The second author's research 
was partially supported by the ``Israel Science
Foundation'', administered by the Israel Academy of Science and Humanities.
This is  the second author's publication no.\ 684}


\nothing{
\date{May 29, 1998, rev.\ in June, \ref{1.4} and \ref{1.5} added on July 2,
Section 7 revised on September 11, 1998, Sections 2 to 4
revised on October 2, 1998}}
\date{December 6, 1998}

\address{Heike Mildenberger,
Mathematisches Institut der Universit\"at Bonn,
Beringstr.~1,
53115~Bonn, Germany}

\email{heike@math.uni-bonn.de}

\address{
Saharon Shelah,
Mathematical Institute,
The Hebrew University of Jerusalem,
 Givat Ram,
91904~Jerusalem, Israel
}
 
\email{shelah@math.huji.ac.il}

\begin{abstract}
We show: There are pairs of universes $V_1 \subseteq V_2$
and there is a notion of forcing $P \in V_1$
such that the change mentioned in the title occurs when going from
$V_1[G]$ to $V_2[G]$ for a $P$-generic filter $G$ over $V_2$. 
We use forcing iterations with 
partial memories.
Moreover, we implement
highly transitive automorphism groups into the forcing orders.
\end{abstract}

\maketitle
\tableofcontents

\setcounter{section}{-1}
\section[Introduction]{Introduction}

In \cite{Mi2} it is shown that some cardinal characteristics can be
changed without changing $\omega$-sequences or cardinalities, that is 
we can have two models $V_1 \subseteq V_2$ of \zfc\ such that 
 $(^\omega V_1)^{V_2} \subseteq V_1$ and such that $V_1$ and
$V_2$ have
the same cardinalities and such that, e.g., $\gd^{V_2} <
\gd^{V_1}$ ($\gd$ is the dominating number, the
minimum size of a subset ${\mathcal D} \subseteq \omega^\omega$ such that
every function $f \in \omega^\omega$ is eventually 
dominated by some member of ${\mathcal D}$). 
Since in such a situation the covering theorem for
$(V_1,V_2)$ fails, there is consistency strength of at least a measurable
cardinal. In \cite{Mi2} a change of a cofinality of a regular cardinal in
$V_1$ was the main step when changing all the entries of
Cicho\'{n}'s Diagram (for information on cardinal characteristics and 
Cicho\'{n}'s Diagram  see e.g.\ \cite{Blasshandbook,BJ,vanDouwen,Vaughan})
without changing cardinalities or the reals.
In this work we show that we do not need to change cofinalities in
order to change $\gb$, $\cov{\meager}$, $\cov{\lebesgue}$, $\unif{\meager}$ or 
$\unif{\lebesgue}$ and both additivities without changing
cardinalities or the reals. These are all
entries of Cicho\'{n}'s Diagram that are not norms of transitive relations.
In order to cover all these cases we use two different procedures.

\smallskip

In Section 1, we show how to change $\gb$, $\unif{\meager}$ and
$\cov{\lebesgue}$ and both additivities
starting from a bare set-theoretic situation.
We use an iteration with partial memory.

In \cite{Mi2} it is shown that $\gd$, 
$\cof{\meager}$ and $\cof{\lebesgue}$ cannot be changed if their values in 
$V_1$ are 
regular in $V_2$ and if $V_1$ and $V_2$ have the same cardinalities.
At the end of Section~\ref{S1}, we shall
show that if $V_1$ and $V_2$ have the same cofinalities, then
 these characteristics (and some more, whose definition exhibits a 
certain syntax) cannot be changed either when starting 
from a singular value in $V_1$.

\smallskip

In Sections~\ref{S2} to \ref{S5}, we 
show  how to change $\unif{\lebesgue}$.
We work with partial random forcing as in 
\cite{Shelah592, Shelah619}, however, as we need 
special  instances of the methods presented there, we 
(try to) 
make our present work self-contained. We include some comments 
on the connections to \cite{Shelah592, Shelah619} and give references to 
items we use almost literally, so that the reader may also read these.
In Section~\ref{S6} we shall present a variation of the techniques for
a case with countable cofinality.

\smallskip

In Section~\ref{S7}, we show how to obtain the  set-theoretic assumptions 
made in Theorems \ref{1.1} and \ref{2.1} from Gitik's work in 
\cite{Gitik, Gitik89}.

\smallskip
The authors would like to thank Andreas Blass for reading a section 
and  commenting.

\smallskip

\begin{notation} 
Our notation is fairly standard, see \cite{Jech, Kunen}. 
However, we adopt the Jerusalem convention that the stronger forcing 
condition is the larger one. 
We often use $V^P$ for $V[G]$, where
$G$  is any $P$-generic filter over $V$.
For two forcing notions $P, Q$ we write
$P \lessdot Q$ if $P$ is a complete suborder of $Q$.
A forcing notion $P$ is called $\sigma$-linked if $P = \bigcup_{n \in \omega}
P_n$ such that each $P_n$ is linked, that is any two $p,q \in P_n$ 
are compatible. Martin's axiom for less than $\lambda$
dense subsets of a $\sigma$-linked partial order is denoted by
$\MA_{<\lambda}(\sigma\mbox{-linked})$. We speak of $\omega^\omega$, 
the set of all functions from $\omega$ to $\omega$, as the reals. 
For $f,g \in \omega^\omega$ we write
$f \leq^\ast g$ if $\exists n \: \forall k \geq n \: f(k) \leq g(k)$.
The ideal of Lebesgue null sets is denoted by $\lebesgue$,
and the ideal of meagre sets is denoted by $\meager$. 
The bounding number, $\gb$, is the smallest size of a subset 
$B\subseteq \omega^\omega$
such that for any $f \in \omega^\omega$ there is some 
$b \in B$ such that $b \not\leq^\ast f$. Let $\mathcal I$ be 
an ideal on the reals.
The uniformity of ${\mathcal I}\subseteq {\mathcal P}(\omega)$,
$\unif{\ideal}$, is the
smallest size of a subset of the reals that is not a member of $\mathcal I$.
The covering number of $\mathcal I$, $\cov{\ideal}$, is the smallest size of
a subfamily of $\mathcal I$ whose union covers the reals. 
The additivity of $\mathcal I$, $\add{\ideal}$, is  the smallest size
of a subset of $\mathcal I$ whose union is not in $\mathcal I$.
\end{notation}

\section{Changing the Uniformity of Category}\label{S1}

In this section, we show how to change $\unif{\meager}$.
Since $\add{\meager} \leq \gb \leq \unif{\meager}$ and $\add{\lebesgue} 
\leq \cov{\lebesgue} \leq \unif{\meager}$
(for proofs of these
inequalities, see
\cite{Fremlin}, e.g.), and 
in the beginning, that is in $V_1[G]$,
everything is large because of an instance of Martin's axiom,
the other four mentioned characteristics drop as well. 

\begin{theorem}\label{1.1}
Assume that we have 
\begin{myrules}
\item[a)] $V_1 \subseteq V_2$, both models of $\zfc$, $(^\omega V_1)^{V_2}
\subseteq V_1$, 

\item[b)] $\mu$ is a cardinal in $V_2$,
$C \subseteq \mu$, $C \in V_2$, ${\mathcal I} \in V_2$ is an
$\aleph_1$-complete proper ideal on ${\mathcal P}(C)$,

\item[c)] $\exists \lambda \leq \mu$ such that
 $\forall B \in V_1$, if $V_1 \models |B| < \lambda$, 
then $B \cap C \in {\mathcal I}$,

\item[d)] $V_1 \models \lambda > \aleph_0$ and $\lambda$ is regular.
\end{myrules}

Then for some $P$
\begin{myrules}
\item[$\alpha$)] $V_1 \models P$ is a finite support iteration of
$\sigma$-linked forcing notions, and the 
cardinality of $P$ is $\mu^{<\lambda}$,
\item[$\beta$)]  $P$ is c.c.c.\ in $V_2$.
\end{myrules}
For $G \subset P$ generic over $V_2$ we have
\begin{myrules}
\item[$\gamma$)]
 $(^\om V_1[G])^{V_2[G]} \subseteq V_1[G]$,
\item[$\delta$)]
 $V_1[G]$ and $V_2[G]$ have the same cardinals if $V_1$ and $V_2$ have,
\item[$\varepsilon$)]
 $V_1[G]$ and $V_2[G]$ have the same cofinality function if $V_1$ and $V_2$
have,
\item[$\zeta$)] $V_1[G] \models \MA_{<\lambda}(\sigma\mbox{-linked})$,
\item[$\eta$)] in $V_2[G]$ there is $\langle r_i \such i \in C \rangle$,
$r_i \in (^\omega 2)^{V_1[G]}=(^\omega 2)^{V_1[G]}$, such that $\forall s \in 
(^\omega 2)^{V_1[G]} \; \exists B \subseteq \mu$, $B \in V_1$, $|B|^{V_1} < 
\lambda$ 
(so $C \cap B \in {\mathcal I}$) $\forall i \in C \setminus B$ , $r_i$ 
is Cohen over $V_2[s]$.
\end{myrules}
\end{theorem}
\proof
In $V_1$ we build a finite support iteration $$\langle P_i, \Name{Q}_j
\such j< \alpha^\ast, i \leq \alpha^\ast  \rangle$$
 of length
$\alpha^\ast = \mu + \mu^{< \lambda}$ as follows. 
For $\beta < \mu$ we let $Q_\beta = (^{<\omega}2, 
\vartriangleleft)$, 
the Cohen  forcing. 

\smallskip

For $\beta < \mu^{<\lambda}$ we shall choose
$\Name{Q}_{\mu +\beta}$ such that it is a 
name built from only part of  $P_{\mu +\beta}$.
We first need some definitions in order to specify good parts of the past.
This forcing technique has also been applied in \cite{Shelah546}, \cite{Shelah592},
\cite{Shelah619}  and their predecessors  and in \cite{ShelahThomas}.
The part 
\cite[3.3 to 3.7]{ShelahThomas} contains some lemmas showing that
there are complete embeddings from specified suborders of the iteration
that are not just initial segments. 
The organisation of our forcing will 
be slightly different from that in \cite{ShelahThomas} inasmuch as
we have the initial Cohen part here at once.
%

\smallskip

The support of a condition $p \in P_\beta$ is $\supt(p) =\{ \gamma \in \beta 
\such p(\gamma) \neq 1_{\Name{Q}_\gamma} \}$, 
where $1_{\Name{Q}_\gamma}$ is a name
for the weakest element in $\Name{Q}_{\gamma}$. In addition to having
finite supports we shall require that the supports hereditarily stem only
from a part of the ``past'' $P_\beta$. These parts of the
past can be called memories.

\smallskip

First we explain how to choose sequences $\langle a_\beta 
\such \beta \in \mu^{<\lambda} \rangle$ which will allow us to
define suitable memories.
  Given a sequence 
$\langle a_\beta \such \beta \in \mu^{<\lambda} \rangle = \bar{a}$ of 
subsets of an ordinal, we say $c$ is 
$\bar{a}$-closed, if 
$$
c \subseteq \alpha^\ast \mbox{ and } \forall \beta \in c \:\; a_\beta 
\subseteq c.
$$

We regard $\mu^{<\lambda}$ as an ordinal and as a set of sequences
of length less than $\lambda$. The set of all subsets of a set $A$
of size less than $\lambda$ is denoted by $[A]^{<\lambda}$.
For $x \in \mu^{<\lambda}$ we can also regard $x $ as a function from 
some ordinal  less than $\lambda$ to $\mu$ and then write $\rge(x)$
for its range, which is a subsets of $\mu$. This will be used for 
referring to  a part of the Cohen reals. 

\smallskip

 We show that there is some $\langle a_\beta \such
\beta < \mu^{<\lambda} \rangle$ such that 
\begin{enumerate}
\item $\forall b \in [\mu^{<\lambda} ]^{<\lambda} \: \exists \beta \; 
  b \subseteq a_\beta$, 
\item $a_\beta \subseteq \beta$, 
\item $|a_\beta| < \lambda$,
\item $\gamma \in a_\beta \rightarrow a_\gamma \subseteq a_\beta$ (i.e.\ 
each $a_\beta$ is $\bar{a}$-closed).
\end{enumerate}

This can be seen as follows:
Let $\langle b_\beta \such \beta \in \mu^{<\lambda}\rangle$ enumerate 
$[\mu^{<\lambda}]^{<\lambda}$, where $b_\beta \subseteq \beta$. By induction on 
$\beta$ we now choose $a_\beta$. Suppose $a_\gamma$ is chosen for 
$\gamma < \beta$. Then we set
\begin{eqnarray*}
a_\beta^1 &=&\bigcup_{j \in b_\beta} a_j \cup b_\beta,\\
a_\beta^{n+1} &=& \bigcup_{j \in a_\beta^n} a_j \cup a_\beta^n,\\
a_\beta &=& \bigcup_{n \in \omega} a_\beta^n.
\end{eqnarray*}
This is still in $[\mu^{<\lambda}]^{<\lambda}$ because $\lambda $ is regular and
$\cf(\lambda) > \aleph_0$.
Now it is easy to see that $\bar{a}$ fulfils 1.\ to 4., and 
we fix such a sequence. 

\nc{\pluss}{\oplus}
In order to take care of the initial Cohen part, we need shifts and
write $\mu \pluss a_\beta$ for $\{\mu + \gamma \such \gamma \in a_\beta \}$.

For each $\beta \in \mu^{<\lambda}$ we define a 
suborder $P^*_{\mu \pluss a_\beta}
$ of $P_{\mu +\beta}$ inductively by
\begin{equation*}
\begin{split}
P^*_{\mu \pluss a_\beta} = \{ p \in P_{\mu +\beta} \such &
\supt(p) \cap \mu \subseteq \bigcup\{\rge(x) \such x \in a_\beta \} \wedge 
\\ &
\supt(p) \cap [\mu,\mu+\mu^{<\lambda}) \subseteq \mu \pluss a_\beta  \wedge
\\ &
 \forall \gamma \in \supt(p) \cap [\mu,\mu + 
\mu^{<\lambda}) \; p(\gamma) \mbox{ is a } P^*_{\mu \pluss a_\gamma}\mbox{-name}\}.
\end{split}\end{equation*}

\smallskip

If $b \subseteq \alpha \leq\mu^{<\lambda}$ then 
$p\restriction (\bigcup\{\rge(x)\such x \in b\}\cup
 \mu \pluss b))$ denotes the $
\mu + \alpha$-sequence defined by 
\begin{equation*}\begin{split}
&
(p\restriction (\bigcup\{\rge(x)\such x \in b\}\cup
 \mu \pluss b)))(\gamma) \\
&= \left\{ \begin{array}{ll}
p(\gamma) & \mbox{ if } \gamma \in (\bigcup\{\rge(x)\such x \in b\}\cup
 \mu \pluss b)) , \\
1_{\Name{Q}_\gamma} & \mbox{ else.}
\end{array}
\right.
\end{split}\end{equation*}

Now we have for all $\alpha \in \mu^{<\lambda}$:
If $b \subseteq \alpha$ is $\bar{a}$-closed, then $P^*_{\mu \pluss b}
\lessdot P_{\mu +\alpha}$.
If $p \in P_{\mu +\alpha}$, then $(p \restriction 
(\bigcup\{\rge(x)\such x \in b\}\cup \mu \pluss b))) \in P^*_{\mu \pluss b}$
 and for $q \geq p \restriction (\bigcup\{\rge(x)\such x \in b\}\cup
 \mu \pluss b)$ (in the Jerusalem notation)  we have that
 $q \cup p \restriction (\alpha \setminus (\bigcup\{\rge(x)\such x \in b\}\cup
 \mu \pluss b)) \in P_{\mu +\alpha}$.
For proofs, see \cite{ShelahThomas}.

\smallskip

We choose $\Name{Q}_{\mu +\beta}$ such that $|\dom(\Name{Q}_{\mu +\beta})| 
< \lambda$,
$\Name{Q}_{\mu +\beta}$ is a $P^*_{\mu \pluss a_\beta}$-name, $1 \Vdash_{P^*_{\mu \pluss
a_\beta}} \mbox{``}\Name{Q}_{\mu +\beta} \mbox{ is $\sigma$-linked''}$, 
and with some bookkeeping such that
$\Name{Q}_{\mu +\beta}$ ranges cofinally often 
over all $P^*_{\mu \pluss a_\gamma}$-names for $\sigma$-linked forcings 
for every $\gamma \in \mu^{<\lambda}$.
In order to allow 
such a bookkeeping, we assume that 
$\forall b \in [\mu^{<\lambda}]^{<\lambda} \;
\exists ^{\mu^{<\lambda}} \beta \; b \subseteq a_\beta$, 
which can easily be reached 
by starting with suitable $\langle b_\beta \such \beta 
\in \mu^{<\lambda} \rangle$.

\smallskip

Now we are in a position to check all the items of the theorem:

\smallskip

$\alpha$) follows immediately from our definition of $P$.

\smallskip

$\beta$) If $P = \bigcup_{n \in \omega} P_n$ 
witnesses $\sigma$-linkedness in $V_1$ then it does so in $V_2$ as well.
Thus  in $V_2$, $P$ is a finite support iteration of $\sigma$-linked
forcing notions and hence c.c.c.

\smallskip

$\gamma$) $(^\omega V_1[G])^{V_2[G]} \subseteq V_1[G]$  follows from
$(^\omega V_1)^{V_2} \subseteq V_1$ and the countable chain condition
of $P$ in$V_2$.
(There are also proofs in \cite[\S 37]{Jech} and more explicit in
 \cite{Cummings}.)

\smallskip

$\delta$) and $\varepsilon$) $V_i$ and $V_i[G]$ have the same cofinalities.

\smallskip

$\zeta$) Let $Q$ be in $V_1[G]$ be a $\sigma$-linked notion of forcing 
such that $Q \subseteq \lambda' < \lambda$. Let 
${\mathcal D} = \{ D_\alpha \such \alpha < \lambda' \}$ be a set
of dense sets in $Q$. Since the supports are 
finite and since we have c.c.c., 
there is some $A \subseteq \mu + \mu^{<\lambda}$ of size less than $\lambda$
such that there is a name for $(Q,{\mathcal D})$ that contains
only conditions whose support is in $A$.
Then we take $\alpha \in \mu^{<\lambda}$ such that

\begin{equation*}
\begin{split}
x= \bigcup \{ \rge(x) \such x \in a_\alpha \} & \supseteq A \cap \mu \mbox{ and }
 \\
y= \mu \pluss a_\alpha & \supseteq A \cap [\mu,\mu + \mu^{<\lambda}).
\end{split}
\end{equation*}

and have that 
that $
{\mathcal D}, Q \in V_1^{P^*_{\mu \pluss a_\alpha}}$. Hence a $Q$-generic $G \subseteq Q$ is
added at some stage in our iteration.

\smallskip

$\eta$) Let $\langle r_i \such i \in \mu\rangle$ be the Cohen reals added by
$P_\mu$. We show that  $\{ r_i \such i \in C \}$
is as claimed. Let $s \in (2^\omega)^{V_1[G]}$. Say $s$ was 
added by forcing with $\Name{Q}_{\mu + \beta}$ (the case 
when $s$ was added before stage $\mu$ is similar), a $P_{\mu \pluss
a_\beta}$-name.
We take $B = a_\beta$. Then $B \in V_1$, $B \subseteq \mu$, 
and $|B|^{V_1} < \lambda$.  As $C \cap B \in {\mathcal I}$, we have 
$C \setminus B \neq \emptyset$. For $i \in C \setminus B$ 
$r_i$ is Cohen over $V_1[s]$. Proof:  
For $Q_i =(^{< \omega} 2, \vartriangleleft)$ we have
$$
Q_i \ast P^*_{\mu \pluss a_\beta}
= Q_i\times P^*_{\mu \pluss a_\beta}.$$

\begin{remark} This equation is very crucial: Note that
there is ``no time-dependence'', i.e.\ the location of
$i$ in $\mu + \mu^{<\lambda}$ as compared to the location of
$ x \cup y$ does not have any influence.
Neither $Q_i$ nor $ P^*_{\mu \pluss a_\beta}$
is the ``later'' forcing, because neither of them is
influenced by the extension performed by the other.
 All the work with the partial memory 
was done in order to get this equation. Counting cardinalities of
unions of supports of conditions appearing in nice names 
seems not to suffice for it.\end{remark}

The analogue of the crucial equation is true for the
subforcing of $P^*_{\mu \pluss a_\beta}$ that hat $s$ as a generic.
 Now in 
product forcing, the factors commute, hence we have
 $V_1[r_i][s] = V_1[s][r_i]$.
\proofendof{\ref{1.1}}

Putting things together we get
\begin{corollary}\label{1.2}
(1) The following are equiconsistent (even (B) $\Rightarrow$ (A), (A) 
$\Rightarrow$ (B) in some c.c.c.\ forcing extension):
\begin{myrules}
\item[(A)($\alpha$)] there are $V_1$, $V_2$, $\mu$, 
                   $\theta$, $\lambda$, $\sigma$, $C$, such that: \\
                   $V_1 \subseteq V_2$, \\
                   $V_1 \models \lambda \mbox{ regular} > \aleph_0$,\\
                   $(^\omega V_1)^{V_2} \subseteq V_1$,\\
                   $\mu \geq \theta$, $\mu  \geq \lambda > \sigma \geq 
                   \aleph_1$,\\
                   $C \subseteq \mu$,\\
                   $|C|^{V_2} = \theta$,\\
                   $\forall B \in V_1 \: (|B|^{V_1} < \lambda 
               \rightarrow |B \cap C|^{V_2} < \sigma)$
  \item [\phantom{(A)}($\beta$)] $V_1$ and $V_2$ have the same cardinals.
  \item [\phantom{(A)}($\gamma$)] $V_1$ and $V_2$ have 
the same cofinality function on ordinals.
  
\item[(B)($\alpha$)] like (A)($\alpha$) but in addition\\
   ($\ast_1$) $V_1 \models \MA_{<\lambda}(\sigma\mbox{-linked})$\\
   ($\ast_2$) in $V_2$ there are $\langle r_i \such i \in C \rangle$,
   $r_i \in 2^\omega$ and a submodel $V$
    such that $\forall s \in 2^\omega \: \exists B 
   \in [C]^{<\sigma} $ such that $\langle r_i \such i \in C 
   \setminus B \rangle$ is Cohen over $V[s]$.
\item[\phantom{(B)}($\beta$)] as ($\beta$) above.
\item[\phantom{(B)}($\gamma$)] as ($\gamma$) above.
\end{myrules}

(2) We can leave out ($\beta$) or (($\beta$) and ($\gamma$)) in both 
(A) and (B).

(3) If we strengthen (A)($\alpha$) by adding\\
 $(^{\omega_1} V_1)^{V_2}
\subseteq V_1$, then we can get $\MA_{<\lambda}(ccc)$ in (B).

\end{corollary}

\proof (A) is as the premise of \ref{1.1} with ${\mathcal I} =
\{ C'\subset C \such C' \in V_2, |C'|^{V_2} < \sigma \}$.
Note that $\sigma$ as in (A)($\alpha$) is uncountable because we have 
the condition $(^\om V_1)^{V_2} \subseteq V_1$.
For (3), take all names for c.c.c\ forcing notions, not only the
for the $\sigma$-linked ones. The additional premise ensures that
(the new) $P$ has the c.c.c.\ in $V_2$ as well. \proofendof{\ref{1.2}}

We get the following conclusion for cardinal characteristics in (B) of
\ref{1.2}:

\begin{theorem} \label{1.3}
In (B) of \ref{1.2} we have 

\medskip

a) $\gb^{V_1} \geq \lambda$, $\gb^{V_2} \leq \sigma$ (and in 
the construction from the proof of \ref{1.1}, we have $\gb^{V_1} = \lambda$.
Moreover, if  
$\forall B \in ([[\mu]^{<\lambda}]^{<\sigma})^{V_2} \: 
\exists B' \in ([[\mu]^{<\lambda}]^{<\lambda})^{V_1} \: 
B \subseteq B'$, then 
the construction from \ref{1.1} gives $\gb^{V_2} = \sigma$).

\medskip

b) $\unif{\meager}^{V_1} \geq \lambda$, $\unif{\meager}^{V_2} \leq 
\sigma$,

\medskip

c) $\cov{\lebesgue}^{V_1} \geq \lambda$, $\cov{\lebesgue}^{V_2}
\leq \sigma$.
\end{theorem}

\proof The $V_1$-part of a), b), and c): 
$\MA_{<\lambda}(\sigma\mbox{-linked})
$ implies that the three cardinal characteristics (and 
$\add{\meager}$, $\add{\lebesgue}$) 
are $\geq \lambda$, because 
all of them can be increased by $\sigma$-linked notions of forcing 
(see e.g.\ \cite{BJ}). 

 \smallskip

In order to  show $\unif{\meager},\gb \leq \sigma$, we 
take $\{ r_i \such i \in C' \}$, $C'\subset C$, $|C'| = \sigma$.             
This set is unbounded and not meagre in $V_2$, because for
any $s \in V_2$ (either in $\omega^\omega$ or as a name for 
a meagre ($F_\sigma$-)set) there is some $B_s \in [C]^{<\sigma}$ 
such that for $i \in C'\setminus B_s \neq \emptyset$ we have the 
$r_i$ is Cohen  over
$V_2[s]$, hence it is not bounded by $s$ nor in a meagre set coded by $s$.

 \smallskip

Proof of $\cov{\lebesgue} \leq \sigma$:
This follows from Rothberger's inequality
$\cov{\lebesgue} \leq \unif{\meager}$ (see \cite{Rothberger38, Fremlin}).
In order to give a proof not using this inequality, 
 we can take $\{ r_i \such i \in C'\}$ 
as above. We set $M(r_i) = \{ m \such r_i \mbox{ is Cohen over } V[m]\}$.
Then (by Fubini) we have that $M(r_i)$ is a Lebesgue null set and for 
$s \in (2^\omega)^{V_2}$ we have there is some 
$B_s \in [C']^{<\sigma}$ such that for $i \in C' \setminus B_s$ 
the real $r_i$ is Cohen over $V[s]$, hence $s \in M(r_i)$, so
$\{ M(r_i) \such i \in C' \}$ covers $(2^\omega)^{V_2}$.

\smallskip

Regarding the part of a) in parentheses: Any $\lambda$ of the Cohen
reals added in the beginning are unbounded and show that $\gb^{V_1}
\leq \lambda$. Under the additional premises, we have that
$\gb^{V_2} \geq \sigma$: Suppose that $M \subset (^\omega 2)^{V_2}$ 
and
$|M|^{V_2} < \sigma$. We take $M_1 \subseteq \mu$ and 
$M_2 \subseteq \mu^{<\lambda}$ such that each member of 
$M$ has a name containing 
only conditions from $\{C_i \such i \in M_1 \} \cup 
\{ P^*_{\mu \pluss a_\beta} \such
\beta \in M_2 \}$. Then $B=
\{ \{ i\} \such i \in  M_1 \} \cup \{ a_\beta \such \beta \in M_2 \} \in
([\mu^{<\lambda}]^{<\sigma})^{V_2}$. Hence there is some
$B' \in ([\mu^{<\lambda}]^{<\lambda})^{V_1}$ such that $B' \supseteq B$.
We take $\beta$ such that $a_\beta \supseteq \bigcup B'$.
Hence at some later stage Hechler forcing
over $V^{P^*_{\mu \pluss a_\beta}}$ will be done in the iteration
and add a real that dominates all reals in $M$. 

\smallskip

\noindent{\bf Remark on the violation of covering.} 
Assume that for some first order sentence $\phi= \phi(P,\in)$, where
$\in$ is a two place predicate and $P$ is a unary predicate, we have that
\begin{eqnarray*}
&\vdash & \forall x Px \rightarrow \phi,\\
& \phi & \mbox{is preserved by increasing } P.
\end{eqnarray*}
Then we define $$
{\rm inv}^\phi = \min \{ |A| \such ({H}(\aleph_1), \in, A) 
\models \phi \}.$$
${H}(\mu)$ is the set of all sets that are
hereditarily of cardinality less that $\mu$. Now, if we
have two models $V_1$, $V_2$ of set theory such that
\begin{equation*}
\begin{split}
& V_1 \subseteq V_2, \mbox{ and }\\ 
&V_1 \mbox{ and } V_2 \mbox{ have the same cardinals and the same 
${H}(\aleph_1)$}\\
&\mbox{\qquad \qquad (which is the same as having the same reals), and } \\
&C \mbox{ is of minimal cardinality such that }({H}(\aleph_1), \in, C) 
\models \phi \mbox{ and }\\
&({\rm inv}^{\phi})^{V_1} = \lambda > |C| \geq ({\rm inv}^{\phi})^{V_2},
\end{split}
\end{equation*}
then we have that $C$ is not
covered by any set in $V_1$ of cardinality less than $\lambda$.

\bigskip

\noindent{\bf Remark on changing $\gd$, $\cof{\meager}$ and 
$\cof{\lebesgue}$.} 
Assume that for some first order sentence $\phi= \phi(\in)$, where
$\in$ is a two place predicate, we have that
\begin{gather*}
\forall x y z \in H(\aleph_1) \;\;
(\phi(x,y) \wedge \phi(y,z) \rightarrow \phi(x,z)) \wedge
\\
\forall x \in H(\aleph_1) \;\;\exists y \in H(\aleph_1) \;\;
\phi(x,y)
\end{gather*}
Then we define for $B \subseteq H(\aleph_1)$, $B \in V$:
$$
{\rm inv}^V_{\phi,B} = \min \{ |A| \such \mbox{for all } x \in
B \; \mbox{ exists } y \in A \mbox{ such that } ({H}(\aleph_1), \in) 
\models \phi(x,y) \}.$$ 
Note that $\gd$, $\cof{\meager}$ and 
$\cof{\lebesgue}$ are characteristics of this type.

Now we have:

\begin{theorem}\label{1.4}
If $V_1$ and $V_2$ are two models of $\zfc$, such that
$V_1 \subseteq V_2$ and 
such that they have the same cofinalities and the same reals, and if
$B \in V_1$, $B \subseteq H(\aleph_1)$, then 
$${\rm inv}^{V_1}_{\phi,B} \leq {\rm inv}^{V_2}_{\phi,B}.$$
\end{theorem}

\begin{corollary}\label{1.5}
If $V_1$ and $V_2$ are two models of $\zfc$, $V_1 \subseteq V_2$ and 
they have the same cofinalities and the same reals 
then their dominating numbers and their cofinalities of the ideals of
Lebesgue null sets and meagre sets coincide.
\end{corollary}

\noindent{\em Proof of Theorem~\ref{1.4}.\/}
Given $V_1$ and $V_2$ and $\phi$ we carry out an induction 
over ${\rm inv}^{V_1}_{\phi,B}$ simultaneously for all $B
\subseteq H(\aleph_1)$, $B \in V_1$.

If  ${\rm inv}^{V_1}_{\phi,B} =1$, then the premise
$H(\aleph_1)^{V_1} = H(\aleph_1)^{V_2}$ and the requirements on $\phi$
immediately yield the claim.

\smallskip
 
Now suppose that the claim is proved for all $\phi$, $B$ such that
${\rm inv}^{V_1}_{\phi,B} < \kappa$ and that we have some
$\phi$, $B$ such that
${\rm inv}^{V_1}_{\phi,B} = \kappa$.

\smallskip

First case: $\kappa$ is regular in $V_1$ and hence in $V_2$.
In this case, Blass' Prop.~2.3 of \cite{Mi2} applies. For completeness' sake we
repeat the argument here:
Suppose that ${\rm inv}^{V_2}_{\phi,B} = \mu \leq \kappa$

Let ${Z} = \{ z_\alpha \such \alpha < \kappa \}$ 
witness ${\rm inv}^{V_1}_{\phi,B} = \kappa$,  and 
${Z}' = \{ z'_\alpha \such \alpha < \mu \}$ 
witness ${\rm inv}^{V_2}_{\phi,B} = \mu$.  
Since ${\mathbb R}^{V_2} \subseteq {\mathbb R}^{V_1}$, 
in $V_2$ there is a function $h \colon \mu \to \kappa$ such that 
for $\alpha < \kappa $:
$$
H(\aleph_1) \models \phi(z'_\alpha, z_{h(\alpha)}).$$
If $\mu $ were less than $\kappa$, then $\rge(h)$ would be bounded in 
$\kappa$, 
say by a bound $\beta \in \kappa$.

Then $\forall a \in {\mathbb R}^{V_1} \; \exists \alpha \in \mu 
\;\phi( a,  z'_\alpha)\wedge
\phi(z'_\alpha, z_{h(\alpha)})$.
 Hence
$\{z_\alpha \such \alpha \leq \beta \}$ 
were a witness for $ {\rm inv}^{V_1}_{\phi,B} \leq \card(\beta) < \kappa$,
which contradicts the premise.

\smallskip

Second case: $\kappa$ is singular in $V_1$ and hence in $V_2$.

Let $\kappa = \lim_{i \to \cf(\kappa)} \kappa_i$ and $\kappa_i < \kappa$.

Let ${Z} = \{ z_\alpha \such \alpha < \kappa \}$ 
witness ${\rm inv}^{V_1}_{\phi,B} = \kappa$.

Set 
\begin{eqnarray*}
Z_i & = & \{ z_\alpha \such \alpha \in \kappa_i\}\mbox{ and }\\
B_i & = & \{ b \in B \such \exists z \in Z_i \; \phi(b,z)\} 
\end{eqnarray*}

Now we have that
\begin{eqnarray*} 
 {\rm inv}^{V_1}_{\phi,B_i}& \leq &\kappa_i, \mbox{ and }
\\
\sup_{i \in \cf(\kappa)}
{\rm inv}^{V_1}_{\phi,B_i}& = &\kappa.
\end{eqnarray*}
The second equation is easy to see:
If $\sup_{i \in \cf(\kappa)}
{\rm inv}^{V_1}_{\phi,B_i} = \theta < \kappa$ then we would have that
 ${\rm inv}^{V_1}_{\phi,B} = \theta \cdot \cf(\kappa) < \kappa$.

\smallskip

By induction hypothesis
$$ {\rm inv}^{V_1}_{\phi,B_i} \leq {\rm inv}^{V_2}_{\phi,B_i}.$$

Since any witness for the computation of ${\rm inv}^{V_2}_{\phi,B}$
is a union of witnesses of the computation of
${\rm inv}^{V_2}_{\phi,B_i}$, we get that
${\rm inv}^{V_2}_{\phi,B} \geq \sup\{{\rm inv}^{V_2}_{\phi,B_i}
\such i \in \cf(\kappa)\}
= \kappa$
\proofendof{\ref{1.4}}

%
\section{Changing the Uniformity of Lebesgue Measure}\label{S2}

In this and the next three sections, we show 
how to change $\unif{\lebesgue}$  (and $\cov{\meager}$, which
comes for free, because of the inequality
 $\cov{\meager} \leq \unif{\lebesgue}$, see~\cite{Fremlin}) 
under our given side conditions. In this section we start to define the 
forcings we are going to use and look at automorphisms of forcings.
We carry out the proof of the changing procedure up to some point 
in the proof of item $\eps)$ of our main Theorem~\ref{2.1} at which 
techniques about transferring information about $\om$-tuples of conditions
(in \cite{Shelah592} called ``whispering'') are needed. We try
 to give some motivation for this fact by proving a lemma about
a pure Cohen situation (Lemma~\ref{2.12}),  of which a weakened
 analogue for iterations of partial random reals and
small c.c.c.\ forcings will be used later.  This weakened analogue is 
the statement  $(**)_{\bar{Q}}$ introduced in \ref{2.11} and proved 
only by the end of Section~\ref{S5}.

\smallskip
   
These technical parts are then carried through in Sections~\ref{S3}, \ref{S4}
and \ref{S5}.

\begin{theorem}\label{2.1}
Assume that we have 
\begin{myrules1}
\item[a)] $V_1 \subseteq V_2$, both models of $\zfc$, $(^\omega V_1)^{V_2}
\subseteq V_1$ [and ($\beta$) or (($\gamma$) + ($\beta$)) from \ref{1.2}(A)], 

\item[b)] $C \in V_2$, $|C| < \lambda$, $C\subseteq \mu$, $\lambda \leq
\mu$,

\item[c)] $\forall B \in V_1$, if $V_1 \models |B| < \lambda$, 
then ${\rm sup}(C \setminus B) = \mu)$,

\item[d)] $\cf^{V_1}(\mu)  > \aleph_0$ and $\cf^{V_1}(\lambda)  > \aleph_0$.

\item[e)] In $V_1$, there are uncountable cardinals $\chi \geq 2^\mu$ and
$\kappa $ such that 
$\kappa < \chi$ and $2^\kappa \geq \chi$.

\end{myrules1}

Then for some c.c.c.\ $P$ in $V_1$ we have
\begin{myrules1}
\item[$\alpha$)] $V_1 \models P$ is a finite support iteration of
$\sigma$-linked forcing notions, 
\item[$\beta$)] $P$ is c.c.c.\ in $V_2$, and
 
\end{myrules1}
for $G 
\subset P$ generic over $V_2$ we have
\begin{myrules1}
\item[$\gamma$)]
 $(^\om V_1[G])^{V_2[G]} \subseteq V_1[G]$,
[and ($\beta$) or (($\gamma$) + ($\beta$)) from \ref{1.2}(A)], 

\item[$\delta$)] $\unif{\lebesgue}^{V_2[G]} \leq |C|^{V_2[G]}$

\item[$\varepsilon$)]
$\unif{\lebesgue}^{V_1[G]} \geq \lambda$.
\end{myrules1}
\end{theorem}
\proof We work in $V_1$ (and often write $V$ 
instead of $V_1$). For $\chi \geq 2^\mu$ we let $g_\chi \colon
\chi \to [\mu]^{<\lambda}$ 
increasing with $\chi$, that is for $\chi \leq \chi'$ we have that $g_{\chi'}
\restriction \chi = g_\chi$,  \marginparr{$g_\chi$}
and
$$ \forall B \in [\mu]^{<\lambda} \: \: 
\exists^\chi \alpha < \chi \: \: g_\chi(\alpha) = B.$$ 
For $\xi < \mu$ let \marginparr{$E^\chi_\xi$}\marginparr{$A_{\chi + \xi}^\chi$}
\begin{eqnarray*}
E_\xi & = & E_\xi^\chi = 
\{ \alpha < \chi \such \xi \not\in g_\chi(\alpha) \}, \mbox{ and}\\
A_{\chi + \xi}^\chi & =& E^\chi_\xi \cup [\chi, \chi+ \xi).
\end{eqnarray*}

\bigskip

We take $\mu$ and $\lambda$ as in the premises of \ref{2.1}. 
We also fix 
$\kappa \geq \aleph_1$ and some $\chi \geq 2^\mu$ as above
such that $\cf(\chi) > \mu$
(used in \ref{2.11} on page~\pageref{pigeonhole})
and $2^\kappa \geq \chi$ and such that $\kappa < \chi$
(for our special iteration where all $Q_\alpha$ of cardinality
$<\kappa$ are already countable, $\kappa \leq \chi$ would suffice, see at
\ref{5.2} and the remarks in \ref{2.11}, if you like to work with
weaker premises).
Note for use in \ref{5.5}:
The definition of $g_\chi$ and $E_\xi$, $A^\chi_{\chi + \xi}$ makes sense 
also 
if $2^\kappa < \chi$.
 
\marginparr{$\chi$,$\kappa$}

\smallskip

\begin{definition}\label{2.2}
1) ${\mathcal K}$ is the class of sequences
$$\bar{Q} = \langle P_\alpha,\Name{Q}_\beta, A_\beta, \mu_\beta,
\name{\tau}_\beta \such \alpha\leq \alpha^*, \beta < \alpha^* \rangle$$
satisfying:

\begin{myrules}

\item[(A)] $\langle P_\alpha, \Name{Q}_\beta \such \alpha
\leq \alpha^*, \beta < \alpha^* \rangle$ is a finite support iteration 
of c.c.c.\ forcings. We call $\alpha^* =\lgg(\bar{Q})$ the 
length of $\bar{Q}$, and
$P_{\alpha^*}$ is the limit. 

\item[(B)] $\name{\tau}_\alpha \subseteq \mu_\alpha < \kappa$ is a name of
the  generic
of $\Name{Q}_\alpha$, 
  i.e.\ over $V^{P_\alpha}$ from $G_{\Name{Q}_\alpha}$  we can compute
$\name{\tau}_\alpha$ and vice versa.

\item[(C)] 
$A_\alpha \subseteq \alpha$.

\item[(D)] $\Name{Q}_\alpha$ is a $P_\alpha$-name of a c.c.c.\ forcing 
notion that is computable from $\langle 
\name{\tau}_\gamma[\Name{G}_{P_\alpha}]$
$\such \gamma \in A_\alpha \rangle$.

\item[(E)] $\alpha^* \geq \chi$ and for $\alpha < \chi$ we have that 
$Q_\alpha = (^{\om} 2, \vartriangleleft)$ (the Cohen forcing) and 
$\mu_\alpha = \aleph_0$ (identify $^{<\om}2 $ with $\om$).

\item[(F)] For each $\alpha< \alpha^*$ one of the following holds 
(and the case is determined in $V$):

\begin{myrules}
\item[($\alpha$)] 
$|\Name{Q}_\alpha| < \kappa$, $|A_\alpha| < \kappa$ and 
(just for notational simplicity) 
the set of elements of $Q_\alpha =
\Name{Q}_\alpha[G_{P_\alpha}]$ is $\mu_\alpha < \kappa$ (but the order 
not necessarily the order of the ordinals) and $Q_\alpha$ is separative
(i.e. $\alpha \Vdash \beta \in \Name{G}_{Q_\alpha} \Leftrightarrow
 Q_\alpha \models \beta \leq \alpha$)

\item[($\beta$)] $Q_\alpha = {\rm Random}^{V[\name{\tau}_\gamma[G_{P_\alpha}]
\such \gamma \in A_\alpha]}$ and $|A_\alpha| \geq \kappa$.
\end{myrules}
\end{myrules}

\medskip

2)For the proof of \ref{2.1} we shall be using the following instance of 
1): For $\chi$, $\mu$, $A^{\chi}_\alpha$ as above
 we define a finite support iteration
$$
\bar{Q}^\chi = \langle P_\alpha^\chi, \Name{Q}_\beta^\chi, A^\chi_\beta,
\aleph_0, \name{\tau}_\beta \such \alpha \leq \chi + \mu,
\beta < \chi + \mu  \rangle,$$
 $P^\chi = P^\chi_{\chi+\mu}$.
For $\alpha < \chi$ we let $Q_\alpha^\chi = (^{<\omega} 2, \vartriangleleft)$, 
the Cohen forcing. For $\alpha = \chi + \xi$, $\xi < \mu$, we let
$$
\Name{Q}^\chi_\alpha = {\rm Random}^{V[\name{\tau}_\beta^\chi 
\such \beta \in A_\alpha^\chi]},
$$
where $\name{\tau}_\beta^\chi$ is 
$\Name{Q}^\chi_\beta$-generic over $V^{P_\beta}$.
\end{definition}

Thus, the $\bar{Q}^\chi$ from b) is a member of ${\mathcal K}$
(and of 
\cite[Def.\ 2.2]{Shelah592} and \cite[Definition 1.4]{Shelah619}) of a 
special form:
$A_\alpha= \emptyset$ if $\alpha < \chi$, and $A^\chi_{\chi +\xi}=
E_\xi \cup [\chi, \chi+\xi)$ for $\xi< \mu$.

\smallskip

The reader may wonder why we do not really fix $\chi$. The reason 
is that in Section~\ref{S5}
we use a L\"owenheim Skolem argument and work simultaneously with
$\chi$, $\chi^+$, $\chi^{++}$, \dots, $\chi^{+(n-1)}$, $n$ the size of some 
heart of a $\Delta$-system, in order to expand
$\bar{Q}^\chi$ to a richer structure that will be used for the proof of 
part $\varepsilon$) of \ref{2.1}.

\smallskip

The Lebesgue measure is denoted by 
${\rm Leb}$ and for a tree $T \subseteq 2^{<\om}$ 
we define $\lim(T) = \{ f \in 2^\om \such \forall n \in \om \;
f \restriction n \in T \}$.
Similar to \cite[2.2]{Shelah592}, we specify dense suborders of
Random and call them Random again:

\begin{definition}\label{2.3} 
a)\\
 ${\rm Random}^{V[\name{r}_\alpha \such \alpha \in A]} =  \{ p \such$
\parbox[t]{0.67\textwidth}
{there is in $V$ a Borel function ${\mathcal B}^p = {\mathcal B}$ with variables
ranging among $\{\mbox{true}, \mbox{false}\}$ and 
range perfect subtrees $r$ of
$^{< \omega} 2$ with ${\rm Leb}(\lim(r)) > 0$ such that
$\forall \eta \in r \: {\rm Leb}(\lim r^{[\eta]} >0)$
(where $r^{[\eta]} = \{ \nu \in r \such
 \nu \trianglelefteq \eta \vee \eta \trianglelefteq \nu \}$) and 
there are pairs 
$(\gamma_\ell,\zeta_\ell)$ for $\ell \in \omega$, such that
 $\gamma_\ell \in A$, 
$\zeta_\ell \in \omega $, and such that $p = {\mathcal B}^p((\mbox{truth 
value}(\zeta_\ell \in \name{r}_{\gamma_\ell}))_{\ell \in \omega}) \}$.}

\smallskip

b) In this case we let $\supt(p) = \{ \gamma_\ell \such \ell \in \omega \}$.

\smallskip

c) $P'_\alpha = \{ p \in P_\alpha \such $\parbox[t]{0.8 \textwidth}
{$\forall \gamma \in \dom(p), \mbox{ if } |A_\gamma| < \kappa, \mbox{ then }
p(\gamma) \in \mu_\gamma$\\
 (not just a name for a member of $\mu_\gamma$),\\
and if $|A_\gamma| \geq \kappa$, then $p(\gamma) \in {\rm Random}^
{V[\name{r}_\delta \such \delta \in A_\gamma]} \}.$}

\smallskip

d) For $A \subseteq \alpha$, we set
$$P'_A =\{ p \in P_\alpha \such \dom(p) \subseteq A \wedge 
\forall \gamma (\gamma \in \dom(p)
\rightarrow \supt(p(\gamma)) \subseteq A)\}.$$

\smallskip

e) $A \subseteq \alpha$ is called $\bar{Q}$-closed or called
$\langle A_\gamma \such \gamma \in \alpha^*\rangle$-closed if
$$\forall \alpha \in A \; (|A_\alpha| < \kappa \rightarrow A_\alpha 
\subseteq A).$$

\end{definition}

So, in our situation of Definition \ref{2.2},
 where all non-empty $A_\alpha$ have size $\chi
\geq \kappa$,  any $A \subseteq \chi + \mu$ is $\langle A_\alpha \such \alpha
< \chi + \mu \rangle$-closed.

\begin{fact}\label{2.4}
Let $\bar{Q}^\chi$ be in $\mathcal K$ from Definition~\ref{2.2}.

1) If $\alpha \leq \chi + \mu$
and $\Name{X}$ is a $P_\alpha$-name of a subset of $\theta <\chi+\mu$ 
then there is a 
set $A \subseteq \alpha$ such that
$|A| \leq \theta$ and 
$\Vdash_{P_\alpha} \mbox{``} \Name{X} \in V[\name{\tau}_\gamma
\such \gamma \in A]$''. Moreover for each $\zeta < \theta$ there is in $V$ a
Borel function ${\mathcal B}_\zeta(x_0,x_1, \dots)$ with domain and range the 
set $\{{\rm true}, {\rm false}\}$ and $\gamma_\ell \in A$, $\zeta_\ell 
< \mu_\ell$ for $\ell \in \om$ such that
$$\Vdash_{P_\alpha} \mbox{``} \zeta \in \Name{X}
\mbox{ iff } {\rm true} = {\mathcal B}_\zeta((\mbox{truth value}(\zeta_\ell 
\in \Name{\tau}_{\gamma_\ell}[G_{Q_{\gamma_\ell}}]))_{\ell \in \om})\mbox{''}.
$$

\smallskip
2) For $\bar{Q} \in {\mathcal K}$ and $A \subseteq \alpha$ every real in 
$V[\tau_\gamma \such \gamma\in A]$ 
has the form 
$$({\mathcal B}_n((\mbox{truth value}(\zeta_\ell 
\in \name{\tau}_{\gamma_\ell}[G_{Q_{\gamma_\ell}}]))_{\ell \in \om} 
))_{n \in \om}.
$$ with ${\mathcal B}_n$ as in 1), and ``true'' interpreted by 1 and ``false'' 
interpreted by 0.
\end{fact}

\proof 1) Let $\Name{X}$ be a name for a subset of $\theta$.
Let $\rho$ be a regular cardinal, and let
the relation  $<^\ast_{\rho}$ be a well-ordering of
$H(\rho)$ such that $x \in y$ implies that $x <^\ast_{\rho} y$.
Take $\rho$ such that $(\bar{Q},\theta, \Name{X}) \in H(\rho)$; let
$M$ be an elementary submodel of ${\mathcal H}(\rho)=(H(\rho), \in, 
<^\ast_{\rho})$
to which $\{ \bar{Q}, \Name{X},\theta \} $ belongs and such that
$\theta \subseteq {H}(\rho)$.  

Thus,  $\Vdash_{P_{\alpha^\ast}} \mbox{``} M[\Name{G}_{P_{\alpha^\ast}}] 
\cap H(\rho)= M\mbox{''}$. 
Since $V^{P_\alpha} = V[\name{\tau}_\beta \such \beta \in
\alpha] $ we have that $M[\Name{G}_{P_{\alpha^\ast}}] =
M[\langle \tau_\beta \such \beta \in \alpha \cap M \rangle]$.
So $X \in M[\langle \tau_\beta \such \beta \in \alpha \cap M \rangle]$,
and we may choose a name for $\Name{X}$ of the form
$\Name{X}=
\{ (\zeta, p)\such \zeta \in \mu, p \in C_\zeta \}$ 
where $C_\zeta$ is a maximal antichain
in  $V[\name{\tau}_\gamma \such \gamma \in \alpha \cap M]$ 
and from that we can 
build a Borel function ${\mathcal B}_\zeta$ in $V$ such that
$$ \Vdash_{P_\alpha} \; \mbox{``}
\zeta \in \Name{X} \Leftrightarrow {\mathcal B}_\zeta(\langle 
\mbox{truth value}(\xi_\ell \in \name{\tau}_{\beta_\ell}) \such \ell \in
\om \rangle) = 1 \mbox{''},
$$
where all the 
$\beta_\ell \in \alpha \cap M$.

\smallskip

Hence we have that
$\Vdash_{P_\alpha} \mbox{``} \Name{X} \in V[\name{\tau}_\gamma
\such \gamma \in M \cap \alpha]$''.

\smallskip

2) is a special case of 1) with $\theta=\om$. We may clue the
${\mathcal B}_n$, $n\in \om$, together to one Borel function is this case,
and write all the arguments into all ${\mathcal B}_n$.
\proofendof{2.4}

\smallskip

We are going to combine the techniques of \cite{Shelah592} and of 
\cite{Shelah619}.
We use automorphisms of $P_{\alpha^\ast}$ that stem from permutations
of $\lgg(\bar{Q}) = \alpha^\ast$.

\begin{definition}\label{2.5}
1) For $\bar{Q} \in {\mathcal K}$ of the special form of
\ref{2.2} Part 2), $\alpha < \alpha^*$, we let
\begin{eqnarray*}
AUT(\bar{Q}\restriction \alpha) & = &
\{ f \colon \alpha \to \alpha \such f \mbox{ is bijective, and },
\\&&
(\forall \beta \in \alpha)(\forall \gamma \in [\chi, \alpha))
\\&&
((\beta < \chi \leftrightarrow f(\beta) < \chi)
\wedge (\beta \in A_\gamma  \leftrightarrow f(\beta) \in A_{f(\gamma)}))\}.
\end{eqnarray*}

2) 
We
let for $f \colon \alpha \to \alpha$ the function $\hat{f} \colon
P'_{\alpha} \to P'_{\alpha}$ be defined by
$p_1 = \hat{f}(p_0)$ if $\dom(p_1) =\{ f(\beta) \such \beta \in \dom(p_0)\}$,
$p_1(f(\beta))= {\mathcal B}^\beta_{p_0}((\mbox{truth value}(f(\zeta_\ell) 
\in \name{\tau}_{f(\gamma_\ell)}))_{\ell \in \omega})$, where
$p_0(\beta)= {\mathcal B}^\beta_{p_0}((\mbox{truth value}(\zeta_\ell \in 
\name{\tau}_{\gamma_\ell}))_{\ell \in \omega})$.
(Here, we write ${\mathcal B} $ for $({\mathcal B}_\zeta)_{\zeta \in \mu}$ when
$Q_\beta = \mu$.)

We can also naturally extend $\hat{f}$ onto the set of all 
$P'_{\alpha}$-names and name this extension $\hat{f}$ as well.
\end{definition}

Now we have for $\bar{Q} \in {\mathcal K}$:

\begin{lemma}\label{2.6} (cf.\ \cite[Fact 1.6. parts 4) and  5)]{Shelah619})

\smallskip

1) For $f \in AUT(\bar{Q} \restriction \alpha)$ we have that
$\hat{f}$ is an automorphism of $P'_\alpha$.

\smallskip

2) Let $\otimes_{(\bar{Q}, A)}$ be the following:
\begin{equation}\tag*{$\otimes_{(\bar{Q}, A)}$}\begin{split}
&\mbox{For every $\alpha \in A \cap 
[\chi, \chi + \mu)$ and for every 
countable}\\
& \mbox{$B \subseteq \alpha$ there is some
$f \in AUT(\bar{Q} \restriction \alpha)$ such that}\\
&f \restriction (A \cap B) = {\rm id},\\
&     f''(B) \subseteq A,\\
&     f''(B \cap A_\alpha) \subseteq A \cap A_\alpha.
\end{split}\end{equation}

If $A$ is $\bar{Q}$-closed and $\otimes_{(\bar{Q}, A)}$, then 
$P'_A \lessdot P'_{\lgg(\bar{Q})}$,  and $\forall q \in P'_{\lgg(\bar{Q})}$ 
we have

\begin{myrules1}
\item[(a)] $q \restriction A \in P'_A$,
\item[(b)] $P'_{\lgg(\bar{Q})} \models q \restriction A \leq q $,
\item[(c)] if $q \restriction A \leq p \in P'_A$, then $q' = p 
\cup q \restriction (\lgg(\bar{Q}) \setminus A)$ belongs to 
$P'_{\lgg(\bar{Q})}$
and is the lub of $p,q$.
\end{myrules1}
\end{lemma}

\proof 1) is easy. 2) is carried out as in \cite{Shelah619}, but since we
promised to write the proofs in a  self-contained style, we write 
down a proof here:

\smallskip

We prove by induction on $\beta \leq \lgg(\bar{Q})$ that for $A' =
A \cap \beta$ and $ q \in P'_{\beta}$, clauses a), b), and c)
hold. 

In successor stages $\beta = \alpha +1$, if $\alpha \not\in A$ or $A_\alpha
=\emptyset$ it is trivial. So assume that $\alpha \in A$ and 
 $A_\alpha \neq \emptyset$.
By induction hypothesis, $P'_{A \cap \alpha} \lessdot
P_\alpha$ and the analogues of a), b) and c) hold for stage $\alpha$.
It is enough to show

\begin{myrules}
\item[($\ast$)] if in $V^{P'_{A\cap \alpha}}$, 
$\mathcal I$ is a maximal antichain in ${\rm
Random}^{V^{P'_{A\cap \alpha \cap A_\alpha}}}$, then in
$V^{P'_\alpha}$ the set $\mathcal I$ is a maximal antichain in 
${\rm Random}^{V^{P'_{A_\alpha}}}$.
\end{myrules}

By the c.c.c.\ this is equivalent to

\begin{myrules}
\item[($\ast$)']
if $\zeta^* < \omega_1$, $\{ p_\zeta \such \zeta  < \zeta^*\} \subseteq 
P'_{A \cap(\alpha +1)}$,
$p  \in P'_{A \cap \alpha}$, and
\begin{equation*}
\begin{split}
p \Vdash_{P'_{A \cap \alpha}} &
\mbox{``} \{ p_\zeta(\alpha) \such \zeta < \zeta^* \mbox{ and }
p_\zeta \restriction \alpha \in G_{P'_{A \cap \alpha}} \}\\ 
&\mbox{ is a predense subset of }
{\rm Random}^{V^{P'_{A \cap \alpha \cap A_\alpha}}}\mbox{''},
\end{split}
\end{equation*}
then 
\begin{equation*}
\begin{split}
p \Vdash_{P'_\alpha} &\mbox{``}
\{ p_\zeta(\alpha) \such \zeta < \zeta^* \mbox{ and }
p_\zeta \restriction \alpha \in G_{P'_\alpha} \} \\
&\mbox{ is a predense subset of }
{\rm Random}^{V^{P'_{A_\alpha}}}\mbox{''}.
\end{split}
\end{equation*}
\end{myrules}

Assume that ($\ast$)' fails. So we can find $q$ such that
\begin{eqnarray*}
p &\leq & q \in P'_\alpha,\\
q &\Vdash_{P'_\alpha} & ``\{ p_\zeta(\alpha) \such \zeta < \zeta^* 
\mbox{ and }
p_\zeta \restriction \alpha \in G_{P'_\alpha} \} 
\\ 
&&\mbox{ is not  a predense subset of Random}^{V^{P'_{A_\alpha}}}\mbox{''}.
\end{eqnarray*}

So for some $G_{P'_{A_\alpha}}$-name $\name{r}$ 
$$q 
\Vdash_{P'_{A_\alpha}}
 ``\name{r} \in {\rm Random}^{V^{P'_{A_\alpha}}}(=\Name{Q}_\alpha)
\mbox{ and is
incompatible with every } p_\zeta(\alpha) \in \Name{Q}_\alpha\mbox{''}.
$$
Possibly increasing $q$ w.l.o.g.\ $\name{r}={\mathcal B}((\mbox{truth value}
(\eta_\gamma \in \name{\tau}_\gamma))_{\gamma \in w})$ 
with a suitable countable
$w \subseteq A_\alpha$. Now we choose
\begin{eqnarray*}
B &=& \dom(q) \cup \bigcup_{\zeta < \zeta^*} \dom(p_\zeta \restriction \alpha)
\cup \bigcup \{\supt(q(\beta)) \such \beta \in \dom(q) \}\\
&& \cup  \bigcup\{ \supt(p_\zeta(\beta)) \such 
\beta \in \dom(p_\zeta \restriction \alpha) \mbox{ and } \zeta < \zeta^* \} \cup w.
\end{eqnarray*}
Since $B$ is a countable subset of $\alpha$ and since we have
$\otimes_{(\bar{Q},A)}$ there is an $f \in AUT(\bar{Q} \restriction \alpha)$
such that
\begin{eqnarray*}
f \restriction (B \cap A) &=& \mbox{ the identity},\\
f^{''}(B) & \subseteq & A,
\\
f^{''}(B \cap A_\alpha) &\subseteq & A \cap A_\alpha.
\end{eqnarray*}
As $\hat{f}$ is a automorphism of $P'_\alpha$ and is the identity
  on $P_{A \cap B}$ we have that
\begin{eqnarray*}
\hat{f}(p)&=& p,\\
\hat{f}(p_\zeta) &=& p_\zeta,\\
p & \leq &  \hat{f}(q) \in P'_{A \cap\alpha},\\
\hat{f}(\name{r}) &=& {\mathcal B}((\mbox{truth value}
(\eta_\gamma \in \name{\tau}_{f(\gamma)}))_{\gamma \in w}),
\\
f{''}(w) &\subseteq & f{''} (B \cap A_\alpha) \subseteq A \cap A_\alpha,\\
\mbox{ hence } & \Vdash_{P'_{\alpha}} & \hat{f}(\name{r}) \in 
{\rm Random}^{V^{P'_{A \cap A_\alpha}}},\\
\hat{f}(q) &\Vdash _{P'_{A\cap \alpha }} & \mbox{`` in $\Name{Q}_\alpha$, 
 $\hat{f}(\name{r})$
and $p_\zeta(\alpha)$ are incompatible for $\zeta < \zeta^*$}.
\end{eqnarray*}
and thus get a 
contradiction to the fact that we started with a maximal antichain.
\proofendof{\ref{2.6}}

\begin{lemma}\label{2.7}
For $A = E_\xi \cup [\chi, \chi+ \xi)$, and for $\bar{Q}$ as in 
Definition~\ref{2.2} Part 2), we have that $\otimes_{(\bar{Q}, A)}$ is true.
\end{lemma}

\proof Let $\alpha \in A$ and $B \subseteq \alpha$ be countable.
W.l.o.g., we treat here the case when $\alpha \geq \chi$.
We have to show that there in an $f$ such that
\begin{myrules}
\item $f \colon \alpha \to \alpha$ bijective,
\item $f\restriction \chi  \colon \chi \to \chi$ bijective,
\item $\forall \beta, \gamma < \alpha \; 
(\beta \in A_\gamma \leftrightarrow f(\beta) \in A_{f(\gamma)})$,
\end{myrules}
(These first three items ensure that $f \in AUT(\bar{Q} 
\restriction \alpha)$, and
next we write the conditions in $\otimes_{(\bar{Q},A)}$:)
\begin{myrules}
\item $f\restriction 
((E_\xi \cap B) \cup ([\chi,\xi) \cap B)) = id$,
\item $f''(B) \subseteq E_\xi \cup [\chi,\alpha)$,
\item $\forall \alpha \in [\chi,\chi + \xi) \;\;
f''(B \cap (E_{\alpha - \chi} \cup [\chi,\alpha))) \subseteq
(E_\xi \cap E_{\alpha -\chi}) \cup [\chi,\alpha)$.
\end{myrules}

Next we require that the $f$ preserves slightly more

\begin{myrules}
\item $f \restriction [\chi,\alpha) = id$ and hence
\item $\forall \beta \in [\chi,\alpha] \: 
f\restriction E_{\beta -\chi} 
\colon E_{\beta -\chi} \to E_{\beta -\chi}$.
\end{myrules}

So, $f$ has to map $(B \setminus E_\xi) \cap E_{\alpha -\chi}$ into
$E_\xi \cap E_{\alpha -\chi}$ and
$((B \setminus E_\xi) \setminus E_{\alpha -\chi})\cap \chi$ into
$E_\xi \setminus E_{\alpha -\chi}$.

\smallskip

For $\gamma \in \chi$, $\alpha'\in \xi + 1$
 we write $tp_{\alpha'}(\gamma) =
\{ \beta \in \alpha' \such \gamma \in E_\beta \} =
\{ \beta \in \alpha' \such g(\gamma)\not\ni \beta\} $.
All subsets $T \subseteq \alpha'$ such that $|\alpha' \setminus T|
<\lambda$ are realised  as the type of $\chi$ elements 
 because for  each
$B \in [\mu]^{<\lambda}$ we have $\chi$ many $\gamma$ such that
$g_\chi(\gamma) = B$. 
Since $\alpha - \chi < \xi$, the relation $E_\xi$ does not
play a r\^ole in $tp_{\alpha +1 -\chi}(\gamma)$ and so 
we have that for all such $\alpha +1 -\chi$-types $T$
\begin{eqnarray*}
 |\{\gamma \such tp_{\alpha +1 - \chi}(\gamma)=  T \}|
&  = & \\
|\{\gamma \such tp_{\alpha +1- \chi}(\gamma) = T
\wedge \gamma \in E_\xi \}|
& = &\\
 |\{\gamma \such tp_{\alpha +1- \chi}(\gamma) = 
T \wedge \gamma \not\in
E_\xi \}|
& = &  \chi.
\end{eqnarray*}

Hence there is a bijection $f'$ of $\chi$ preserving the
$\alpha + 1-\chi$-types and being the identity on 
$(E_\xi \cap B) \cup[\chi, \alpha)$ 
but mapping $(B\cap \chi) \setminus E_\xi$ into $E_\xi$.
Then $f=f' \cup id_{[\chi,\alpha)}$ is
as required.
\proofendof{\ref{2.7}}

Now we return to the conclusion of Theorem~\ref{2.1}:
\begin{myrules1}
\item[($\gamma$)]
If $G \subseteq P$ is generic over $V_2$, then
\begin{myrules1}
\item[-] $V_1[G]$ and $V_2[G]$ have the same reals, indeed
$(^\omega V_1[G])^{V_2[G]} \subseteq V_1[G]$
\item[-] 
$V_1[G]$ and $V_2[G]$ have the same cardinals if $(V_1,V_2)$ have
\item[-] $V_1[G]$ and $V_2[G]$ have 
the same cofinality function  if $(V_1,V_2)$ have.
\end{myrules1}
\end{myrules1}

Since Cohen forcing and random forcing are $\sigma$-linked,
the proof of Theorem~\ref{1.1} applies here as well. \proofendof{(\gamma)}

Next we show

\begin{myrules}

\item[($\delta'$)] $V_2 \models \; \Vdash_{P_{\chi+\mu}}\;
\mbox{``}\{ \tau_{\chi +i} \such i \in C \}$ is not null.''

\end{myrules}

\proof
Let $\Name{N}\in V_2$ be a $P_{\chi + \mu}$-name for a Borel null set.
Since $(^\om V_1)^{V_2} \subseteq V_1$ we may assume that $\Name{N}
\in V_1$. By \ref{2.4}(2),
for some Borel function ${\mathcal B} \in V_1$ for some countable
$X = \{ x_\ell \such \ell \in \omega \}\subseteq \chi, Y
=\{ y_\ell \such \ell \in \omega \} \subseteq \mu$, $\zeta_\ell, 
\ell \in \omega$,
$\zeta'_\ell, \ell \in \omega$, we have that
$$\Name{N}=
{\mathcal B}((
\mbox{truth value}(\zeta_\ell \in 
\name{\tau}_{x_\ell}))_{\ell \in \omega},
(\mbox{truth value}(\zeta'_\ell \in \name{\tau}_{\chi 
+ y_\ell}))_{\ell \in \omega}).$$

Let $i(\ast) < \mu$ be such that $i(\ast)
> \sup(Y)$. (Here we use that
$\cf^{V_1}(\mu) > \aleph_0$.)
Since $\cf^{V_1}(\lambda) > \aleph_0$, we have that
$B:= \bigcup_{\xi \in X} g_\chi(\xi) \in ([\mu]^{<\lambda})^{V_1}$.
Since $\sup(C \setminus B) = \mu$, there is some $i \geq i(\ast)$,
$i \in C \setminus B$. We claim, that $r_{\chi + i}$ is random (in the
sense of $V_1$ and hence also in the sense of $V_2$
as Random and all maximal (countable) antichains of
the random forcing are the same in $V_1$ and in $V_2$) over 
an extension of $V_1$, in which $\Name{N}[G]$ has a name. 
Then the proof will be finished, because then 
$r_{\chi + i} \not\in \Name{N}[G]$ in $V_1[G]$ and also in $V_2[G]$.
By our construction, we have
$$
\tau_{\chi+i} \mbox{ is the } {\rm Random}^{V[\name{\tau}_\alpha 
\such \alpha \in E_i \vee
\chi \leq \alpha < \chi + i ]}\mbox{-generic over } V_1^{P_{\chi+ i}}.
$$
Since $i \in C \setminus B$, we have that $\forall \xi \in X \: 
g_\chi(\xi) \neq i$, hence
$\forall \xi \in X \: \xi \in E_i$, so $X \subseteq E_i$. Moreover $\chi +
Y \subseteq [\chi, \chi + i)$, as $i \geq i(\ast) \geq \sup(Y)$.
Since, by Lemmas~\ref{2.6} and~\ref{2.7}, 
$P_{A_{\chi + i}} \lessdot P_{\lgg(\bar{Q})}$ the name $\Name{N}$ is
evaluated in the right manner in $V_1^{P_{A_{\chi + i}}}$.
Thus the claim is proved. \proofendof{(\delta')}
\begin{myrules}
\item[($\delta$)] $V_2[G] \models \unif{\lebesgue} \leq |C|$.
\end{myrules}
This follows from ($\delta'$).

\medskip
Now comes the part whose proof will be finished only at the end of
Section~\ref{S5}.
\begin{myrules}
\item[($\varepsilon$)] $V_1[G] \models \unif{\lebesgue} \geq \lambda$.
\end{myrules}

\proof Suppose that not. In $V_1$
there is $i(\ast) < \lambda$ and $p \in P_{\chi + \mu}$ such that
$$
p \Vdash_{P_{\chi +\mu}} 
\mbox{``}\name{\eta}_i \in \; ^\omega 2 \mbox{ for } i < i(\ast)
\wedge \{ \name{\eta}_i \such i < i(\ast) \} \mbox{ is not null.''}
$$
A name of  a real in $V_1[G]$ is given  by
$$\name{\eta}_i = {\mathcal B}_i((\mbox{truth value}(\zeta_{i,\ell} \in 
\name{r}_{j_{i,\ell}}))_{\ell \in \omega})$$ for
suitable $\langle \zeta_{i,\ell}, j_{i,\ell} \such \ell \in 
\omega \rangle$,
$\zeta_{i,\ell} \in \omega$, $j_{i,\ell} \in \chi + \mu$.

We set 
\begin{eqnarray*}
X &=& \{ j_{i,\ell} \such i \in i(\ast), 
\ell \in \omega \} \cap \chi,\\
Y &=& \{j_{i,\ell} \such i \in i(\ast), \ell \in \omega \} \cap 
[\chi, \chi+ \mu).
\end{eqnarray*}

We show the main point:

\smallskip

In $V_1[G]$, $(^\omega 2)^{V[\{\name{\tau}_\xi \such \xi 
\in X \cup Y\}]}$ 
is a Lebesgue null set.

\smallskip

Since   $\exists^{\chi} \alpha \; g_\chi(\alpha) = Y - \chi$  
we can fix such an $\alpha \in \chi\setminus X$ \marginparr{$\alpha$} that is 
not in $E_\xi$ for every $\xi \in Y-\chi$.
It is important to note that therefore the premises of \ref{2.8} and
or \ref{2.11} can be fulfilled 
for our any $X,Y$ as above, with a suitable choice of $\alpha$.

\begin{lemma}\label{2.8}

In $V_1^{P_{\alpha^*}}$, the set
$(^\omega 2)^{V_1[\tau_\xi \such \xi \in X \cup Y]}$ has Lebesgue measure $0$,
and a witness for a definition for a measure zero superset
can be found in $V^{P_{\alpha+1}}$ (a forcing name is already in 
$V^{P_{\alpha}}$) for  any
$\alpha \in \chi\setminus X$ \marginparr{$\alpha$} that is 
not in $E_\xi$ for every $\xi \in Y-\chi$.

\end{lemma}

\proof Explanation: 
This proof will be finished only with the proof of Lemma~\ref{2.11},
which will, as we already mentioned, only be finished by the end of
 Section~\ref{S5}. The proof of this lemma
requires reworking of almost the whole
\cite{Shelah592}. The lemma is also
stated in  \cite[1.11 and 1.12]{Shelah619}, where a proof
assuming the knowledge of \cite{Shelah592} is given.

\smallskip
 First we introduce some paradigm null sets
(see also \cite[2.4 and 2.5]{Shelah592}):
 
\begin{definition}
1) Suppose that $\bar{a} = \langle a_\ell \such \ell \in \omega
\rangle$ and $\bar{n} = \langle n_\ell \such \ell \in \omega 
\rangle$ are such that for $\ell \in \omega$
\begin{myrules}
\item[(a)] $a_\ell \subseteq \; ^{n_\ell}2$,
\item[(b)] $n_\ell < n_{\ell +1} < \omega$,
\item[(c)] $\begin{displaystyle}
\frac{|a_\ell|}{2^{n_\ell}} > 1-\frac{1}{10^\ell}.
\end{displaystyle}$
\end{myrules}
Then we set $N[\bar{a}] = \{ \eta \in \; ^\omega  2
\such \exists^\infty \ell \: \forall \nu \in a_\ell \; \nu 
\not\trianglelefteq \eta \}$.

\smallskip

2) For $\bar{a}$ as above and $n \in \omega$, we let
${\rm tree}_n(\bar{a}) = \{ \nu \in \;^{<\omega} 2
\such n_\ell \geq \max( n, \lgg(\nu)) 
\rightarrow \nu \restriction n_\ell \in a_\ell \}$.
%
%
\end{definition}\label{2.9}

Then $N[\bar{a}] = ^\omega\!2 \setminus \bigcup_{n \in \omega}
 \lim {\rm tree}_n(\bar{a})$ and ${\rm Leb}(N[\bar{a}]) =0$.
The definitions $N[\bar{a}]$ and  
$ \lim {\rm tree}_n(\bar{a})$ may be intepreted in any model $V$ such that
$\bar{a} \in V$. We indicate the model of set theory in which we
evaluate them by superscripts.

\begin{definition}\label{2.10}
For $\beta < \chi$ we identify $Q_\beta$, the Cohen forcing, with
\begin{eqnarray*} \{ \langle (a_\ell, n_\ell) \such \ell < k \rangle &\such &
k \in \omega, n_\ell < n_{\ell +1} < \omega,
 a_\ell \subseteq \; ^{n_\ell} 2, \frac{|a_\ell|}{2^{n_\ell}}
> 1 - \frac{1}{10^\ell} \}.
\end{eqnarray*}
If $G_{Q_\beta}$ is $Q_\beta$-generic, let
\begin{eqnarray*}
\bar{a}^\beta = \name{\bar{a}}^\beta[G_{Q_\beta}] 
&=&
\{ (\ell,a) \such \exists k \geq
\ell + 1 \; \exists \langle (a_j, n_j) \such j < k \rangle \in G_{Q_\beta} \\ 
&&
\exists j<k \; (\ell,a) = (j,a_j)\},
\end{eqnarray*}
and define 
$\name{\bar{n}}^\beta[G_{Q_\beta}]$ analogously. We let
$\name{\bar{a}}^\beta = \langle \name{a}^\beta_\ell \such \ell \in \omega 
\rangle$ and $
\name{\bar{n}}^\beta = \langle \name{n}^\beta_\ell \such \ell \in \omega 
\rangle$ be the names for the corresponding objects.
\end{definition}

\begin{lemma}\label{2.11}
If $\betta \in \chi \setminus X$ is such that
$\forall \xi \in Y-\chi \: \betta \not\in E_\xi$,
then 
\begin{equation*}
(^\omega \! 2)^{V[r_\xi \such \xi \in X \cup Y]} 
\subseteq (N[\bar{a}^{\betta}])^{V[G]}.
\end{equation*}
\end{lemma}

\noindent
{\em Beginning of the proof.}
In this section, we shall only show that 

\begin{equation}\tag*{$(**)_{\bar{Q}}$}\label{star}
\begin{split}
& \mbox{ in } V[G], \mbox{ for } E \in [\chi]^{\kappa^+} 
\mbox{ we  have}\\
&\bigcap_{\betta \in E} {\rm tree}_{\ell^\ast}
(\bar{a}^{\betta}) \mbox{ does not  contain a perfect tree}.
\end{split}\end{equation}

is a sufficient condition for 
\ref{2.11}. For 
certain members $\bar{Q}$ of $\mathcal K$,
\ref{star} will be proved in the next three 
sections. 
Let $\betta \in \chi \setminus X$ be such that
$\forall \xi \in Y-\chi \: \betta \not\in E_\xi$.

We show by induction on $\gamma \geq \chi$ that
\begin{equation}\tag*{$(**)_{\bar{Q}\restriction \gamma}$}\label{stargamma}
\begin{split}
&\mbox{ in $V^{P_\gamma}$, for } E \in [\chi]^{\kappa^+} 
\mbox{ we  have}\\
&\bigcap_{\betta \in E} {\rm tree}_{\ell^\ast}
(\bar{a}^{\betta})
 \mbox{ does not  contain a perfect tree}.
\end{split}\end{equation}
implies: 
\begin{equation}\tag*{$(\ast)_{\bar{Q}\restriction \gamma}$}
\begin{split}
& \forall X \subseteq \chi \;\; \forall Y \subseteq [\chi,\chi+\mu)
\;\; \\
& \forall \betta \in \chi \setminus X (
\;\; \forall \xi \in Y-\chi \: \betta \not\in E_\xi \longrightarrow \\
&
(^\omega \! 2)^{V[r_\xi \such \xi \in (X \cup Y)\cap \gamma]} 
\subseteq (N[\bar{a}^{\betta}])^{V^{P_\gamma}}).
\end{split}
\end{equation}
\smallskip

Preliminary remarks:
Assuming $\neg (\ast)_{\bar{Q}\restriction\gamma}$ we get a 
$P_\gamma$-name $\name{b}$ referring only to
$r_\xi$, $\xi \in (X \cup Y)\cap \gamma$ such that
$$p \Vdash_{P_\gamma}  \name{b} \not\in
N[\name{\bar{a}}^\betta].$$

\smallskip
Since  $\forall \xi \in Y-\chi$ $\betta \not\in E_\xi$,  we have for all 
$\xi' = \chi + \xi \in Y $, $\betta \not \in E_{\xi}\cup[\chi,
\chi + \mu) = A^\chi_{\xi'}$.
Since all $r_{\xi'}$, $\xi' \in Y$ are ${\rm Random}^{V^{P_{A_{\xi'}}}}$-generic
there are automorphisms  $f_\zeta \in AUT(\bar{Q})$, $\zeta \in \chi$, leaving 
$\name{b}$ and every point from $[\chi,\chi+\mu)$ fixed  and moving
$\betta$ to $\betta_\zeta \not\in \{ \betta_{\zeta'} \such \zeta' < \zeta\}$.
 Hence  we get
  $$p_\zeta = \hat{f_\zeta}(p) \, \Vdash_{P_\gamma} \, \name{b} \not\in
\bigcup_{\zeta \in \chi}N[\name{\bar{a}}^{\betta_\zeta}]$$
for $\chi \geq \kappa^+$ pairwise different $\betta_\zeta$'s.
\smallskip

Now we start the induction.\\
For $\gamma = \chi$ the proof is easy, because
$( ^\omega 2)^{V[r_\xi \such \xi \in (X \cup Y) \cap \chi ]}$ contains
only Cohen reals: If there is one real $\name{b}[G_\gamma]$  not in
$(\bigcup_{\zeta \in \kappa^+}N[\bar{a}^{\betta_\zeta}])^{V^{P_\gamma}}$, 
then  this real is Cohen and  gives rise to a perfect tree of Cohen reals
not in 
$(\bigcup_{\zeta \in \kappa^+}N[\bar{a}^{\betta_\zeta}])^{V^{P_\gamma}}$.
So we have that $\neg (\ast)_{\bar{Q}\restriction \gamma}$
 implies $\neg$ \ref{stargamma}. 

\smallskip

Now let $\gamma > \chi$ be a limit.
Assuming $\neg (\ast)_{\bar{Q}\restriction\gamma}$ we get a 
$P_\gamma$-name $\name{b}$ referring only to
$r_\xi$, $\xi \in (X \cup Y)\cap \gamma$ such that
$$p \Vdash_{P_\gamma}  \name{b} \not\in
N[\name{\bar{a}}^\betta].$$
By automorphisms leaving $\name{b}$ and moving
$\betta$ to $\betta_\zeta$ and $p$ to $p_\zeta$ we get
  $$p_\zeta   \Vdash_{P_\gamma}  \name{b} \not\in
\bigcup_{\zeta \in \chi}N[\name{\bar{a}}^{\betta_\zeta}]$$
for $\chi$ pairwise different $\betta_\zeta$'s.

\smallskip

Because of the induction hypothesis we may assume that $p \Vdash_{P_\gamma}
\name{b} \not\in V^{P_\delta}$ for $\delta < \gamma$, and hence
by the properties of c.c.c.\ iterations that $\cf(\gamma) = \aleph_0$.

\smallskip

So for each $\zeta < \chi$ there are $p_\zeta$, $m_\zeta$ such that
$$p \leq p_\zeta  \in P_\gamma, \:\;
p_\zeta \Vdash \name{b} \in 
\lim {\rm tree}_{m_\zeta}(\name{\bar{a}}^{\betta_\zeta}).
$$ 
By properties of c.c.c.\ forcing notions 
$\langle \{ \zeta <\chi \such p_\zeta \in P_\delta \} \such \delta \in \gamma
\rangle $ is an increasing sequence of subsets
of $\chi$ of length $\gamma \leq \mu$. In the beginning
on the proof of \ref{2.1} we chose 
$\mu<\chi$\label{pigeonhole}.
So for some $\gamma_1 < \gamma$
there is $E \in [\chi]^{\kappa^+}$ such that
$p_\zeta \in P_{\gamma_1}$ for $\zeta \in E$ and
$m_\zeta = m$ for $\zeta \in E$. Note that for all but $< \kappa^+$
of the ordinals $\eta \in E$ we have that
$$p_\eta \Vdash |\{ \zeta \in E \such p_\zeta \in G_{P_{\gamma_1}}
\}| = \kappa^+.
$$ Fix such an $\eta$, and let $G_{P_{\gamma_1}}$ be $P_{\gamma_1}$-generic
over $V$ so that $p_\eta \in G_{P_{\gamma_1}}$. 
In $V[G_{P_{\gamma_1}}]$, let $E' =
\{ \zeta \in E \such p_\zeta \in G_{P_{\gamma_1}} \}$, so $|E'| = 
\kappa^+$. 
Let $T^* = \bigcap_{\zeta \in E'} {\rm tree}_m(\bar{a}^{\betta_\zeta})$.
In $V^{P_{\gamma}}$, $T^*$ is a subtree of $^{<\omega}2$ and by
$(**)_{\bar{Q}\restriction \gamma}$, $T^*$ contains no perfect subtree. Hence
$\lim(T^*)$ is countable, so absolute:
$T^*$ is a $P_{\gamma_1}$-name and
$(\lim(T^*))^{V[G_{P_\gamma}]} = (\lim(T^*))^{V[G_{P_{\gamma_1}}]}$. But 
$p_\eta \Vdash \name{b} \in \lim(T^*)$, hence
$p_\eta \Vdash \name{b} \in V^{P_{\gamma_1}}$, a contradiction.

\smallskip 
Assume now that $\gamma = \delta +1$ and 
that  $\neg (\ast)_{\bar{Q}\restriction\gamma}$.
 Choose
$p_\zeta = p'_\zeta * \name{q}_{\delta}(\zeta)$ as in the preliminary
remark such that 
$p_\zeta \in P_\delta$, $\name{q}_{\delta}(\zeta)
\in Q_\delta$, 
 and additionally such that the $\name{q}_{\delta}(\zeta)$ all coincide
(because we may assume that
$f_\zeta$, chosen as in the preliminary remarks, does not move $\delta$), say that all
 $\name{q}_{\delta}(\zeta) = \name{q}_\delta$.
Choose $E$, $p_\eta$, $G_{P_\gamma}$ analogous to the  above.
We have $E' =
\{ \zeta \in E \such p'_\zeta \in G_{P_{\delta}} \}
=
\{ \zeta \in E \such p'_\zeta * \name{q}_\gamma \in G_{P_{\delta}} \}$,
and similarly to the above,
together with $(\ast\ast)_{\bar{Q}\restriction\gamma}$
we get the contradiction $p_\eta \Vdash \name{b} 
\in V^{P_\delta}$.

\smallskip
Since we have covered the cases $\gamma = \chi$ and $\gamma    > \chi$ 
limit and $\gamma >\chi$ successor, we have finished the proof that
\ref{star} implies the statement in Lemma~\ref{2.11}. 

\medskip
 
Our proof of \ref{star} will in some parts be
 similar to \cite{Shelah592}.
However, the difference to \cite{Shelah592} is that
the our $A^\chi_\alpha, \alpha \in [\chi,\chi+\mu)$ 
(from \ref{2.2} Part2)) are large in cardinality, 
namely the same as the iteration length, and hence
some techniques of \cite{Shelah592} are not applicable here.
We also take the
technique of automorphisms of $\bar{Q}$ taken from \cite{Shelah619}, 
and additionally, like there as well, we are going to work 
$\bar{Q}^\chi$ for many $\chi$'s at the same time.
Tomek Bartoszy\'nski \cite{Bahandbook} gives a simplified exposition of
some of the results of \cite{Shelah592}, 
that the reader might want to consult first.

The proof of \ref{2.11} will be finished only at the end of Section~\ref{S5}.


\smallskip

In the next lemma, which stems from Winfried Just, we show $(**)_{\bar{Q}}$
in the special case that all the  $p_\zeta$ are  Cohen.
It   serves as a motivation for
the rest of our work: it shows that the 
main point is to get something similar to the premise no.\ 3 
of Just's lemma for the partial random conditions. We may (and later  do)
weaken the conclusion of Just's lemma: Instead of requiring the 
intersection to be empty we derive only
that the intersection does not contain a perfect tree, that is
$(**)_{\bar{Q}}$.

\begin{lemma}\label{2.12}
[Winfried Just \cite{Just97}]
Suppose that $\{ p_\zeta \such \zeta \in Z \}$ is a set of
conditions in $P_{\chi + \mu}$ such that
\begin{myrules}
\item[1.] $Z$ is infinite.
\item[2.] $\{ \dom(p_\zeta) \such \zeta \in Z \}$ forms a 
$\Delta$-system with root $u$.
\item[3.] $\exists q \; \forall \zeta \in Z \; p_\zeta \restriction u = q$.
\item[4.] $\betta_\zeta \in \dom(p_\zeta) \setminus u$ for all $\zeta$,
$p_\zeta(\betta_\zeta)$ is  Cohen.
\item[5.] $\exists k^*, n^*$ such that $\forall \zeta \in Z$, if $p_\zeta 
(\betta_\zeta) = \langle(n_\ell^\zeta, a_\ell^\zeta) \such \ell \in
k_\zeta \rangle$ then $k_\zeta = k^*$ and $n^\zeta_{k_{\zeta-1}} = n^*$.
\end{myrules}
We set $\Name{E} = \{ 
\zeta \in Z \such p_\zeta \in G \}$. Then we have for every $\ell^\ast
\in \omega$ that
$$q  \Vdash  \bigcap_{\zeta \in \Name{E}} \lim {\rm tree}_{\ell^\ast}
(\name{\bar{a}}^{\betta_\zeta}) = \emptyset.$$
\end{lemma}

\proof Suppose that not. Then there exist some $\ell^\ast$ and
 some $q_1 \geq q$ and some
name $\name{b}$ for an infinite branch such that
$$q_1  \Vdash  \name{b} \in  
\bigcap_{\zeta \in \Name{E}} \lim {\rm tree}_{\ell^\ast}
(\name{\bar{a}}^{\betta_\zeta}).$$

Let $n> \max\{k^*-1, n^*\}$ and such that
$2^{-n} < 10^{-k^*}$. There are some $r \geq q_1$ and some $\nu$
 such that
$$r \Vdash \name{b} \restriction n = \nu.$$

Now take some  $\zeta$ such that $\dom(p_\zeta) \cap \dom(r) =u$. Since $Z$ is
infinite and all conditions are bounded in size by $k^*, n^*$, such a
$\zeta$ exists.
Finally we set $n^\zeta_{k^\ast} = n$ and $a^\zeta_n = 2^n \setminus \{ \nu \}$ and
$$p_\zeta^+ = p_\zeta \restriction (\dom(p_\zeta) \setminus \{\betta_\zeta\})
\cup \{ (\betta_\zeta,\langle n_\ell^\zeta,a_\ell^\zeta \such \ell 
\leq n \rangle ) \}.$$
Since $\nu \not\in a^\zeta_n$, we get
$$p^+_\zeta \Vdash \name{b} \in \lim {\rm tree}_{\ell^*}
(\name{\bar{a}}^{\betta_\zeta}) \rightarrow \name{b} \restriction n \neq \nu.$$

However, $p_\zeta^+$ and $r$ are compatible. Contradiction.
\proofendof{\ref{2.12}}


\section{About 
Finitely Additive Measures}\label{S3}

In order to prove the existence of a condition $p^\otimes$ that forces 
that many of the $p_\ell$'s (where the
$p_\ell$, $\ell \in \om$ are the first
$\om$ of some thinned out part of the $p_\zeta$ from \ref{2.11}) are in
$G_{\alpha^*}$ we use names $(\Name{\Xi}_\alpha^t)_{t\in {\mathcal T},\alpha \in 
\chi + \mu}$  for 
finitely additive measures. We shall have that for every $\alpha < \chi+\mu$,
$\Vdash_{P_\alpha} \mbox{``}\Name{\Xi}_\alpha^t$ is a  finitely additive
measure on ${\mathcal P}(\om)$''. The superscript $t$ ranges
over some set of blueprints (see \ref{4.1}) and indicates the type of the 
$\om$  conditions $p_\ell$  that are taken care of by $\Name{\Xi}_\alpha^t$,
and there are some coherence requirements regarding different $\alpha$'s.
The $\Name{\Xi}_\alpha^t$ are an  item in the class of
forcing iterations ${\mathcal K}^3$ that we are going to define in \ref{4.2}.
Certain members of $\mathcal K$ can be expanded to members of ${\mathcal K}^3$,
and these expandible members of ${\mathcal K}$ are the notions of forcing for 
which we show $(**)_{\bar{Q}}$ is Sections \ref{S4} and 
\ref{S5}. 

\smallskip

For the expansion of a $\bar{Q}$ in $\mathcal K$ to a member of ${\mathcal K}^3$
some requirements linking the $A_\alpha$ and the $\Name{\Xi}^t_\alpha$ 
need to be 
fulfilled (called ``whispering'' in \cite[Def.~2.11 (i)]{Shelah592}). 
By increasing the $A_\alpha$ these can be satisfied. Another way
is to use the requirements only at finitely many points that are determined 
at a later stage in a proof. We shall work according this latter method:
In our case, where we have also automorphisms as in
\ref{2.4}, we shall first specify
som $\langle p_\ell \such \ell \in \omega \rangle$, and only
thereafter we shall define sufficiently many
$\Name{\Xi}_\alpha^t$ (see \ref{5.5}).

\smallskip

Anyway, the ``sufficiently many $\Name{\Xi}_\alpha^t$'' 
need the same lemmas about extensions
of finitely additive 
measures to longer iterations that are also used to proof  that  
our class ${\mathcal K}^3$ of forcings has enough 
members. These will be Lemmas~\ref{4.5}, \ref{4.6}, and \ref{4.7}. 

\smallskip

This short section collects  some facts about finitely additive measures, 
that can be presented separately before we return to the iterated 
forcings in ${\mathcal K}$ and come to the mentioned lemmas.
All statements of this section, however only few of their proofs, 
can also be found in \cite{Shelah592}.

\begin{definition}\label{3.1}
1) ${\mathcal M}$ is the set of functions $\Xi$ from some Boolean subalgebra
$P$ of ${\mathcal P}(\om)$  including the finite sets to $[0,1]_{\mathbb R}$
such that
\begin{enumerate}
\item[$\bullet$] $\Xi(\emptyset) =0$, $\Xi(\om)=1$,
\item[$\bullet$] $\Xi$ is finitely additive, that is: 
If $Y, Z  \in P$ are disjoint, then $\Xi(Y \cup Z) = \Xi(Y) + \Xi(Z)$.
\item[$\bullet$]  $\Xi(\{n\}) = 0$ for $n \in \om$.
\end{enumerate} 

Members of $\mathcal M$ are called partial finitely additive measures.

\smallskip

2) ${\mathcal M}^{\rm full}$ is the set of $\Xi \in {\mathcal M}$ whose domain is
${\mathcal P}(\om)$, and the members of ${\mathcal M}^{\rm full}$ are called 
finitely additive measures. 

\smallskip

3) We write ``$\Xi(A)=a$'' (or $>a$ or whatever) if $A \in \dom(\Xi)$ and
$\Xi(A)=a$ (or $>a$ or whatever).
\end{definition}

For extending finitely additive measures we are going to use:

\begin{theorem}\label{3.2}
[Hahn Banach] 
Suppose that $\Xi$ is a partial
finitely additive measure on a algebra $P$ and that 
$X \not\in P$. Let $a \in [0,1]$ be such that 
$$\sup\{ \Xi(A) \such A \subseteq X, A \in P \} \leq a \leq
\inf\{ \Xi(B) \such  B \supseteq X, B \in P\}.$$
Then there exists a finitely additive measure $\Xi^\ast$ 
extending $\Xi$ and such that $\Xi^\ast(X) = a$.
\proofend
\end{theorem}

\begin{proposition}\label{3.3}

Let $\alpha^*$ be an ordinal.
Assume that $\Xi_0 \in {\mathcal M}$ and that for $\alpha < \alpha^*$, 
$A_\alpha 
\subseteq \om$ and $0 \leq a_\alpha \leq b_\alpha \leq 1$, 
$a_\alpha$, $b_\alpha$ reals. Then we have that

\begin{itemize}
\item (1) $\Rightarrow $ (2)
\item (2) $\Rightarrow $ ((3.A) with all $b_\alpha =1$)
\item (3.A) $\Leftrightarrow $ (3.B),
\end{itemize}

where  
\begin{myrules}

\item[(1)] If $A^* \in \dom(\Xi_0)$, $\Xi_0(A^*) > 0 $ and $n \in \om$ and
$\alpha_0 < \dots < \alpha_{n-1} < \alpha^*$ then $ A^* \cap 
\bigcap_{\ell < n} A_{\alpha_\ell} \neq \emptyset$.

\item[(2)] $\forall \varepsilon >0$, $\forall A^* \in \dom(\Xi_0)$ such 
that  $\Xi_0(A^*) >0$, $n \in \om$, $\alpha_0 <
\dots < \alpha_{n-1} < \alpha^*$  we can find a finite non-empty
$u\subseteq A^*$ such that for $\ell \in n$ 
$$ a_{\alpha_\ell} - \varepsilon \leq \frac{|A_{\alpha_\ell} \cap u|}
{|u|}.$$

\item[(3.A)]
There is $\Xi \in {\mathcal M}^{\rm full}$ extending $\Xi_0$ such that 
$\forall \alpha < \alpha^* \; \Xi(A_\alpha) \in [a_\alpha,b_\alpha]$.

\item[(3.B)]
for all $\varepsilon >0$, for all $k \in \om$, 
for all $ \langle A^*_0, \dots A^*_{m-1}\rangle$
partition of $\om$ and $A^*_i \in \dom(\Xi_0)$ such 
that  $\Xi_0(A^*_i) >0$, $n \in \om$, $\alpha_0 <
\dots < \alpha_{n-1} < \alpha^*$  we can find a finite non-empty
$u\subseteq \om \setminus k$ such that for $\ell \in n$ and $i \in m$
\begin{eqnarray*} 
a_{\alpha_\ell} - \varepsilon \,\leq & \frac{|A_{\alpha_\ell} \cap u|}
{|u|} & \leq \, b_{\alpha_\ell} + \eps,\\
\Xi_0(A^*_i) - \varepsilon \, \leq & \frac{|A^*_i \cap u|}
{|u|} & \leq \, \Xi_0(A^*_i) + \eps.
\end{eqnarray*}

\end{myrules}
\end{proposition}
\proof
(1) $\Rightarrow$ (2): Given $\eps, A^*, \alpha_0, \alpha_1, \dots 
\alpha_{n-1}$ we take $k \in A^* \cap \bigcap_{\ell<n} A_{\alpha_\ell}$ 
and $u =\{k\}$.

\smallskip

(2) $\Rightarrow$ (3.B) with $b_\alpha =1$: 
Given $\eps, k, A^*_0, \dots A^*_{m-1}$, pairwise disjoint
with positive $\Xi_0$ measure, $ 
\alpha_0, \alpha_1, \dots \alpha_{n-1}$
then we can find finite $u_i$, $i <m$ such that
\begin{eqnarray*}
u_i &\subseteq & \om\setminus k,\\
u_i & \subseteq & A^*_i,\\
\frac{|u_i|}{|\bigcup_{i \in m} u_i|} & \in & (\Xi_0(A^*_i)-\eps, \Xi_0(A^*_i) 
+\eps), \\
a_{\alpha_\ell} -\eps &\leq & \frac{|A_{\alpha_\ell} \cap u_i|}{|u_i|}.
\end{eqnarray*}
It is now easy to check that $u=\bigcup_{i<m} u_i$ is as required.

\smallskip

(3.B) $\Rightarrow$ (3.A): 
This is the special case of a symmetrized variant of (\ref{3.6} with
$a^\alpha_\ell = 1 $ iff $\ell \in A_\alpha$ and $a^\alpha_\ell
=0$ else). This is the most important implication. Its proof is
not circular, it just more economic to do \ref{3.4},
\ref{3.5}, and \ref{3.6} first.

\smallskip

(3.A) $\Rightarrow$ (3.B):
Fix $\eps'$ such that $2\ell m \eps' \leq \eps$. 
We put for $i < m $ and $\ell < n$  the first
$$
\left\lceil \frac{\Xi(A^*_i \cap A_{\alpha_\ell})}{\eps'}\right\rceil
$$
elements of $A^*_i \cap A_{\alpha_\ell}$ into $u$ (and nothing else).
It is important to see that the tasks for the different $A_{\alpha_\ell}$
can be simultaneously fulfilled. Best look for each $i<m$
at the atoms in the
Boolean algebra generated by the $A_{\alpha_\ell} \cap A^*_i$, 
$\ell <n$. 

For a real $x$, $\lceil x \rceil$ is the least integer greater 
than or equal $x$.
Then it is an easy computation that the $\frac{|A_i^* \cap u|}{|u|}$ and
the $\frac{|A_{\alpha_\ell} \cap u|}{|u|}$ are in the right intervals of
width $2 \eps$.
\proofendof{\ref{3.3}}

\smallskip

In order to convey information to later stages of our
forcing iteration, we are going to use averages. These are integrals
of functions from $\om$ to with respect to finitely additive measures. 
If the average of some function is large 
then we can go back to some finite subset 
of $\om$ where the function takes large values.

\begin{definition}\label{3.4}
1) For $\Xi \in {\mathcal M}^{\rm full}$ and a sequence $\bar{a} =
\langle a_\ell \such \ell \in \om \rangle$ of reals in $[0,1]_{\mathbb R}$
(or just $\sup_{\ell \in \om }|a_\ell| < \infty$) we let
\begin{equation*}\begin{split}
\av{\Xi}(\bar{a}) = & 
\sup\left\{ \sum_{k <k^*} \Xi(A_k) \inf(\{a_\ell \such \ell \in A_k\})
\such \langle A_k \such k <k^* \rangle \mbox{ is a partition of } \om \right\}
\\
= &
\inf\left\{ \sum_{k <k^*} \Xi(A_k) \sup(\{a_\ell \such \ell \in A_k\})
\such \langle A_k \such k <k^* \rangle \mbox{ is a partition of } \om \right\}.
 \end{split}
\end{equation*}

(Think of $A_k = \{ \ell \such a_\ell \in [\frac{k}{2^n},\frac{k+1}{2^n})\}$
and $n \to \infty$, then it is easy to see that both are equal.)

\smallskip

2) For $\Xi \in {\mathcal M}$, $A \subseteq \om$ such that $\Xi(A) > 0$ 
define $\Xi_A(B) =
\Xi(A\cap B)/ \Xi(A)$ and $\av{\Xi}(\langle a_k \such k \in B \rangle ) =
\av{\Xi_B}(\langle a_k' \such k \in \om \rangle)$ with
$$ a_k' = \left\{ \begin{array}{ll}
a_k,& \mbox{ if } k \in B,\\
0, & \mbox{ if } k \not\in B.
\end{array}
\right.
$$
\end{definition}

\begin{proposition}\label{3.5}
Assume that $\Xi \in {\mathcal M}^{\rm full}$ and 
$a_\ell^i \in [0,1]_{\mathbb R}$
for $i < i^* \in \om$, $\ell \in \om$, $B \subseteq \om$, $\Xi(B) > 0$ and
$\av{\Xi_B}(\langle a_\ell^i \such \ell < \om \rangle) = b_i$ for
$i<i^*$, $m^* < \om$ and lastly $\eps > 0$. Then for some finite $u 
\subseteq B \setminus m^*$ we have: If $i<i^*$ then
$$
b_i -\eps < \frac{\sum\{a_\ell^i \such \ell \in u \}}{|u|} < b_i + \eps.$$
\end{proposition}

\proof Let $j^* \in \omega$
and $\langle B_j \such j<j^* \rangle$ be a partition of $B$ 
such that for every $i < i^*$ we have 
$$
\left(\sum_{j < j^*} \sup \{ a_\ell^i \such \ell \in B_j \} \Xi(B_j)\right) 
-\left(\sum_{j < j^*} \inf\{ a_\ell^i \such \ell \in B_j \} \Xi(B_j)\right) 
< \frac{\eps}{2}.
$$
Now choose $k^*$ large enough such that there are $k_j$ satisfying $k^*
= \sum_{j < j^*} k_j$ and for $j < j^*$ 
$$
\left|\frac{k_j}{k^*} - \frac{\Xi(B_j)}{\Xi(B)}\right| < \frac{\eps}{2}.$$
Let $u_j \subseteq B_j \setminus m^*$, $|u_j| = k_j$ for $j < j^*$. 
Now let $u = \bigcup_{j < j^*} u_j$ and calculate
\begin{eqnarray*}
\sum_{\ell \in u} \frac{a_\ell^i}{|u|} &=& \sum_{j< j^*} \sum_{\ell \in u_j} 
\frac{a_\ell^i}{|u|} \leq \sum_{j < j^*} \sup\{ a_\ell^i \such \ell \in B_j\} 
\frac{k_j}{k^*}\\
&\leq &\sum_{j < j^*} \sup\{|{a_\ell^i \such \ell \in B_j}\} 
\left(\frac{\Xi(B_j)}{\Xi(B)} + \frac{\eps}{2j^*}\right) \leq b_i 
+ \frac{\eps}{2} + \frac{\eps}{2} = b_i + \eps;
\end{eqnarray*}
\begin{eqnarray*}
\sum_{\ell \in u} \frac{a_\ell^i}{|u|} &=& \sum_{j< j^*} \sum_{\ell \in u_j} 
\frac{a_\ell^i}{|u|} \geq \sum_{j < j^*} \inf\{ a_\ell^i \such \ell \in B_j\} 
\frac{k_j}{k^*}\\
&\geq & \sum_{j < j^*} \inf\{|{a_\ell^i \such \ell \in B_j}\} \left(
\frac{\Xi(B_j)}{\Xi(B)} -\frac{\eps}{2j^*} \right)
\geq b_i - \frac{\eps}{2} - \frac{\eps}{2} = b_i - \eps.
\end{eqnarray*}
\proofendof{\ref{3.5}}

\begin{fact}\label{3.6}
Assume that $\Xi$ is a partial finitely additive measure and $\bar{a}^\alpha
=     \langle a_k^\alpha \such k \in \om \rangle$ is a sequence of reals 
for $\alpha  < \alpha^*$ such that $\limsup_{k\to \om} |a_k^\alpha |
< \infty $ for each $\alpha$. Then $(B) \Rightarrow (A)$.
\begin{myrules}
\item[(A)]
 There is $\Xi^* \supseteq \Xi$, $\Xi^* \in {\mathcal M}^{\rm full}$ such that
$\av{\Xi^*}(\bar{a}^\alpha) \geq b_\alpha$ for $\alpha < \alpha^*$.
\item[(B)]
For every partition $\langle B_0, \dots B_{m^* -1}\rangle$ of $\om$
with $B_m \in \dom(\Xi)$ and $\eps > 0$, $k^* > 0$ and $\alpha_0 <  \dots < 
\alpha_{n-1} < \alpha^*$ there is a finite $u \in \om \setminus k^*$ such 
that \begin{myrules}
\item[(i)] $\Xi(B_m) - \eps 
< \frac{|B_m\cap u|}{|u|} < \Xi(B_m) + \eps$.
\item[(ii)] $\frac{1}{|u|} \sum_{k \in u} a_k^{\alpha_\ell} > 
b_{\alpha_\ell} - \eps$ for $\ell < n$.
\end{myrules}
\end{myrules}
\end{fact}

\proof

We take $$\Delta =
[\{ \mbox{partitions } \langle B_0, \dots B_{m^*-1} \rangle \mbox{ of }
\dom(\Xi) \} \times (0,1] \times \om \times [\alpha^*]^{<\om}]^{<\om}.$$
and take a filter ${\mathcal F} \subseteq {\mathcal P}(\Delta)$ such that for each
$$
\bar{c} \in 
\{ \mbox{partitions } \langle B_0, \dots B_{m^*-1} \rangle \mbox{ of }
\dom(\Xi) \} \times (0,1] \times \om \times [\alpha^*]^{<\om}
$$ 
we have that
$$ 
\{ F \in \Delta \such \bar{c} \in F \} \in {\mathcal F}.
$$
For each $F \in \Delta$ we choose $u(F)$ fulfilling
the tasks (B) simultaneously for all $\bar{c} \in F$, i.e.\
(i) and (ii) of (B) hold for $u(F)=u$, $\bar{c}(0) =
\langle B_0, \dots B_{m^*-1} \rangle$, $\bar{c}(1) 
= \eps$, $\bar{c}(2) = k^\ast$, $\bar{c}(3) =
\{ \alpha_0, \dots, \alpha_{n-1}\}$.

Then we take an ultrafilter ${\mathcal U} \supseteq {\mathcal F}$ and set for
$A$ in the algebra $\mathcal A$ generated by $\left\{\{k \such a_k^{\alpha} 
\in [q,q']\}
\such \alpha < \alpha^*, 0\leq q \leq q' \leq 1\right\} \cup
\dom({\Xi})$: 
$$\Xi^*(A) = \mbox{ the standard part of }
\left( 
\left. 
\left\langle \frac{|u(F) \cap A|}{|u(F)|}\such F \in 
\Delta \right\rangle 
\right/ 
{\mathcal U}\right).$$
By the Hahn Banach Theorem, there is an extension of
$\Xi^* $ to ${\mathcal P}(\om)$.
\proofendof{\ref{3.6}} 

An important application of \ref{3.3} (and the
hard part thereof, which is only proved in \ref{3.6}) is:

\begin{claim}\label{3.7a}
Suppose that $Q_1, Q_2$ are forcing notions in $V$, $\Xi_0 \in {\mathcal 
M}^{\rm full} $ in $V$, $\Vdash_{Q_\ell} \mbox{``}\Name{\Xi}_\ell$ is a 
finitely additive measure  extending $\Xi_0$ for $\ell = 1,2,$''.
Then $\Vdash_{Q_1\times Q_2} $ ``there is a finitely additive measure 
extending $\Name{\Xi}_1$ and $\Name{\Xi}_2$ (and hence $\Name{\Xi}_0$)''.
\end{claim}

\proof 
We are going to show, that 
$\Vdash_{Q_1 \times Q_2}$ `` $\Name{\Xi}_1$ (in the r\^ole of 
$\Xi_0$ of \ref{3.3}) and $\{A^*_\alpha \such A^*_\alpha  
\in V^{Q_2} \cap {\mathcal P}(\om)\}$ (in the r\^ole of $\langle
A^*_\alpha \such \alpha < \alpha^* \rangle$ of \ref{3.3})
fulfil (3.B) of \ref{3.3}''.

\smallskip

First we show that 
$$\Vdash_{Q_1 \times Q_2}  \dom(\Name{\Xi}_1)
\cap \dom(\Name{\Xi}_2) = \dom({\Xi}_0) = \check{V} \cap {\mathcal P}(\om).$$

So assume that we have an $Q_1$-name $\Name{X}$ and a $Q_2$-name
$\Name{Y}$ such that
$\Vdash_{Q_1 \times Q_2} \Name{X} = \Name{Y}$.

Let $Z = \{ n \in \om \such \exists 
p \in Q_1 \; p \Vdash_{Q_1} n \in \Name{X} \}$. The set $Z $ is in $V$ and
$\Vdash_{Q_1} \Name{X} \subseteq Z$. It is easy to see that
$\Vdash_{Q_2} Z \subseteq \Name{Y}$. So we get
$$\Vdash_{Q_1\times Q_2} \Name{X} \subseteq Z
\subseteq \Name{Y} = \Name{X},$$
and our first claim is proved.

\smallskip

Now we check (3.B). Let $\eps$, $k$, $\langle 
A^*_i \in V^{Q_1} \such i < m \rangle$ 
a partition of $\om$ and $\alpha_\ell$, $\ell < n$ be given. 
W.l.o.g.\ the $A_{\alpha_\ell} \in V^{Q_2}$ are a partition of $\om$
as well.

If for some $i, \ell$
$$\Vdash_{Q_1 \times Q_2}  \name{A_i^*} 
\cap \name{A_{\alpha_\ell}} \mbox{ is finite},$$

then $A_i^*$ and $A_{\alpha_\ell}$ can be separated by some $A \in V$.
This is shown in a manner similar to the proof of the first claim.

We  choose a separator $A^{i,\ell} \in V$ for each $i, \ell$ such that
$\Vdash_{Q_1 \times Q_2}  \name{A_i^*} 
\cap \name{A_{\alpha_\ell}} \mbox{ is finite}$ and let $A^j$, $j <j^*$
be the partition of $\om$ in $V$ that is generated by all the 
$A^{i, \ell}$.

\smallskip

Then, we set $\eps' = \frac{\eps}{m n  j^*}$ and 
 put   for each  $i, \ell, j$ such that
$$\Vdash_{Q_1 \times Q_2}  \name{A_i^*} 
\cap \name{A_{\alpha_\ell}} \cap \check{A^j} \mbox{ is  infinite},$$
in the 
forcing extension $V^{Q_1 \times Q_2}$, the first
$$\left\lceil \frac{{\Xi}_1(A^*_i \cap A^j ) 
\times {\Xi}_2(A_{\alpha_\ell} \cap
A^j)}{\eps' \times \Xi_0(A^j)} \right\rceil
$$
elements of  $A_i^* \cap A_{\alpha_\ell}\cap A^j$ (and no further points)
into $u$.

\proofendof{\ref{3.7a}}

%
%

\section{The First Part of the Proof of $(**)_{\bar{Q}}$: Introduction of
${\mathcal K}^3$}\label{S4}

In order to prove $(**)_{\bar{Q}}$, we need that for suitable 
 $\bar{Q}  =\langle P_\alpha, 
\Name{Q}_\beta, A_\beta, \name{\tau}_\beta, \mu_\beta,
 \such \beta < \lgg(\bar{Q}), \alpha \leq \lgg(\bar{Q})
\rangle$  from $\mathcal K$ (see Definition \ref{2.2})
 we have almost (in the sense explained in the proof of~\ref{5.5})
an expansion of the form
$$\bar{Q}^{exp}  =\langle P_\alpha, 
\Name{Q}_\beta, A_\beta, \name{\tau}_\beta, \mu_\beta,
 \eta_\beta,
(\Name{\Xi}_\alpha^t)_{t \in {\mathcal T}} \such \beta < \lgg(\bar{Q}),
\alpha \leq \lgg(\bar{Q})
\rangle$$ such that $\bar{Q}^{exp}$ is in
a special class ${\mathcal K}^3$, which we shall define in
Definition~\ref{4.2}.

\smallskip

In order to introduce ${\mathcal K}^3$, we shall
first define and (try to) explain the set 
$\mathcal T$ of blueprints (Definition~\ref{4.1}).
For each blueprint $t$ and $\alpha < \alpha^*$ the $\name{\Xi}^t_\alpha$
will be $P_\alpha$-name for some finitely additive measure on ${\mathcal P}(\om)$ 
that conveys some information about $\om$-tuples $\langle p_k \such k \in \om
\rangle$ of conditions that fit well to the blueprint $t$, 
from stage $\alpha$ to later stages in the iteration.

\smallskip
     
Let us tell more about the ideas of the proof of $(**)_{\bar{Q}}$:
In Lemma~\ref{2.12},  if the $p_\zeta$ are not all Cohen, 
the premise 3 is hard to fulfil. Think 
of $\kappa^+$ 
many $p_\zeta$ being given, so that we can do many 
thinning out procedures and have them similar, i.e.\ similar partial random 
conditions and Cohen conditions. Then we 
keep only the first $\omega$ of the $\zeta$'s and the first 
$\omega$ conditions $\langle p_\zeta \such \zeta \in \om \rangle$.
We try to strengthen them a little bit (to $p'_\zeta$)
and then get that the strengthened conditions allow to define one 
condition $p^\otimes \geq p^*$ such that 
$$
p^\otimes \Vdash \mbox{``}\bigcap_{\zeta \in \Name{E}
=\{ \zeta \such p'_\zeta \in \Name{G}\}} 
{\rm tree}_{\ell^*} (\name{\bar{a}}^{\alpha_\zeta}) 
\mbox{ has finitely many
branches''}
$$
and hence cannot contain a perfect tree.
There are some requirements on $\langle p_\zeta \such \zeta \in \om \rangle$,
as they have to predict some probabilities about the branches of
the ${\rm tree}_{\ell^*}  (\name{\bar{a}}^{\alpha_\zeta})$
and about the subset of the $\{ p'_\zeta \such \zeta \in \om \}$,
that lies in $G$.

\smallskip

The technical means to allow these predictions is the
use of finitely additive measures and the properties
(e) to (i) in the definition of ${\mathcal K}^3$.
These items in the definition have long premises by themselves. However
the premises are sufficiently often fulfilled if we start with $\kappa^+$
many  $p_\zeta$,
thin out, and choose an appropriate $t \in {\mathcal T}$.

 \smallskip

We embark with the definition of a blueprint $t$. The set of all
blueprints is denoted by $\mathcal T$. The reader may think that
$t$  describes some relevant information about 
the chosen tuples $\langle p_\zeta \such 
\zeta \in \om \rangle$.
Later it will turn out that sequences described by the same $t$ are 
compatible forcing conditions (though we have finite supports and
are not interested in taking
the union of countably many conditions). 
This will be used in Lemma~\ref{4.8}.

\smallskip

In the case of iterations where 
all Cohen forcings are just those forcings
 in an initial segment of the iteration
(as in \ref{2.2} Part 2)), we can dispense with
the parameter $\m$ in the next definition.
This simplification  is not worthwhile because the
generality allows another application of the method:
In Section~\ref{S6}, we shall work with a type of iteration where Cohens 
are added cofinally often.

\smallskip

However, we could  simplify \ref{4.2} slightly  and leave out (f) there
in the special case that 
the $f_\zeta$ of \ref{2.11} move only one $\alpha$ in the Cohen part
and leave the indices at which partial randoms are attached fixed.
We do not simplify because we hope for future applications. 

\begin{definition}\label{4.1}
We fix a $\kappa$ such that $2^\kappa \geq \chi$ (from \ref{2.2}).
The set $\mathcal T$ of blueprints is the set of tuples
$$t = (w^t, \n^t, \m^t, \bar{\eta}^t, h_0^t, h_1^t, h_2^t, \bar{n}^t)$$
such that
\begin{myrules}
\item[(a)] 
$w^t \in [\kappa]^{\aleph_0}$.
(What is the purpose?
Think of the latter as $[\chi]^{\aleph_0}$ disguised.
Suppose that $|\dom(p_\zeta)| = \n^t$ for all $\zeta$,
$\dom(p_\zeta) = \{\gamma_\zeta^i \such i < \n^t\}$, 
$\langle \gamma^i_k \such k \in
\omega \rangle \in \chi^\omega$ for each fixed $i <  \n^t$,
but $\chi \leq 2^\kappa$ and we can fix an injection and keep
as relevant information certain parts of $\kappa$ coming from
 of certain $f \in 2^\kappa$. Look at the $w^t$ in Subclaim~\ref{5.3}.)

\item[(b)] $0 < \n^t < \omega$, $0 \leq \m^t \leq \n^t$.
($\n^t$ will be the cardinality of the heart of the $\Delta$-system
built from many $p_\zeta$ and $\m^t$ will be the cardinality of 
the part of the heart that is lying below $\chi$.)

\item[(c)] 
$\bar{\eta}^t = \langle \eta^t_{n,k} \such
n < \n^t, k \in \om \rangle$, $\eta_{n,k}^t \in \,^{w^t}2$
($\eta_{n,k}^t $ codes the $n$th element of the support 
 of $p_k$ for $k \in \omega$  and these $k$ 
are the first $\om$ of the $\zeta$).

\item[(d)] $h_0^t$ is a partial function from $[0,\n^t)$ to $\kappa$
\footnote{We do carry out the simplification suggested in
a footnote in \cite{Shelah592} and take $\kappa$ instead of
$^\omega \kappa$ here. This does not bring any
disadvantages, because when choosing $\langle p_\zeta
\such \zeta \in \omega \rangle$ we have initially $\kappa^+$ many
$p_\zeta$, and hence can thin out such that for each
$\zeta$, $| \dom p_\zeta|$ is the same, say $\n^t$, and that for auch
$n< \n^t$, the $p'_\zeta(\mbox{$n$th element of }\dom(p'_\zeta))=
h^t_0(n)$ are independent of
$\zeta$, if they lie in some notion of forcing with
conditions in  some $Q_\alpha$ with $|Q_\alpha| <\kappa$.}
($\dom(h_0^t)$ is the part of those $\alpha$
in  the heart of the $\Delta$-system where 
$\Name{Q}_\alpha$ is the Cohen forcing. In the somewhat
simpler case of \ref{2.2} Part 2), this
domain coincides with the part of the heart that lies below $\chi$.)

\item[(e)] $h_2^t$ is a function from $[0,\n^t) \setminus \dom(h_0^t)$ to
$^{<\om} 2$. (Think of $h_2^t$ giving some information of
a partial random condition attached at some point of the heart.)

\item[(f)]
$h_1^t$ is a function from $[0,\n^t)$ into the rational interval
 $[0,1)_{\mathbb Q}$, such that $\{ n \such h_1^t(n) \neq 0 \} 
\subseteq \dom(h_2^t)$. Furthermore we have that
$\sum_{n < \n^t} \sqrt{h_1^t(n)} < \frac{1}{10}$. 
(Think of $h_1^t$ giving some   information about the Lebesgue measure of
the limit of the
a partial random condition attached 
at some point of the heart intersected with
$\dom(h_2^t)$.)

\item[(g)] $\eta^t_{n_1, k_1} = \eta^t_{n_2, k_2} \Rightarrow n_1 = n_2$
(This is some compatibility requirement, which is useful in
\ref{4.5}.)

\item[(h)] For each $n <\n^t$
 we have that $\langle\eta^t_{n,k} \such k 
\in \om \rangle$ is either constant or with no repetitions
(that is: either
 in the heart of the system or among the moved parts of the
domains of the $\langle p_k \such k \in \omega \rangle$).

\item[(i)] $\bar{n}^t = \langle n^t_k \such k \in \om \rangle$
where $n_0^t = 0$, $n^t_k < n^t_{k+1} < \om$ 
and the sequence $\langle n^t_{k+1}  -n^t_k \such k \in \om \rangle$ 
goes to infinity. 
(This last ingredient does not describe
$p_\ell$ but is just an additional part handling the
finitely additive measures $\name{\Xi}^t_\alpha$. The sequences
$\bar{n}^t$
 shall allow to
compute intersections of sets of branches from {\em lim tree}, and for
these computations (see \ref{5.3}) the
$p_\ell$ are grouped together for $\ell \in [n^t_k, n^t_{k+1})$.)
\end{myrules}
\end{definition}

There are $\kappa^\om$ many blueprints. 
(Remember we also require that $2^\kappa \geq \chi$, otherwise the
choice of the 
$\eta$ in the following definition would fail.)

\smallskip

\noindent{\bf Explanation:} We continue the explanations 
begun in the parentheses in order to explain how the conditions 
shall work together:

As mentioned, $(**)_{\bar{Q}}$ 
follows from the fact that in $V^{P_{\alpha^*}}$, if $E \in
[\chi]^{\kappa^+}$ and $m \in \om$, then
$\bigcap_{\alpha \in E} {\rm tree}_m(\bar{a}^\alpha)$ is a tree 
with finitely many branches.
Suppose some $p$ forces the contrary. We take $p_\zeta \geq p$ such that
$p_\zeta \Vdash$ ``$\betta_\zeta \in E$''
for $\zeta \in \kappa$
and such that $\betta_\zeta \not\in \{\betta_\xi \such \xi < \zeta \}
$.

\smallskip

We can assume that the $p_\zeta$ are in some given dense set
(will be ${\mathcal I}_{\bar{\eps}}$ of \ref{5.1} in our case)
and that the $\langle p_\zeta \such \zeta \in \kappa^+\rangle$ form
a $\Delta$-system with some additional thinning demands, putting
$\kappa^+$ many objects into less than $\kappa$ many
pigeonholes. (See our earlier remarks about working with $\kappa^+$ many
$\zeta$ and the proof of Lemma~\ref{5.2}.)
\smallskip

We assume that $\dom(p_\zeta) =
\{ \gamma_{\n,\zeta} \such \n < \n^t \}$,
$\gamma_{\n,\zeta}$ is increasing in $\n$ and $\gamma_{\n,\zeta}
< \chi$ iff $\n < \m^t$ and that $\betta_\zeta$ is one of the
$\gamma_{\n,\zeta}$.
We let $p_\zeta'$ be $p_\zeta$ except that $p_\zeta(\betta_\zeta)$ is
increased a little.

\smallskip

It suffices to find some $p^\otimes \geq p$ such that
$p^\otimes \Vdash $ ``$\Name{A} = \{ \zeta \in \om
\such p_\zeta' \in \name{G } \}$ is `large enough' such that 
$\bigcap_{\zeta \in \Name{A}} {\rm tree}_m(\name{\bar{a}}^{\betta_\zeta})$
has only finitely many branches''.

\smallskip

The `large enough' is interpreted in terms of a $\Xi^t_\alpha$-measure.

\smallskip

The $\n < \n^t$ such that $Q_{\gamma_{\n,\zeta}}$ is a forcing notion of
cardinality $<\kappa$ (in our forcings, then it is just the Cohen forcing)
do not cause problems because $h_0^t(\n)$ tells us exactly what the
condition is. Still there are many cases of such $\langle
p_\zeta \such \zeta \in \om \rangle$ which fall into the same $t$,
and we will get contradictory demands if $\gamma_{\n_1,\zeta_1}
=\gamma_{\n_2,\zeta_2}$ and $\n_1 \neq \n_2$.
But the $w^t$, $\bar{\eta}^t$ are built in order
to prevent this. That is we have to assume that $2^\kappa
\geq \chi$ in order to be able 
to choose $\langle \eta_\alpha \such \alpha \in \chi 
\rangle$, $\eta_\alpha \in 2^\kappa$ 
with no repetitions and such that for
 $v \subseteq \chi$, $|v| \leq \aleph_0$ 
(in the applications, we shall have
 $v= \{ \alpha_{\n,\zeta} \such \zeta \in \om \}$)
there is some $w=w^t \in [\kappa]^{\aleph_0}$ such that
$\langle \eta_\alpha \restriction w \such \alpha \in v \rangle$ is
without repetitions.

\smallskip

So the blueprint $t$ describes such a situation giving much information, 
though the number of blueprints is $\kappa^\omega$.

\smallskip

If $Q_{\alpha_{\n,\zeta}}$ is partial random, we get many different
possibilities
for $p_\zeta(\gamma_{\n,\zeta})$, too many to apply a pigeonhole
principle.
We want that many of them will lie
in the generic set.  Using $(h_1^t(\n), h_2^t(\n))$ we 
know that in the interval $(^\om 2)^{[h_2^t(\n)]}$ the 
set $\lim(p_\zeta(\gamma_{\n,\zeta}))$ is of relative
measure $\geq 1-h_1^t(\n)$.
Still there are too many (possibly incompatible) 
$p_\zeta(\gamma_{\n,\zeta})$
and finally, in \ref{5.2} and\ref{5.3},
 the existence of many compatible candidates is ensured by
the finitely additive measures.

\smallskip

The $\bar{n}^t= \langle n_k^t \such k \in \om \rangle$ 
are going to be used
in the end of Section \ref{S5}, where we show that
$\{ \zeta \such p'_\zeta \in G \}$ is large by showing that
for infinitely many $k$ we have that
 $$\frac{|\{ \zeta \such n^t_k \leq \zeta < n^t_{k+1} \mbox{ and }
p'_\zeta \in G \}|}{n^t_{k+1} - n^t_k}$$
is large, say $> \eps>0$.

\smallskip

The $n_k^t$ will be chosen such that they are increasing fast enough with $k$
and $\langle p'_\zeta(\gamma_{\n,\zeta} ) \such \zeta \in [n_k^t,n_{k+1}^t)
\rangle$ will be chosen such that for each $\eps > 0$
there is some $s\in \om$ such that for $k$ large enough:
if the above fraction is above $\eps$ then
$$
^k 2 \cap \bigcap \{ {\rm tree}_m(\bar{a}^{\beta_\ell}) \such n_k^t \geq \ell
< n_{k+1}^t
\mbox{ and } p'_\ell \in G \}$$
has $< s$ members, hence the tree has fewer than $s$ branches.

\smallskip

{\em Comment on simplifications:\/}
Now we finally define the kind of iteration we use
for the proof of $(**)_{\bar{Q}}$. 
The reader who is longing for some simplification may  omit
the condition (f)
in \ref{4.2}, \ref{4.5} and \ref{5.3} and work
just with conditions $p_\zeta$ that do not
differ at any index in the iteration where a partial 
random real is attached to it, but only at those indices
where a forcing of size less than $\kappa$ is attached, 
or even work with with $p_\zeta$ that differ only
at $\betta_\zeta < \chi$ (from \ref{2.11}). A look at the beginning of
\ref{5.2}, where the $p_\zeta$ and $p'_\zeta$ are chosen, 
and a look $AUT(\bar{Q})$ shows that the restriction to this simplified
 situation is always possible when forcing with
a member of the restricted class described in Definition~\ref{2.2}
Part 2.

\begin{definition}\label{4.2}
${\mathcal K}^3$ is the class of sequences 
$$\bar{Q} = \langle P_\alpha, \Name{Q}_\beta, A_\beta, \mu_\beta,
\name{\tau}_\beta, \eta_\beta, (\Name{\Xi}^t_\alpha)_{t \in {\mathcal T}}
\such \alpha \leq \alpha^*, \beta < \alpha^* \rangle$$ 
(we write $\alpha^* = \lgg(\bar{Q})$) 
such that
\begin{myrules1} 

\item[(a)] 
$$\bar{Q} = \langle P_\alpha, \Name{Q}_\beta, A_\beta, \mu_\beta,
\name{\tau}_\beta, 
\such \alpha \leq \alpha^*, \beta < \alpha^* \rangle$$
is in $\mathcal K$ from Definition~\ref{2.2}.

\item[(b)] 
$\eta_\beta \in \, {^\kappa  2}$ and for $\beta < \alpha < \alpha^*$ we have that
$\eta_\beta \neq \eta_\alpha$. 

\item[(c)] ${\mathcal T}$ is the set of all blueprints, and
$\Name{\Xi}^t_\alpha$ is a $P_\alpha$-name for a finitely additive measure 
in $V^{P_\alpha}$, increasing with $\alpha$. 

\item[(d)] We say the $\langle \alpha_\ell \such \ell \in \om \rangle$
satisfies $(t,\n)$ for $\bar{Q}$, if

\smallskip

(Think of $p_{\ell}$ being the first $\omega$ of the 
$p_\zeta$ and
$\langle \alpha_\ell
\such \ell \in \omega \rangle =
\langle \gamma_{\n,\zeta} \such \zeta \in 
\omega \rangle$, and in particular, 
$\langle \alpha_\ell \such \ell \in \omega \rangle =  \betta_\ell
\such \ell \in \omega$ from \ref{2.10}. 
($\alpha_\ell$ is for some $\n$ always
the $\n$th element in $\dom(p_\ell)$)
Further think that the following items also mean that  
$\langle p_{\ell} \such \ell \in \om
\rangle$ being sufficiently described by $t \in {\mathcal T}$)

\begin{enumerate}
\item $\langle \alpha_\ell \such \ell \in \om \rangle \in V$,

\item $t \in {\mathcal T}$, $\n < \n^t$,

\item $\alpha_\ell < \alpha_{\ell +1} < \alpha^*$,

\item $\n < \m^t \Leftrightarrow \forall \ell (\alpha_\ell
< \chi) \Leftrightarrow \exists \ell (\alpha_\ell < \chi)$
(the moved positions $\alpha_\ell$ are in the Cohen part),

\item $\eta^t_{\n,\ell} = \eta_{\alpha_\ell} \restriction w^t$.
($\eta_{\alpha_\ell}$ describes where $\alpha_\ell$ really is, and
$\eta^t_{\n,\ell}$ describes a part of it of size $\om$. For a given $t$,
the $\n$ such that $\bar{Q}$ satisfies $(t,\n)$ is unique by \ref{4.1}
(g).),

\item If $\n \in \dom(h_0^t)$ then $\mu_{\alpha_\ell} < \kappa$ and
$\Vdash_{P_{\alpha_\ell}} ``|Q_{\alpha_\ell} | < \kappa$ and
$(h_0^t(\n))(\ell) \in \Name{Q}_{\alpha_\ell}$'',

\item If $\n \in \dom(h_1^t)$ then $\mu_{\alpha_\ell} \geq \kappa$, so
$\Vdash_{P_{\alpha_\ell}} \mbox{``} Q_{\alpha_\ell}$ 
has cardinality  $\geq \kappa$'' (hence it is partial random),

\item If $\langle \eta^t_{\n,k} \such k \in \om \rangle$ is constant, then 
$\forall \ell \: \alpha_\ell = \alpha_0$,

\item  If $\langle \eta^t_{\n,k} \such k \in \om \rangle$ is not
constant, then 
$\forall \ell \: \alpha_\ell < \alpha_{\ell+1}$.
\end{enumerate}

\item[(e)] If $\bar{\alpha} = \langle \alpha_\ell \such \ell \in \om 
\rangle$ satisfies $(t,\n)$ for $\bar{Q}$, $\bigwedge_{\ell \in \om}
(\alpha_\ell < \alpha_{\ell +1})$, $\n \in \dom(h_0^t)$ and
$$
C = \{ k \in \om \such
\forall \ell \in [n_k,n_{k+1}) \; h^t_0(\n)(\ell) \in 
G_{Q_{\alpha_\ell}} \},
$$ 
then
$$\Vdash_{P_{\alpha^*}} \name{\Xi}^t_{\alpha^*}(\Name{C}) = 1.$$

\item[(f)] If $\bar{\alpha} = \langle \alpha_\ell \such \ell \in \omega 
\rangle$ satisfies $(t,\n)$ for $\bar{Q}$, $\bigwedge_{\ell \in \om}
(\alpha_\ell < \alpha_{\ell +1})$, $\n \in \dom(h_1^t)$, 
$\name{\bar{p}} = 
\langle \name{p}_\ell \such \ell \in \om \rangle$ is such that
$\name{p}_\ell $ is a $P_{\alpha_\ell}$-name for a member of 
$\Name{Q}_{\alpha_\ell}$, and for every $\ell$,
\begin{equation}
\tag{$*$}
\Vdash_{P_{\alpha_\ell}}
1 - h_1^t(\n) \leq \frac{{\rm Leb}(\{ \eta 
\in {^\om 2} 
\such h_2^t(\n) \vartriangleleft \eta \in \lim(\name{p}_\ell) \})}
{2^{\lgg(h_2^t(\n))}}
\end{equation}
and if $\eps >0$ is such that
$$
C = \left\{ k \in \om \: \left| \:
\frac{|\{\ell \in [n^t_k,n^t_{k+1}) \such \name{p}_\ell 
\in G_{Q_{\alpha_\ell}}\}|}{n_{k+1}^t - n_k^t} 
\geq (1-h_1^t(\n)) (1-\varepsilon)\right.\right\},
$$ 
then
$$
\Vdash_{P_{\alpha^*}} \Name{\Xi}^t_{\alpha^*}\left(
\Name{C}\right) = 1.$$

\item[(g)] If $\bar{\alpha} = \langle \alpha_\ell \such \ell \in \om 
\rangle$ satisfies $(t,\n)$ for $\bar{Q}$, $\bigwedge_{\ell \in \om}
\alpha_\ell = \alpha$, $\n \in \dom(h_1^t)$, 
$\name{r}$ and $\name{\bar{r}} = \langle \name{r}_\ell 
\such \ell \in \om \rangle$ are $P_\alpha$-names  for members of
$Q_\alpha$ such that 

\begin{equation}\tag{$**$}
\begin{split}
&\;\; \mbox{ in } V^{P_\alpha}:  \forall r' \in 
Q_\alpha \: \mbox{ if } r' \geq r, \mbox{ then } \\
&\;\; {\rm Av}_{{\Xi^t}_\alpha} \left(\langle a_k(r') \such k \in \om
\rangle\right)  \geq  1 - h_1^t(\n), \mbox{ where }
\\
&\;\; a_k(r')  =  a_k(r',\bar{r}) 
= \left(\sum_{\ell \in [n_k,n_{k+1})} 
\frac{{\rm Leb}(\lim(r')\cap \lim(r_\ell))}
{{\rm Leb}(\lim(r'))}\right) \cdot \frac{1}{n^t_{k+1} -n^t_k},
\end{split}\end{equation}
then

\begin{eqnarray*} 
\Vdash_{P_{\alpha^*}} 
\mbox{``if } \name{r} \in Q_{\alpha},&& \mbox{ then }\\
1-h_1^t(\n) &\leq& {\rm Av}_{\Name{\Xi}_{\alpha^*}^t} \left(\left\langle
\left.  \frac{|\{\ell \in [n^t_k,n^t_{k+1}) \such \name{r}_\ell 
\in G_{Q_{\alpha_\ell}}\}|}{n_{k+1}^t - n_k^t} 
\: \right| \: k \in \om \right\rangle \right)\mbox{''}.
\end{eqnarray*}

\item[(h)] $P'_{A_\alpha} \lessdot P_\alpha$,

\item[(i)] For $t \in {\mathcal T}, \alpha \in \alpha^*$: If $\Vdash_{P_\alpha}
|Q_\alpha| \geq \kappa$, then $\Name{\Xi}^t_\alpha \restriction 
{\mathcal P}(\om)^{V^{P_{A_\alpha}}}$ is a $P_{A_\alpha}$-name.\footnote{This 
is where the information is whispered, showing that $\name{Q}_\alpha$,
the random forcing
over $V[\tau_\beta \such \beta\in A_\alpha ]$, behaves in the sense
of $\Xi^t_\alpha$ instead of the Lebesgue measure in a certain sense generic: 
$r_\alpha$ hits sets of large $\Xi^t_\alpha$ measure.}

\end{myrules1}
\end{definition}

\begin{definition}\label{4.3}
1. For $\bar{Q} \in {\mathcal K}^3$ and for $\alpha^*< \lgg(\bar{Q})$ let
$$\bar{Q} \restriction \alpha^*
= \langle P_\alpha, \Name{Q}_\beta, A_\beta, \mu_\beta, \name{\tau}_\beta,
\eta_\beta, (\name{\Xi}^t_\alpha)_{t \in {\mathcal T}} \such \alpha \leq 
\alpha^*, \beta <\alpha^*\rangle.
$$
2. For $\bar{Q}^1, \bar{Q}^2 \in {\mathcal K}^3$ we say:
$$\bar{Q}^1 < \bar{Q}^2 \mbox{ if }
\bar{Q}^1 = \bar{Q}^2 \restriction \lgg(\bar{Q}^1).$$
\end{definition}

In the next three steps, we show that ${\mathcal K}^3$ is sufficiently rich:
That is, if we have some $\bar{Q}$ in ${\mathcal K}^3$ then we 
can find an extension. The successor step and the 
limit step of cofinality $\om$ require some work, whereas the
 limits of larger cofinality are easy because no new reals
are introduced in these limit steps.

\begin{fact}\label{4.4}
(1) If $\bar{Q} \in {\mathcal K}^3$, $\alpha \leq \lgg(\bar{Q})$, 
then $\bar{Q} \restriction \alpha \in {\mathcal K}^3$.
\\
(2) $({\mathcal K}^3, \leq)$ is a partial order.
\\
(3) If a sequence $\langle \bar{Q}^\zeta \such \zeta < \delta \rangle$
is increasing, $\cf(\delta) > \aleph_0$, then there is a unique $\bar{Q} 
\in {\mathcal K}^3$ which is the least upper bound, $\lgg(\bar{Q}) =
\bigcup_{\zeta <\delta} \lgg(\bar{Q}^\zeta)$ and $\bar{Q}^\zeta \leq \bar{Q}$ 
for all $\zeta < \delta$.
\end{fact}

\proof Easy.

\begin{lemma}\label{4.5}
Suppose that $\bar{Q}_n < \bar{Q}_{n+1}$, $\bar{Q}_n \in {\mathcal K}^3$,
 $\alpha_n = \lgg(\bar{Q}_n)$, $\delta = \sup(\alpha_n)$.
Then there is some $\bar{Q} \in {\mathcal K}^3$ such that $\lgg(\bar{Q})= \delta$
and $\bar{Q}_n < \bar{Q}$ for $n \in \om $.
\end{lemma}

\proof We have to define $(\Name{\Xi}^t_\delta)_{t \in {\mathcal T}}$, such that 
(e) and (f) of the definitions of ${\mathcal K}^3$ hold. The items (g) and (i) 
do not produce no new tasks in the limit steps, 
and we proved (h) in \ref{2.6} and \ref{2.7}.

So, we  look again at (e) and (f) of \ref{4.2}:

\begin{myrules}

\item[(e)] If $\bar{\alpha} = \langle \alpha_\ell \such \ell \in \om 
\rangle$ satisfies $(t,\n)$ for $\bar{Q}$, $\bigwedge_{\ell \in \om}
(\alpha_\ell < \alpha_{\ell +1})$, $\n \in \dom(h_0^t)$ and
$$C = \{ k \in \om \such
\forall \ell \in [n_k,n_{k+1}) \; h^t_0(\n)(\ell) \in 
G_{Q_{\alpha_\ell}} \},$$ then
$$\Vdash_{P_{\alpha^*}} \Name{\Xi}^t_{\alpha^*}(\Name{C}) = 1.$$

\item[(f)] If $\bar{\alpha} = \langle \alpha_\ell \such \ell \in \om 
\rangle$ satisfies $(t,\n)$ for $\bar{Q}$, $\bigwedge_{\ell \in \om}
(\alpha_\ell < \alpha_{\ell +1})$, $\n \in \dom(h_1^t)$, 
$\name{\bar{p}} = 
\langle \name{p}_\ell \such \ell \in \om \rangle$ is such that
\begin{equation}\tag{$*$}
\Vdash_{P_{\alpha_\ell}}
1 - h_1^t(\n) \leq \frac{{\rm Leb}(\{ \eta 
\in ^\om 2 \such h_2^t(\n) \vartriangleleft \eta \in \lim(\name{p}_\ell) \})}
{2^{\lgg(h_2^t(\n))}},
\end{equation}
and $\eps > 0$ and 
$$
C= \left\{ k \in \om \: \left| \:
\frac{|\{\ell \in [n^t_k,n^t_{k+1}) \such \name{p}_\ell 
\in G_{Q_{\alpha_\ell}}\}|}{n_{k+1}^t - n_k^t} 
\geq (1-h_1^t(\n)) (1-\varepsilon)\right.\right\},$$
then
$$
\Vdash_{P_{\alpha^*}} \Name{\Xi}^t_{\alpha^*}\left(
\Name{C}\right) = 1.$$
\end{myrules}

By \ref{3.2} it suffices to show
\begin{eqnarray*}
 &\Vdash_{P_{\delta}} & \mbox{``if }
\Name{B} \in \bigcup_{\alpha < \delta} \dom(\Name{\Xi}^t_\alpha) =
\bigcup_{\alpha<\delta} ({\mathcal P}(\om))^{V^{P_\alpha}}
\\
&&   \mbox{and } \Xi_\alpha^t(\Name{B})> 0 \mbox{ and }
j^* \in \om \mbox { and }   \Name{C}_j, j<j^*, 
\mbox{ are sets }\\
&& \mbox{from (e) or (f) (whose measure
is required to be 1 there), }\\
&& \mbox{ then } 
\Name{B} \cap \bigcap_{j<j^*} \Name{C}_j \neq \emptyset.
\mbox{''}.
\end{eqnarray*}

Towards a contradiction, assume $q \in P_\delta$ forces the negation. 
So possibly
increasing $q$ we have: For some $\Name{B}$ and for some $j^* \in \om $, 
for each $j< j^*$ we have
$\varepsilon > 0$, and $\n(j) <\n^t$, $\langle \alpha^j_\ell 
\such \ell \in \om \rangle$, $\langle \name{p}^j_\ell \such \ell 
\in \om \rangle$ 
involved in the definition of $\Name{C}_j$ (in (e) or (f) of 
Definition~\ref{4.2}), and $q$ forces:
\begin{equation*}\begin{split}
& \Name{B} \in \bigcup_{\alpha<\delta} \dom(\Name{\Xi}^t_\alpha)
= \bigcup_{\alpha<\delta} {\mathcal P}(\om)^{V^{P_\alpha}},\\
& \bigcup_{\alpha<\delta} 
\dom\left(\Name{\Xi}^t_\alpha(\Name{B})\right) > 0,\\
& \Name{C}_j \mbox{ comes from (e) or (f)},\\
& \Name{B} \cap \bigcap_{j <j^*} \Name{C}_j = \emptyset.
\end{split}
\end{equation*}
There is some $\alpha(*) < \delta$ such that $\Name{B} \in 
\dom(\Name{\Xi}^t_{\alpha(*)})$ is a $P_{\alpha(*)}$-name. The $\Name{C}_j$ 
have $\n(j) < \n^t$, $\langle \alpha^j_\ell \such \ell \in \om \rangle$,
$\langle \name{p}^j_\ell \such \ell 
\in \om \rangle$ as witnesses as required in
(e) or (f) above. W.l.o.g.\ $q \in P_{\alpha(*)}$ and
$q \in G_{P_{\alpha(*)}} \subseteq P_{\alpha(*)}$,
$G_{P_{\alpha(*)}}$ generic over $V$.

\smallskip
 
We can find $k \in \Name{B}[G_{P_{\alpha(*)}}]$
such that 
$\bigwedge_{j<j^*}\bigwedge_{\ell \in [n^t_k, n^t_{k+1})} (\alpha^j_\ell
> \alpha(*))$ and moreover such that 
$n^t_{k+1} -n^t_k$ is large enough compared to $1/\varepsilon$, $j^*$,
in order to allow us to apply the Tchebyshev inequality and
the law of large numbers for $n^t_{k+1} - n^t_k$ random choices.
(The $n^t_k$ come from item (f) of the definition of a blueprint,
and are not the $\n$.)

\smallskip

Let $\{ \alpha^j_\ell \such j < j^*$ and $\ell 
\in [n^t_k,n^t_{k+1})\}$ be listed as
$\{\beta_m \such m <m^* \}$, in increasing order (so $\beta_0 > \alpha(*)$)
(possibly $\alpha^{j_1}_{\ell_1} = \alpha^{j_2}_{\ell_2} \wedge
(j_1,\ell_1) \neq(j_2,\ell_2)$). We now choose by induction on 
$m\leq m^*$ a condition $q_m \in P_{\beta_m}$ above $q$, 
increasing with $m$ and such that $\dom(q_m)
= \dom(q)\cup \{\beta_0,\beta_1, \dots \beta_{m-1}\}$.
We stipulate $\beta_{m^*} = \delta$.

\smallskip
During this definition we throw a dice
 and the probability of 
success (i.e.\ $q
\Vdash \mbox{``} k \in \name{C}_j$'' for $j < j^*$) is positive, 
and hence $q_{m^*}$ will show that our assumption on $q$
is false.

{\em Case A}: $m=0$
 
 Let $q_0=q$.

{\em Case B}: We are to choose $q_{m+1}$ and 
for some $\n <\n^t$ we have $\n \in \dom(h_0^t)$ and $\gamma$ and:
if $j<j^*$ and $\ell \in \om$ then 
$(\alpha^j_\ell = \beta_m \Rightarrow
\n(j)=\n \; \wedge\; p^j_\ell=\gamma (= h_0^t(\n(j))(\ell)) 
\in Q_{\beta_m})$.

In this case $\dom(q_{m+1}) = \dom(q_m) \cup \{ \beta_m\}$, and
$$
q_{m+1}(\beta) = \left\{ \begin{array}{ll}
  q_m(\beta) & \mbox{ if } \beta < \beta_m, \\
  \gamma &  \mbox{ if } \beta = \beta_m.
\end{array}
\right.
$$

The choice of $(j,\ell)$ is immaterial 
as for each $\beta_m$ there is by the definition of
``satisfying $(t,\n)$ for $\bar{Q}$'', item 5, a unique $\n < \n^t$,
such that there is some  $\ell$ 
such that $\eta_{\beta_m} \restriction w^t = \eta^t_{\n,\ell}$ and
conditions (g) of \ref{4.1} and (d) 8 of \ref{4.2} imply that if
$\eta_{n,\ell}^t $ is not
constant then ($\beta_m = \alpha^{i_1}_{\ell_1}
= \alpha^{i_2}_{\ell_2}
\rightarrow \ell_1=\ell_2$). Hence $\gamma = p_\ell^j$ is well-defined.

{\em Case C}:  We are to choose 
$q_{m+1}$ and for some $\n <\n^t$ we have $\n \in \dom(h_1^t)$ and:
if $j<j^*$ and $\ell \in \om$ then $\alpha^j_\ell = \beta_m \Rightarrow
\n(j)=\n$.

Work first in $V[G_{P_{\beta_m}}]$, $q_m \in G_{P_{\beta_m}}$,
$G_{P_{\beta_m}}$ generic over $V$. The sets
$$
\left\{ \lim(\name{p}^j_\ell[G_{P_{\beta_m}}]) 
\: \mid \: \alpha^j_\ell = \beta_m, \ell \in [n^t_k,n^t_{k+1}), 
j < j^* \mbox{)} \right\}
$$
are subsets of 
$(^\om 2)^{[h^t_2(\n)]} = \{ \eta \in ^\om 2
\such h^t_2(\n) \vartriangleleft \eta \}$. We can define an equivalence
relation $E_m$ on $(^\om 2)^{[h^t_2(\n)]}$:
$$
\nu_1 E_m \nu_2  \mbox{ iff }
\left(
\forall (j,\ell) \mbox{ s.th. } \alpha^j_\ell=\beta_m:
\nu_1 \in \lim(\name{p}^j_\ell[G_{P_{\beta_m}}])
\Leftrightarrow \nu_2 \in \lim(\name{p}^j_\ell[G_{P_{\beta_m}}])\right).
$$

Clearly $E_m$ has finitely many equivalence classes, call them $\langle Z_i^m
\such i < i^*_m \rangle$. All are Borel hence are measurable; 
w.l.o.g.\ ${\rm Leb}
(Z^m_i) = 0 \leftrightarrow i \in [i^\otimes_m, i^*_m)$.
For $i < i^\otimes_m$ there is $r=r_{m,i} \in \name{Q}_{\beta_m}
[G_{P_{\beta_m}}]$ such that 
\begin{eqnarray*}
\lim(\name{p}^j_\ell[G_{P_{\beta_m}}]) \supseteq Z^m_i & \Rightarrow&
r \geq \name{p}_\ell^j[G_{P_{\beta_m}}],\\
\lim(\name{p}^j_\ell[G_{P_{\beta_m}}]) \cap Z^m_i = \emptyset & \Rightarrow&
\lim(r) \cap \name{p}_\ell^j[G_{P_{\beta_m}}]= \emptyset.
\end{eqnarray*}

We can also find a rational $a_{m,i} \in (0,1)_{\mathbb R}$ such that 
$$a_{m,i} < \frac{{\rm Leb}(Z^m_i)}{2^{\lgg(h_2^t(\n))}} <
a_{m,i} + \frac{\varepsilon}{2i^*_m}.$$
We can find $q'_m \in G_{P_{\beta_m}}$, $q_m \leq q'_m$ such that $q'_m$
forces all this information (so for $\name{Z}^m_i$, $\name{r}_{m,i}$
we shall have names, but $a_{m,i}$, $i^\otimes_m$, $i^*_m$ are actual objects.)
We then can find rationals $b_{m,i} \in (a_{m,i}, a_{m,i} + \varepsilon/2)$
such that $\sum_{i < i^\otimes_m} b_{m,i} =1$.

Now we throw a dice
 choosing $i_m < i_m^\otimes$ with the probability
 of $i_m = i$ being $b_{m,i}$, and finally we 
choose $q_{m+1} $ as follows
\begin{eqnarray*}
\dom(q_{m+1}) &=& \dom(q_m) \cup \{\beta_m\}, \\
q_{m+1} &=& \left\{ \begin{array}{ll}
q'_m(\beta) & \mbox{ if } \beta < \beta_m,\\
\name{r}_{m,i_m} & \mbox{ if } \beta = \beta_m.
\end{array}
\right.
\end{eqnarray*}

This covers all cases. 
Basic probability
computation (for $n^t_{k+1} - n^t_k$ independent experiments, using $(*)$
of (f)) 
show that for each $j$ coming from clause (f), 
by the law of large numbers the probability of success,
i.e.\ having $q_{m+1} \Vdash_{P_\delta} k \in \name{C}_j \cap \name{B}$, is
$> (1- 1/{j^*})(1- \eps^{-2}\cdot (n^t_{k+1} -n^t_k)^{-1})$.
For $j$  coming from clause (e) we surely succeed.
\proofendof{\ref{4.5}}

In the following lemma, the whispering conditions (i) of \ref{4.2} are
crucial for building ${\mathcal K}^3$.

\begin{lemma}\label{4.6}
1)  Assume that 
\begin{myrules}
\item[(a)] $\bar{Q} \in {\mathcal K}^3$, $\bar{Q} =
\langle P_\alpha, \name{Q}_\beta, A_\beta, \mu_\beta,\name{\tau}_\beta,
\eta_\beta, (\Name{\Xi}^t_\alpha)_{t \in {\mathcal T}} \such 
\alpha \leq \alpha^*, \beta < \alpha^* \rangle$,
\item[(b)] $A \subseteq \alpha^*$, $\kappa \leq |A|$,
\item[(c)] $\eta \in (^\kappa 2)^V \setminus \{\eta_\beta 
\such \beta \in \alpha \}$,
\item[(d)] $P_A \lessdot P_{\alpha^*}$,
$\name{Q}_{\alpha^*}$ is the $P_{\alpha^*}$-name from 
\ref{2.2} (F)($\beta$) and 
$$ \mbox{ if } t \in {\mathcal T} \mbox{ then } 
\Name{\Xi}^t_{\alpha^*} \restriction 
V^{P_A} \mbox{ is a }P_A\mbox{-name}.$$
\end{myrules}

Then there is $\bar{Q}^+ =
\langle P_\alpha, \name{Q}_\beta, A_\beta, \mu_\beta,\name{\tau}_\beta,
\eta_\beta, (\name{\Xi}^t_\alpha)_{t \in {\mathcal T}} \such 
\alpha \leq \alpha^* +1, \beta < \alpha^* +1
\rangle$ from ${\mathcal K}^3$, extending $\bar{Q}$ such that 
$A_{\alpha^*} = A$, $\eta_{\alpha^*} = \eta$.

\smallskip
2) If clauses (a),(b),(c) of part 1) hold then we can find $A'$ such that
$A \subseteq A' \subseteq \alpha^\ast$, $|A'|\leq (|A| +
$ number of blueprints $)^{\aleph_0}$ such that
$\bar{Q}, A',\eta$ satisfy (a),(b),(c),(d).   

\end{lemma}

\proof 
1) As before the problem is to define $\name{\Xi}^t_{\alpha^*+1}$.
We have to satisfy clause (g) of Definition~\ref{4.2} for each fixed $t 
\in {\mathcal T}$. Let $\n^*$ be the unique $\n < \n^t$ such that $\eta 
\restriction w^t = \eta^t_{\n,\ell}$ for some $\ell \in \om$. 
If $\n^* \in \dom(h^t_0)$ or if
$\langle \eta^t_{\n^*,\ell} \such \ell \in \om \rangle$ is not constant or 
if there is no such $\n^*$ then we have nothing to do.

\smallskip

So assume that $\alpha_\ell = \alpha^*$ for $\ell \in \om$ and that 
$\eta^t_{\n^*,\ell} =  \eta \restriction w^t$ for $\ell \in \om$. Let $\Gamma$ 
be the set of all pairs $(\name{r}, \langle \name{r}_\ell \such \ell 
\in \om \rangle)$ which satisfy the assumption $(**)$ of \ref{4.2}(g). 
In $V^{P_{\alpha^* + 1}}$ we have to choose $\name{\Xi}^t_{\alpha^*+1}$ 
taking care of all these obligations.  

\smallskip

We work in $V^{P_{\alpha^*}}$. By the assumption (d),
which says that $\name{\Xi}^t_{\alpha^*} \restriction P_A$
(hence in particular the  $\name{\Xi}_{\alpha^*}(X)$, 
where $X$ is built from the 
$\name{r}$, $\name{r}_\ell$) is a $P_A$-name, and by Claim~\ref{3.7a}
it suffices to prove it for $\name{\Xi}^t_{\alpha^*+1} \restriction (P_A *
\name{Q})$ (as $\name{\Xi}_1$ there) and for 
$\name{\Xi}^t_{\alpha^*+1} \restriction P_{\alpha^*}$
(as $\name{\Xi}_2$ there) separately, and for the latter there is nothing 
to prove.

By \ref{3.6} it is enough to prove condition (B) of \ref{3.6}.
So suppose that fails. Then there are $\langle B_m \such m < m^* \rangle$, a 
partition of $\om$ from $V^{P_A}$ such that $\Xi^t_{\alpha^*}(B_m) >0$ for
$m < m^*$ and $(\name{r}^i, 
\langle \name{r}^i_\ell \such \ell \in \om \rangle ) \in \Gamma$ and $\n(i) =
\n^* < \n^t$ for $i < i^* < \om$ and $\eps^* >0$, $k^* \in \om$
 and $r \in Q_{\alpha^*} $ which forces that there is no finite
$u \subseteq \om \setminus k^*$ with (i) and (ii) of \ref{3.6}(B).
W.l.o.g.\ $r$ forces that $\name{r}^i \in G_{Q_\alpha}$ for $i < i^*$,
otherwise we ignore such an $\name{r}^i$. So $r \geq r^i$ for $i<i^*$.

By our assumption $(**)$ of \ref{4.2}(g) we have that for each 
$i < i^*$ and $r' \geq r$
$$
\av{\Xi^t_{\alpha^*}}(\langle a^i_k(r') \such k \in \om \rangle) \geq
1 - h^t_1(\n),$$
where
$$ a^i_k(r') = \frac{1}{n^t_{k+1} -n^t_k} \sum_{\ell \in n^t_k, n^t_{k+1})}
\frac{{\rm Leb}(\lim(r') \cap \lim(r^i_\ell))}{{\rm Leb}(\lim(r'))}.
$$

Now $V^{P_A}$ plays the r\^ole of the ground model ($V$ in \ref{3.6}) and
${\rm Random}^{V[\tau_\alpha \such \alpha \in A]} =
{\rm Random}^{V^{P_A}}$ is the full random forcing over this ground model.
So by \ref{3.6} is suffices to prove:

\begin{lemma} \label{4.7}
Assume that $\Xi$ is a finitely additive measure, $\langle B_0, 
\dots B_{m^*-1}\rangle$ a partition of $\om$,
$\Xi(B_m) = a_m$, $i^* < \om$ and $r$, $r^i_\ell \in {\rm Random}$ for $i<
i^*$, $\ell \in \om$ are such that
\begin{myrules}
\item[($\ast$)]
for every $r' \in {\rm Random}$ such that $r' \geq r$ and
for every  $i < i^*$ we have 
$$
{\rm Av}_\Xi(\langle a^i_k(r')
\such k \in \om\rangle) \geq b_i
$$ 
where
$$
a_k^i(r') = 
\frac{1}{n^t_{k+1} -n^t_k} \sum_{\ell=n^t_k}^{n_{k+1}^t-1}
\frac{{\rm Leb}(\lim(r') \cap \lim(r^i_\ell))}{{\rm Leb}(\lim(r'))}.
$$
\end{myrules}

Then for each $\varepsilon >0$, $k^* \in \om$ there is a finite $u \subseteq \om \setminus k^*$ and $r' \geq r$ such that 
\begin{myrules}
\item[{\bf(1)}] $a_m - \varepsilon < |u \cap B_m|/|u| < a_m + \varepsilon$, for
$m < m^*$, 
\item[{\bf(2)}] for each $i<i^*$ we have
$$
\frac{1}{|u|} \sum_{k \in u} \frac{|\{\ell \such n^t_k \leq \ell
< n^t_{k+1} \mbox{ and } r' \geq r^i_\ell\}|}{n^t_{k+1} -n^t_k} \geq b_i -
\varepsilon.
$$
\end{myrules}
\end{lemma}

\proof
Let for $i <i^*$, $m <m^*$:
$$c_{i,m}(r') = \av{{\Xi \restriction B_m}}(\langle a^i_k(r') \such k 
\in B_m \rangle ) \in [0,1]_{\mathbb R}.$$
So clearly
\begin{eqnarray*}
b_i &\leq &\av{\Xi}(\langle a^i_k(r') \such k \in \om \rangle )
= \sum_{m < m^*} \av{{\Xi \restriction B_m}}(\langle a^i_k(r') \such k 
\in B_m \rangle ) \cdot \Xi(B_m) \\
&=& \sum_{m<m^*} c_{i,m}(r') \cdot a_m.
\end{eqnarray*}

Since for each $z \in \om \setminus \{0\}$ there are only finitely many
equivalence classes in the equivalence relation $E_z$ where
$$\langle c_{i,m} \such i<i^*, m<m^* \rangle 
\; E_z \;  \langle c'_{i,m} \such i<i^*, m<m^* \rangle$$
iff $$
(\mbox{ for } z' < z, i < i^*, m<m^*) \; c_{i,m} \in \left.
\left[\frac{z'}{z},
\frac{z'+1}{z}\right.\right) \leftrightarrow c'_{i,m} \in \left.
\left[\frac{z'}{z},
\frac{z'+1}{z}\right.\right),$$ 
we have that there is a condition $r^*_z$ such that
each class is is either dense above $r^*_z$ or does not appear above $r^*_z$. 

We apply this with some $z \geq \frac{1}{\eps}$ and get an $r^* \geq r$ 
and a sequence $\langle c_{i,m} \such i<i^*, m<m^* \rangle$ such that

\begin{myrules}
\item[{\bf(a)}] $c_{i,m} \in [0,1]_{\mathbb R}$,
\item[{\bf(b)}] $\sum_{m < m^*} c_{i,m} \cdot a_m \geq b_i, $\\
\item[{\bf(c)}] for every $r' \geq r^*$ there is $r'' \geq r'$ such that 
$$(\forall i < i^*) (\forall m < m^*) 
[c_{i,m} - \eps < c_{i,m}(r'') < c_{i,m} + \eps].$$
\end{myrules}
Let $k^* \in \om$ be given. We now choose $s^* \in \om$ large enough 
and try to choose by induction on $s \leq s^*$ a condition $r_s
 \in {\rm Random}$ and natural numbers $(m_s,k_s)$ (flipping coins 
along the way) such that:
$$\begin{array}{l}
 r_0 = r^*,\\
r_{s+1} \geq r_s \\
c_{i,m} - \eps < c_{i,m}(r_s)< c_{i,m} + \eps \mbox{ for } i < i^*, m< m^*,\\
k_s > k^*, k_{s+1} > k_s\\
k_s \in B_{m_s}.\end{array}$$

In stage $s$, given $r_s$ we define $r_{s+1}$, $i_s$, $m_s$, $k_s$
as follows: We choose $m_s < m^*$ randomly with the probability
of $m_s$ being $m$ being $a_m$. Next we can find a finite set 
$u_s \subseteq B_{m_s}\setminus \max\{k^*+1, k_{s_1} +1 
\such s_1 < s \}$ such that 

\smallskip

$(+) \; \;$ if $i < i^*$ then $c_{i,m_s} - \eps/2 < \frac{1}{|u_s|} 
\sum_{k \in u_s} a^i_k(r_s) < c_{i,m_s} + \eps/2.$

\smallskip

 We define an equivalence relation ${\bf e}_s$ on $\lim(r_s)$ by

\smallskip

$\eta_1 {\bf e}_s \eta_2$ iff $(\forall i<i^*) (\forall k \in u_s) 
(\forall \ell \in 
[n^t_k, n^t_{k+1}) ) [\eta_1 \in \lim(r^i_\ell) \leftrightarrow
\eta_2 \in \lim(r^i_\ell)].$

\smallskip

The number of equivalence classes is finite. If $Y \in \lim(r_s)/{{\bf e}_s}
$ satisfies ${\rm Leb(Y)} > 0$ choose $r_{s, Y} \in 
{\rm Random}$ such that $\lim(r_{s,Y}) \subseteq Y$. 
Now choose $r_{s+1}$ among $\{ r_{s,Y} \such Y \in \lim(r_s)/{\bf e}_s$ 
and ${\rm Leb}(Y) >0 \}$ with the
probability of $r_{s+1} =r_{s,Y}$ being ${\rm Leb}(Y)$. 
Lastly choose $k_s \in u_s$ with all $k \in u_s$ having the same probability.

\smallskip

Now the expected value (in the probability space of the flipping
coins), assuming that $m_s = m$ of 
$$\frac{1}{n^t_{k+1}-n^t_k} \times |\{ \ell \such 
n^t_k \leq \ell < n^t_{k+1}
\mbox{ and } r_{s+1} \geq r^i_\ell \} |$$ 
belongs to the interval $(c_{i,m} - \eps/2, c_{i,m} + \eps/2)$  
because the expected value of   
$$\frac{1}{|u_s|} \sum_{k \in u_s}
\frac{1}{n^t_{k+1}-n^t_k} \times |\{ \ell \such 
n^t_k \leq \ell < n^t_{k+1}
\mbox{ and } r_{s+1} \geq r^i_\ell \} |$$
belongs to this interval (which is straightforward).

\smallskip

Let $r' = r_{s^*}$, $u=\{k_s \such s \leq s^* \}$. Hence the expected value
of
$$\frac{1}{|u|} \sum_{k \in u}
\frac{1}{n^t_{k+1}-n^t_k} \times |\{ \ell \such 
n^t_k \leq \ell < n^t_{k+1}
\mbox{ and } r' \geq r^i_\ell \} |$$
is $\geq \sum_{m < m^*} a_m (c_{i,m} - \eps/2) \geq b_i - \eps/2$.

\smallskip

As $s^*$ is large enough with high probability (though just positive 
probability suffices), the $(r_{s^*}, \{k_s \such s \leq s^* \} )$ are 
as required for $(r', u)$. Note: We do not know the variance, but we 
have an upper  bound for it not depending on $s$. There is also 
a strong law of large numbers that does not require a bound on the
variance (see \cite{Bauer}).
\proofendof{\ref{4.7}, \ref{4.6}, Part 1)}

\smallskip

Ad \ref{4.6}, Part 2: The proof is an easy counting argument,
just enrich $A$ successively such that everything required becomes
an $P_A$-name. \proofendof{\ref{4.6}, Part 2)}

\medskip \noindent{\bf Remark:}
We do not use \ref{4.6} 2) in our work, nor
do we need here that the number of blueprints is small compared to $\chi$
(which is important in \cite{Shelah592}), because we shall never
use that ${\mathcal K}^3$ is not empty. In \ref{5.3}, \ref{5.4}  we 
need only small parts of the properties of elements in ${\mathcal K}^3$.
So we shall keep the parts needed in mind and, in \ref{5.5}
we shall show that an arbitrary member $\bar{Q}$ of the subclass
of  ${\mathcal K}$ given in \ref{2.2} Part 2)
behaves similarly to a member of ${\mathcal K}^3$ as far as
$(**)_{\bar{Q}}$ is concerned.

\smallskip

The following is needed later to show that sufficiently often the 
 clause (g) of Definition~\ref{4.2} is not trivial, 
that is,  the premise $(**)$ there holds.

\begin{lemma}\label{4.8}
Assume 
\begin{myrules} 
\item[{\bf(a)}] $\Xi$ is a finitely additive measure on 
$\om$ and $b \in (0,1]_{\mathbb R}$,
\item[{\bf(b)}] $n^t_k < \om$ for $k \in \om$, $n^t_k < n^t_{k+1}$, and $
\lim(n^t_{k+1} - n^t_k) = \infty$,
\item[{\bf(c)}] $r^*$, $r_\ell \in {\rm Random}$ are such that: 
$(++) \; (\forall \ell \in \om) [\frac{{\rm Leb}(\lim(r^*)
\cap \lim(r_\ell))}{{\rm Leb}(\lim(r^*))} \geq b].$
\end{myrules}

Then for some $r^\otimes \geq r^*$ we have that
\begin{myrules}
\item[$\otimes(r^\otimes)$] For every $r' \geq r^\otimes$ we have 
$\av{\Xi}(\langle a(r',k) \such k \in \om \rangle ) \geq b$ where:
$a_k(r') = a(r',k) = a_k(\lim(r'))$ and for $X \subseteq 2^\om$ we have that
$$
a_k(X) = 
\frac{1}{n^t_{k+1} -n^t_k}\sum_{\ell \in n^t_k, n^t_{k+1})} \frac{{\rm Leb}(X
\cap \lim(r_\ell))}{{\rm Leb}(X)}.$$
\end{myrules}
\end{lemma}

\proof Let 
$${\mathcal I}= \{ r \in {\rm Random} \such 
r \geq r^*, \mbox{ and } \av{\Xi}(\langle a_k(r') \such k \in \om \rangle ) 
< b \}.$$
If $\mathcal I$ is not dense above $r^*$ there is some $^\otimes \geq r^*$
(in Random) such that for every $r \geq r^\otimes$, 
$r \not\in {\mathcal I}$, so $r^\otimes$ is as required.

So suppose that $\mathcal I$ is dense above $r^*$. We take a maximal antichain
$\{s_i: i \leq i^* \}\subseteq {\mathcal I}$. Because ${\mathcal I}$ is dense 
above $r^*$ we have that $\{s_i: i \leq i^* \}$ is a maximal antichain 
above $r^*$. Hence ${\rm Leb}(\lim(r^*)) =
\sum_{i<i^*} {\rm Leb}(\lim(s_i))$. Since Random has the c.c.c.\ we have that $i^*$ is countable and we assume that $i^* \leq \om$. 

For any $j < i^*$ let $s^j= \bigcup_{i \in j} s_i$. 
Note that $\lim(\bigcup_{m<i} s_m  ) =
\bigcup_{m < i} \lim(s_m)$
and 
$$a_k(s^j) = a_k(\bigcup_{m<i} s_m)
= \sum_{i <j} \frac{{\rm Leb}(s_i)}{{\rm Leb}(\bigcup_{m<j} s_m)} 
a_k(s_i).$$
Hence  we compute
\begin{eqnarray*} 
\av{\Xi}(\langle a_k(s^j) \such k \in \om \rangle ) &=& 
\av{\Xi}(\langle a_k(\bigcup_{m<j } s_m) \such k \in \om \rangle )
\\ &=& \sum_{i<j}  \frac{{\rm Leb}(s_i)}{{\rm Leb}(\bigcup_{m<j} s_m)}
\times
\av{\Xi}(\langle a_k(s_i) \such k \in \om \rangle )
\\ &\leq &
\frac{{\rm Leb}(s_0)}{{\rm Leb}(\bigcup_{i<j} s_i)}(b - \eps)
+ \sum_{0<i<j} \frac{{\rm Leb}(s_i)}{{\rm Leb}(\bigcup_{m<j} s_m)} 
\cdot b\\
&=& b - {\rm Leb}(\lim(s_0)) \cdot \eps,
\end{eqnarray*}
where $\eps = b - \av{\Xi}(\langle a_k(s_0) \such k \in \om \rangle)$, 
so $\eps >0$.

\smallskip

Now let $j$ be large enough such that 
$
\frac{{\rm Leb}(\lim(r^*) \setminus \lim(s^j))}{{\rm Leb}(\lim(r^*))}
< {\rm Leb}(\lim(s_0)) \cdot \eps$.
Then
$$ 
\av{\Xi}(\langle a_k(r^*) \such k \in \om \rangle ) =$$ 
$$
\frac{{\rm Leb}(\lim(r^*) \setminus \lim(s^j))}{{\rm Leb}(\lim(r^*))} \cdot
\av{\Xi}(\langle a_k(\lim(r^*) \setminus \lim(s^j)) \such k \in \om \rangle )
$$
$$ + 
\frac{{\rm Leb}(\lim(s^j))}{{\rm Leb}(\lim(r^*))} \cdot
\av{\Xi}(\langle a_k(\lim(s^j)) \such k \in \om \rangle )
$$
$$\leq 
\frac{{\rm Leb}(\lim(r^*) \setminus \lim(s^j))}{{\rm Leb}(\lim(r^*))} \cdot
1 + 
\frac{{\rm Leb}(\lim(s^j))}{{\rm Leb}(\lim(r^*))} \cdot
(b - {\rm Leb}(\lim(s_0)) \cdot \eps)
$$
$$< {\rm Leb}(\lim(s_0)) \cdot \eps + (b -{\rm Leb}(\lim(s_0)) \cdot \eps) = b
$$
contradicting assumption (c).
\proofendof{\ref{4.8}}

Lemma~\ref{4.5} took care of the successor step in the case
of $|A| \geq \kappa$. We close this section with the successor
step for $|A| < \kappa$ (which means empty $A$ for the
iterations from \ref{2.2} Part 2). Everything in this section
applies to \ref{2.2} Part 1). Only at the end of the next section we shall
make use of the particularly good additional features 
of the narrower class in \ref{2.2} Part 2):
Small forcing conditions, orderly separation between
Cohen part and random part etc.  

\begin{claim}\label{4.9}
 Assume that 
\begin{myrules}
\item[(a)] $\bar{Q} \in {\mathcal K}^3$, $\bar{Q} =
\langle P_\alpha, \name{Q}_\beta, A_\alpha, \mu_\beta,\name{\tau}_\beta,
\eta_\beta, (\Name{\Xi}^t_\alpha)_{t \in {\mathcal T}} \such 
\alpha \leq \alpha^*, \beta < \alpha^* \rangle$,
\item[(b)] $A \subseteq \alpha^*$, $\kappa > |A|$, and $\hat{\mu} < \kappa$,
\item[(c)] $\eta \in (^\kappa 2)^V \setminus \{\eta_\beta 
\such \beta \in \alpha \}$,
\item[(d)] 
$\Name{Q}$ is the $P_{\alpha^*}$-name 
for a forcing notion with set of elements $\hat{\mu}$, and is definable in
$V[\langle \tau_\beta \such \beta \in A\rangle ]$ from
$\langle \tau_\beta \such \beta \in A\rangle$ and parameters from $V$.
\end{myrules}

Then there is $$\bar{Q}^+ =
\langle P_\alpha, \name{Q}_\beta, A_\alpha, \mu_\beta,\name{\tau}_\beta,
\eta_\beta, (\name{\Xi}^t_\alpha)_{t \in {\mathcal T}} 
\such \alpha \leq \alpha^* +1,
\beta < \alpha^* +1
\rangle$$
 from ${\mathcal K}^3$, extending $\bar{Q}$ such that 
$\name{Q}_{\alpha^*} = \name{Q}$, $A_{\alpha^*} = A$, $\eta_{\alpha^*} 
= \eta$, $\mu_{\alpha^*} = \hat{\mu}$.
\end{claim}
\proof The Definition~\ref{4.2} gives no requirements on the 
$\name{\Xi}^t_{\alpha^* +1}$ \proofendof{\ref{4.9}}
 

\section{The Last Part of the Proof of $(**)_{\bar{Q}}$}
\label{S5}

In this section we shall finish the proof of $(**)_{\bar{Q}}$ for
${\mathcal K}^3$, and then we shall finish the proof of \ref{2.11}
and \ref{2.1}.

\smallskip

We give an outline of the proof of $(**)_{\bar{Q}}$ for
${\mathcal K}^3$: We assume that we have a counterexample
$p^*, \Name{T}$ (for a perfect tree $\subseteq
(\bigcap_{\zeta \in \Name{E}} \lim {\rm tree}_m(a^\zeta))^{V[G]}$), 
$m$ (for the ${\rm tree}_{m}$),
 $\Name{E}$ to it. We thin out the $p_\zeta$ that
are forced to be in $\Name{E}$. Thus we get a in some sense indiscernible 
set of conditions. Some features the first $\omega$ of these
indiscernibles are described well by a blueprint $t \in {\mathcal T}$, 
and this description  allows us to
define some $p^\otimes \geq p^*$ such that
$p^\otimes$ forces that $T= \Name{T}[G]$ cannot be a perfect tree
because the subset $A \subseteq \Name{E}[G]$ over which 
we build the intersection is `too large',
and thus we have a contradiction. Having $\Xi^t_\alpha$-measure  non zero
ensures infinity, and indeed the measure $\Xi_\alpha^t$ will 
lead to the notion of `too large' that we
are going use (see \ref{5.2} and \ref{5.3}).

\smallskip

Then we show $(**)_{\bar{Q}}$  for the members of the
subclass of $\mathcal K$ that is given in \ref{2.2} Part 2).
We start looking for  finitely additive measures only after
$p_\zeta$, $\zeta \in \om$ and $t \in {\mathcal T}$ (remember:
${\mathcal T}$ is the set
of blueprints for $\kappa$ from \ref{4.1}) 
are chosen and do it only for
one suitable $t$. We  want to have some
$\Xi^t_{\alpha^*}$ that satisfies just the requirements in \ref{4.2}
(with true premises in (e), (f), (g)
for our chosen $\langle \alpha_\ell \such \ell \in \om \rangle$!)
that speak about our $p'_\zeta$, in order to jump into the proofs
of \ref{5.2} and of \ref{5.3}, which 
work with ${\mathcal K}^3$, and go on like there. 

\smallskip

It turns out that only requirements about $p'_\zeta(\chi + \gamma_n)$, 
$n<n^*\in \om$, $n^*$ the size of the part of the heart of a $\Delta$-system
lying above $\chi$, are relevant. 
We shall look at $\bar{Q}^\chi$ for several $\chi$ (and the same
$\kappa$, $\mu$, $\gamma_0, \dots , \gamma_{n^*-1}$)  
and use a L\"owenheim Skolem argument to
provide the $(\Name{\Xi}^t_{\gamma_n})_{n<n^*, t\in {\mathcal T}}$ 
good for these requirements. Besides some elementary embedding, 
we shall use the automorphisms for
the $\bar{Q}$ from \ref{4.2} Part 2) in order to 
make sufficiently many instances of (e), (g),
(i) of \ref{4.2} true. (We already mentioned that (f) is ad libitum.)

\begin{lemma}\label{5.1}
Suppose that 
$\bar{\eps} = \langle \eps_\ell \such \ell \in \om \rangle$ is a sequence of 
positive reals and that $\bar{Q} \in {\mathcal K}^3$ has
length $\alpha$. Recall that $P'_\alpha$ was defined in \ref{2.3}c).
Then  the following ${\mathcal I}_{\bar{\eps}} \subseteq P_\alpha$ is dense:
\begin{eqnarray*}
{\mathcal I}_{\bar{\eps}} = \{ p \in P'_\alpha &\such 
& \mbox{ there are $m$ and $a_\ell,\nu_\ell$ for $\ell <m$ such that:}\\
& (a) & \dom(p) = \{\alpha_0, \dots \alpha_{m-1}\}, \alpha_0 < \alpha_1 
< \dots < \alpha_{m-1} < \alpha,\\
& (b) & \mbox{if } |Q_{\alpha_\ell}| < \kappa, \mbox{ then } p(\alpha_\ell) 
    \mbox{ is an ordinal},\\
& (c) & \mbox{if } |Q_{\alpha_\ell}| \mbox{ is partial random,  then }
\Vdash_{P_{\alpha_\ell}}\mbox{``} p(\alpha_\ell) \subseteq 
(^\om 2)^{[\nu_\ell]}\\
&& \mbox{\phantom{(c)} and } {\rm Leb}(\lim(p(\alpha_\ell))) 
\geq (1-\eps_\ell)/
2^{\lgg(\nu_\ell)} \mbox{''}\}.
\end{eqnarray*}
\end{lemma}

\proof By induction on $\alpha$ for all possible $\bar{\eps}$. 
Use the Lebesgue Density Theorem \cite{Oxtoby}.
\proofendof{\ref{5.1}}

\begin{lemma}\label{5.2}
If $P_\alpha = \Lim(\bar{Q})$, $\alpha= \lgg(\bar{Q})$ 
and $\bar{Q} \in {\mathcal K}^3$, then $(**)_{\bar{Q}}$ from \ref{2.11} holds.
\end{lemma}

\proof
Suppose that $p^* \Vdash _{P_\alpha} \mbox{``} \Name{T}, m, 
\Name{E}$ form a counterexample to $(**)_{\bar{Q}}$''. Let $\bar{\eps} = 
\langle \eps_\ell \such \ell \in \om \rangle$ be such that $\eps_\ell \in 
(0,1)_{\mathbb R}$ and such that $\sum_{\ell \in \om } 
\sqrt{2\eps_\ell} < 1/10$.
For each $\zeta < \kappa^+$ let $p'_\zeta \geq p_\zeta \geq
 p^*$ be such 
that $p'_\zeta \in {\mathcal I}_{\bar{\eps}}$ is witnessed by 
$\langle \nu^\zeta_\alpha \such \alpha \in \dom(p'_\zeta)
\; \wedge\; |Q_\alpha| \geq \kappa \rangle$ and
$$ 
p'_\zeta \Vdash_{P_\alpha} \mbox{``} \betta_\zeta 
\mbox{ is the $\zeta$-th element such that } \Name{T}
\subseteq N[\name{\bar{a}}^{\betta_\zeta}] \mbox{''.}
$$
Call the $p'_\zeta$ now $p_\zeta$ again. 
By thinning out we may assume that there are
$i^\ast$, $v_0$, $v_1$, $\Delta$, $z$,
$\gamma_i^\zeta$, $\nu_i$, $s^\ast$ such that
\begin{enumerate}
\item $\dom(p_\zeta) = \{ \gamma^\zeta_i \such i < i^* \}$ 
with $\gamma^\zeta_i$ increasing with $i$, let
$v_0^\zeta = \{ i < i^* \such |Q_{\gamma^\zeta_i}|< \kappa \}$, 
then $v_0^\zeta = v_0$ is fixed for all $\zeta$, $v_1 = i^* \setminus v_0$,

\item $\dom(p_\zeta) (\zeta< \kappa^+)$ form a $\Delta$-system with 
heart $\Delta \subseteq \dom(p^*)$,

\item $\betta_\zeta \in \dom(p_\zeta)$, $\betta_\zeta = \gamma^\zeta_z$ 
for a fixed $z < i^*$,

\item $(\dom(p_\zeta), \Delta, \chi, < )$ are isomorphic for $\zeta < \kappa^+$,

\item if $i \in v_0$, then $p_\zeta(\gamma^\zeta_i) = \gamma_i$ for 
$\zeta< \kappa^+$,

\item if $i \in v_1$, then $\nu^\zeta_{\gamma_i^\zeta} = 
\nu_i$ (recall $\nu^\zeta_{\gamma^\zeta_i} \in \,^{<\om} 2$ is given by 
the definition of ${\mathcal I}_{\bar{\eps}}$),

\item $p_\zeta(\betta_\zeta) = s^*$ for $\zeta < \kappa^+$, $s^* = \langle
(n_\ell,a_\ell) \such \ell < m^* \rangle$, w.l.o.g.\
$m^* > m$ (where $m$ is from the counterexample to $(**)_{\bar{Q}}$)
and $m^* > 10$ (this is a similar but not the same as in Lemma~\ref{2.12}),

\item for each $i < i^*$ the sequence $\langle \gamma^\zeta_i 
\such \zeta \in \kappa^+\rangle$ is constant or strictly increasing,

\item the sequence $\langle \betta_\zeta \such \zeta \in \kappa^+ \rangle$ 
is with no repetitions (since, if $p_{\zeta_1}$, $p_{\zeta_2}$ are compatible 
and $\zeta_1 < \zeta_2 < \chi$, then $\betta_{\zeta_1} \neq \betta_{\zeta_2}$).
\end{enumerate}

Now we keep only the first $\omega$ conditions $p_\zeta$, 
$\zeta < \om$. For every such $\zeta$ let $p'_\zeta \geq p_\zeta$ be such 
that $\dom(p'_\zeta) = \dom(p_\zeta)$, $p'_\zeta(\gamma) =
p_\zeta(\gamma)$ except for $\gamma = \betta_\zeta$ in which case 
we extend $p_\zeta(\betta_\zeta) = s^*$ in the following way:

\smallskip

We put $\lgg(p'_\zeta(\betta_\zeta)) = \lgg(s^*) +1 = m^* +1$ 
and set $p'_\zeta(\betta_\zeta) = s^* \concat 
\langle (j^0_\zeta, a_\zeta) \rangle $. 

Before we define $(j_\zeta^0, a_\zeta)$ we choose an increasing sequence 
of integers $\bar{s} = \langle s_\ell \such \ell \in \om \rangle$, 
$s_0 = 0$, such that 
$$s_{k+1} - s_k = |(2^{j_k})^{(2^{j_k}(1- 8^{-m^*}))}|,$$
where
$$j^* = 3n_{m^*-1} +1 $$
(recall from 7.\ that
$n_{m^* -1} $ is  the first coordinate of the last pair in $s^*$) and
we let $j_k = j^* + k!!$ and let $j^0_\zeta = j_k$ when 
$\zeta \in [s_k, s_{k+1})$. Now for $\zeta \in [s_k, s_{k+1})$ define 
$a_\zeta$ such that 
$$ 
\{ a_\zeta \such \zeta \in [s_k, s_{k+1}) \} = 
[^{j_k}2]^{2^{j_k}(1-8^{-m^*})}.$$

For $\eps^*>0$ we define a $P_\alpha$-name by
\begin{equation*}
\Name{A}_{\eps^*} = \left\{ k \in \om \:\left|\: 
\frac{|\{ \zeta  \in [s_k, s_{k+1}) \such p'_\zeta 
\in \Name{G}_{P_\alpha}\}|}
{s_{k+1} -s_k} \right.
> \eps^* \right\}.
\end{equation*}

For the proof of \ref{5.2} we need 

\begin{subclaim}\label{5.3}
There is a condition $p^\otimes\geq p^\ast$ 
that forces that for some $\eps^*>0$ the set 
$\Name{A}_{\eps^*}$ is infinite.
\end{subclaim}

\begin{explanation}
The $p^\otimes$ is an analogue to the premise no.~3 of
Just's Lemma~\ref{2.12}. The condition $p^\otimes(\gamma)$ is
roughly spoken ``as compatible as possible with
many, in the sense of the $\Xi^t_\gamma(A_{\eps^*})> 0$,
of the $\langle
p'_\zeta(\gamma) \such \zeta \in \omega \rangle$''.
The coding with the
$\eta^t_{\n,\zeta}$ and the $\eta_\gamma \restriction w^t$,
$w^t$ from \eqref{code},
ensures that $p^\otimes$ is well-defined by the 
definition below. 
\end{explanation}
 
\proof We may choose any $\eps^* < 1 -
\sum_{\ell \in \omega} \sqrt{2 \eps_\ell}$ 
(where the $\bar{\eps}
=\langle \eps_\ell \such \ell \in \omega
\rangle$ was chosen at the beginning of \ref{5.2}). 
First we define a suitable blueprint
$t \in {\mathcal T}$,
$$t=(w^t, \n^t, \m^t, \bar{\eta}^t, h_0^t, h_1^t, h_2^t, \bar{n}^t).$$
We let 
\begin{equation}\label{code} 
\begin{split}
w^t =& \{ \min\{ \beta \in \kappa \such 
\eta_{\gamma_{i(1)}^{\zeta(1)}}(\beta) \neq 
\eta_{\gamma_{i(2)}^{\zeta(2)}}(\beta) \}
\such \zeta(1), \zeta(2) < \om \mbox{ and }\\
& i(1), i(2) < i^*  \mbox{ and } \gamma_{i(1)}^{\zeta(1)} \neq 
\gamma_{i(2)}^{\zeta(2)} \}, 
\end{split}
\end{equation}
where the $\eta_\alpha$ come from the definition of ${\mathcal K}^3$.
($w^t$ is well-defined because $\eta$ is injective.)

\smallskip

Let $\n^t = i^*$, $\dom(h_0^t) = v_0$, $\dom(h_1^t) = \dom(h_2^t) = v_1$ and 
$n^t_\ell=s_\ell$.

We set $\eta^t_{\n,\zeta} = 
\eta_{\gamma_\n^\zeta} \restriction w^t$. Note that the
$\eta^t_{\n,\zeta}$ satisfy the requirements from \ref{4.1}(g) and (h):
By \ref{5.2} item 4., we have that $\gamma_\n^\zeta = 
\gamma_{\n'}^{\zeta'}$ implies $\n = \n'$. Hence we have that 
$\eta^t_{\n,\zeta} = \eta^t_{\n',\zeta'}$ implies that 
$\eta_{\gamma_\n^\zeta}\restriction w^t =\eta_{\gamma_{\n'}^{\zeta'}}\restriction
 w^t$ and hence by the definition of $w^t$, that $\gamma_\n^\zeta = 
\gamma_{\n'}^{\zeta'}$ and hence $\n = \n'$.

\smallskip

 If $\n \in v_0$, then $h_0^t(\n)(\ell) = \gamma_{\n}$ 
so it is constant independent of $\ell$.
\nothing{ We could have restricted the blueprints
to constant functions $h^t_0$ and then get only $\kappa$ blueprints. See our
recurrent remarks about the numbers of the $p_\zeta$ and the number of
blueprints.}

 If $\n \in v_1$ then $h_1^t(\n) = 
\eps_\n$ and $h_2^t(\n) = \nu_\n$. Finally we set
$\m^t = \max\{k \such \forall \zeta \; \gamma^\zeta_k < \chi \} +1$.

\smallskip

Note that by our choice of $t$,
$\langle \gamma_\n^\zeta \such \zeta \in \om \rangle$ satisfies $(t,\n)$ 
for $\bar{Q}$ for every $\n < i^*$.

We now define a condition $p^\otimes$ such that it will be in $P_\alpha$,
$\dom(p^\otimes) = \Delta$, $p^* \leq p^\otimes$. Remember that 
$\dom(p^*) \subseteq \Delta$, because for each $\zeta$ we have that 
$p^* \leq p_\zeta$. 
If $\gamma \in \Delta$ then for 
some $\n < \n^t$, we have that $\bigwedge_{\zeta \in \om} 
\gamma^\zeta_\n = \gamma$. 

\smallskip
Case: $\n \in v_0$.\\
If $\n \in v_0$ we let $p^\otimes(\gamma) = 
h_0^t(\n)$, so in $V^{P_\gamma}$
$$ p^\otimes \Vdash_{Q_\gamma} \mbox{``} \name{\Xi}^t_{\gamma+1}(\{ \zeta \in
\om \such h_0^t(\n) \in G_{Q_\gamma}\}) = 1 \mbox{ if } \n \in \dom(h_0^t)
\mbox{''.}
$$

\smallskip
Case: $\n \in v_1$.\\
If $\n \in v_1$, then we define a $P_\gamma$-name for a member of 
$Q_\gamma$ as follows. Consider 
$\name{r}_\zeta^\n = \name{p}'_\zeta(\gamma)$ for $\zeta < \om$. 
Let $\name{r} = 
\name{p}^*(\gamma) \cap (^\om 2)^{[h_2^t(\n)]}$  if $\gamma \in 
\dom(p^*)$ and otherwise we let 
$\name{r}$ be just $(^\om 2)^{[h_2^t(\n)]}$. 
Now the premise {\bf(c)} $(++)$ of Lemma \ref{4.8} is true 
with $b = 1- 2\eps_\n$.
Thus by
 Lemma~\ref{4.8} there is some $r^*_\gamma \geq r$ such that 
for every $r' \geq r^*_\gamma$ in $Q_\gamma$ we have that 

\smallskip
\begin{equation}
\tag*{$(**)_{r',\bar{\eps}}$}\label{adhoc} 
\begin{split}
&\av{\Xi^t_\alpha}(\langle a^\n_k(r') \such k \in \om 
\rangle ) \geq 1-2h_1^t(\n) = 1 - 2\eps_\n, \mbox{ where}\\
&a^\n_k(r') = \frac{1}{n^t_{k+1} -n^t_k} \sum_{\ell \in [n_k^t,n^t_{k+1})}
\frac{{\rm Leb}(\lim(r') \cap \lim(r^\n_\ell))}{{\rm Leb}(\lim(r'))}.
\end{split}
\end{equation}

Since $\langle \gamma_\n^\zeta \such \zeta \in \om \rangle $ is constant
since, by \ref{adhoc}
the  assumption $(**)$ of condition (g) of \ref{4.2} holds, 
we get that  in $V^{P_\gamma}$ 

$$
r_\gamma^* \Vdash_{Q_\gamma} \mbox{``} \av{\name{\Xi}^t_{\gamma+1}}
\left(\left\langle
\left.
\frac{|\ell \in [n^t_k, n^t_{k+1}) \such p_\ell(\gamma) 
\in \Name{G}_{Q_\gamma}|}{n^t_{k+1} -n^t_k} 
\:\right|\: 
k \in \om \right\rangle 
\right) 
\geq 1- 2\eps_\n\mbox{''}.
$$
For every $\eps'>0$ we have:
If $\av{\Xi}(\langle a_k \such k \in \om \rangle) \geq
1-\eps'$ then for every $\eps > 0$ such that $\eps + \eps' < 1$, 
\begin{eqnarray*}
\Xi(\{ \ell \such a_\ell \leq 1-\eps' -\eps \}) \cdot (1-\eps'
-\eps) + \Xi(\{ \ell \such a_\ell > 1-\eps' -\eps \}) \cdot 1
&\geq&\\
\av{\Xi}(\langle a_\ell \such \ell \in \om \rangle) &\geq&
1-\eps',
\end{eqnarray*}
and hence
$$
\Xi(\{ \ell \such a_\ell \leq 1-\eps' -\eps \}) 
\leq \frac{\eps'}{\eps'+\eps}.
$$
Now we put $\eps' = 2 \eps_\n$ and get for every $\eps > 0$
$$
r_\gamma^* \Vdash_{Q_\gamma} \mbox{``} \name{\Xi}^t_{\gamma+1}
\left\{k \in \om \:\left|\:
\frac{|\ell \in [n^t_k, n^t_{k+1}) \such p_\ell(\gamma) 
\in \Name{G}_{Q_\gamma}|}{n^t_{k+1} -n^t_k}\right.  
\leq 1- 2\eps_\n - \eps \right\} \leq \frac{2\eps_\n}
{2\eps_\n +\eps}\mbox{''}.
$$

\smallskip
We take $\eps = \sqrt{2\eps_\n} - 2\eps_\n$ 
and thus get
$$
r_\gamma^* \Vdash_{Q_\gamma} \mbox{``} \name{\Xi}^t_{\gamma+1}
\left\{k \in \om \:\left|\:
\frac{|\ell \in [n^t_k, n^t_{k+1}) \such p_\ell(\gamma) 
\in \Name{G}_{Q_\gamma}|}{n^t_{k+1} -n^t_k}  
\right. \leq 1-\sqrt{2\eps_\n}\right\} \leq \sqrt{2\eps_\n}.
\mbox{''}
$$

\smallskip

So there is a $P_\gamma$-name $\name{r}^*_\gamma$ of such a condition. 
In this case let $p^\otimes(\gamma) =
\name{r}_\gamma^*$. So we have finished the definition of 
$p^\otimes$, and it clearly has the right domain. 

\smallskip

[Notice for later generalisation: The property (g) is used
 here only for $\gamma$ in the heart of a $\Delta$-system. Moreover, in order 
to establish (g) for $\gamma$ as in \ref{4.6}, the property (i) is needed
only for $\gamma$.]

\smallskip

Now suppose that $\n< \n^t$ is such that $\gamma_\n^\zeta \not\in \Delta$. 
(Note that this case can be avoided by an appropriate choice of $p'_\zeta$,
see our earlier remarks on simplifications.)
Define $\bar{\beta} =\langle \beta_\zeta \such \zeta \in \om \rangle$, 
$\beta_\zeta = \gamma_\n^\zeta$,
$r_\zeta^\n = p'_\zeta(\gamma_\n^\zeta)$.
Then $\bar{\beta}$ satisfies 
$(t,\n)$ for $P_\alpha$. If $\n \in v_1$, by our assumption that
$p'_\zeta(\gamma) \in {\mathcal I}_{\bar{\eps}}$ and $\eps_\n =
h_1^t(\n)$, we get that the premise of clause (f) of 
\ref{4.2} is fulfilled, hence in $V^{P_\alpha}$:

\smallskip

For each $\eps >0$
$$\Vdash_{P_\alpha} \mbox{``} \name{\Xi}_\alpha^t \left(
\left\{ k \:\left|\: \frac{|\{\ell 
\in [n^t_k,n^t_{k+1}): p_\ell(\gamma_\n^\ell) \in \Name{G}_{\gamma^\ell_\n}\}|}
{n^t_{k+1} -n^t_k} \geq (1-\eps_\n)(1-\eps)\right.
\right\} \right) =1\mbox{''.}
$$
Putting both cases of $\n \in v_1$ (the one with $\gamma^\zeta_\n
\in \Delta$ and the latter, complementary one)
together and assuming that
$p^\otimes \in G$ we get
in $V^{P_\alpha}$ for every $\n \in v_1$:
$$
\sqrt{2\eps_\n} \geq \Xi^t_\alpha \left(\left\{k \in \om \:\left|\: 
1-\sqrt{2\eps_\n} \geq \frac{|\{\ell \such n^t_k \leq \ell < n^t_{k+1} 
\mbox{ and } r_\ell^\n \in G_{P_\alpha}\}|}{n^t_{k+1} -n^t_k} \right.
\right\}  \right).$$

\smallskip

Let $$\Name{A}'_{\eps^*} = \{ k \in \om \such \mbox{ if } \zeta \in 
[n^t_k, n^t_{k+1}) \mbox{ and } i \in v_0 \mbox{ then } 
p_\zeta \restriction \{\gamma_i^\zeta\} \in \Name{G}_{P_\alpha}\}.$$

\smallskip

Then, by \ref{4.2} (e), $\Xi^t_\alpha(A'_{\eps^*}) = 1$.

\smallskip

So
\begin{equation*}\begin{split}
&A_{\eps^*} \cup (\om \setminus A'_{\eps^*})  \supseteq \\
&\left\{k \in \om \:\left|\: \mbox{ if } \n \in v_1 \mbox{ then }
\frac{|\{\ell \such n^t_k \leq \ell < n^t_{k+1} 
\mbox{ and } r_\ell^\n \in G_{P_\alpha}\}|}{n^t_{k+1} -n^t_k}
 \geq 1- \sqrt{2\eps_\n}\right.\right\}=\\
& \om \setminus \bigcup_{\n \in v_1} 
  \left\{ k \in \om \:\left|\:
\frac{|\{\ell \such n^t_k \leq \ell < n^t_{k+1} 
\mbox{ and } r_\ell^\n \in G_{P_\alpha}\}|}{n^t_{k+1} -n^t_k} 
\right. < 1 - \sqrt{2 \eps_\n}\right\}.
\end{split}
\end{equation*}

Hence $\Xi^t_\alpha(A_{\eps^*} \cup (\om \setminus A'_{\eps^*})) 
\geq 1-\sum_{\n \in v_0} \sqrt{2\eps_\n} \geq \eps^* >0$, but
$$\Xi_\alpha^t(\om \setminus A'_{\eps^*}) = 
1 - \Xi^t_\alpha(A'_{\eps^*}) = 1-1=0,$$
hence necessarily $A_{\eps^*}$ is infinite.
\proofendof{\ref{5.3}}

Let $p^\otimes$ be as in the Subclaim~\ref{5.3}. Let $G_{P_\alpha}$ 
be a generic subset of $P_\alpha$ to which $p^\otimes$ belongs. So $A =
\Name{A}_{\eps^*}[G]$ be infinite. For $k \in A$, let
$b_k = \{ \zeta \in [s_k, s_{k+1}) \such p'_\zeta \in G \}$. 
We know that $|b_k|>  (s_{k+1} - s_k) \cdot \eps^*$. Let $\Name{T}[G] = T$.

\smallskip

If $k \in A$, then there are $(s_{k+1}-s_k) \cdot \eps^*$ many $\zeta \in
[s_k,s_{k+1})$ such that $p'_\zeta \in G$ and $p'_\zeta \Vdash
\Name{T} \cap \,^{j_k} 2 \subseteq a_\zeta$, hence 
 $T  \cap \,^{j_k} 2 \subseteq \bigcap_{\zeta \in b_k} 
a_\zeta $ as  $\lgg(s^*) = m^* >m$.
To reach a contradiction  it is enough to show that for infinitely many 
$k \in A$ there is a bound 
on the size of $T \cap \,^{j_k} 2$ which does not depend on $k$.

\smallskip

Now $\frac{|b_k|}{s_{k+1}-s_k}$ is at most the probability that if 
we choose a subset of $^{j_k}2$ with $2^{j_k} (1-8^{-m^*})$ elements, 
it will include $T \cap \,^{j_k}2$. If $k \in A$ (and these are infinitely many $k$, because $A$ is infinite) this probability 
has a lower bound  $\eps^*$ not depending on $k$, and this implies that 
$\langle |T \cap \,^{j_k}2| \such k \in \om \rangle$ 
is bounded and that hence $T$ is finite. 

More formally, for a fixed 
$k \in \om$ we have
\begin{eqnarray*}
|b_k| &=& |\{a_\zeta \such \zeta \in [s_k, s_{k+1}) , \zeta \in b_k \}|\\
&\leq & |\{ a_\zeta \such \zeta \in [s_k, s_{k+1}) , T \cap 
\, ^{j_k}2 \subseteq a_\zeta \}|\\
&\leq & |\{ a \subseteq ^{j_k} 2 \such T \cap\, ^{j_k} 2 \subseteq a 
\mbox{ and } |a|= 2^{j_k}(1-8^{-m^*})\}|\\
&=&|\{ a \subseteq ^{j_k}2 \setminus (T \cap \,^{j_k} 2)\such |a| = 2^{j_k} 
\times 8^{-m^*}\}|\\
&=& \left(\begin{array}{c}
2^{j_k} -|T \cap \, ^{j_k} 2|\\
2^{j_k} \cdot 8^{-m^*}\end{array} \right)
\end{eqnarray*}
By definition we have that $s_{k+1} -s _k = \left(
\begin{array}{c}
2^{j_k}\\
2^{j_k}\cdot (1- 8^{-m^*})\end{array} \right)
= \left(
\begin{array}{c}
2^{j_k}\\
2^{j_k}\cdot 8^{-m^*}\end{array} \right)$. 
Hence
$$
\frac{|b_k|}{s_{k+1}-s_k} \leq \frac{\left(\begin{array}{c}
2^{j_k} -|T \cap \,^{j_k}2|\\
2^{j_k} \cdot 8^{-m^*}\end{array} \right)}{\left(
\begin{array}{c}
2^{j_k}\\2^{j_k} \cdot 8^{-m^*}\end{array} \right)} 
= \prod_{i < |T \cap \,^{j_k} 2|} (1- \frac{2^{j_k} 8^{-m^*}}{2^{j_k} -i}).
$$

Let $i_k(*) = \min(|T \cap \, ^{j_k}2|, 2^{j_k -1})$, so 
\begin{eqnarray*}
\eps^* \leq \frac{|b_k|}{s_{k+1}-s_k} 
&\leq & 
       \prod_{i < |T \cap \, ^{j_k} 2|} \left(1- \frac{2^{j_k} 8^{-m^*}}
        {2^{j_k} -i}\right)\\
&\leq & 
       \prod_{i < i_k(*)} \left(1- \frac{2^{j_k} 8^{-m^*}}{2^{j_k}}\right)
        = (1-8^{-m^*})^{i_k(*)}.
\end{eqnarray*}

So we can find a bound on $i_k(*)$ not depending on $k$:
$$i_k(*) \leq \frac{\log(\eps^*)}{\log(1-8^{-m^*})}.
$$ 

Remember $m^* > 10$, so $1-8^{-m^*} \in (0,1)_{\mathbb R}$. 
So for $k$ large enough,
$$ |T \cap \,^{j_k}2| = i_k(*) \leq \frac{\log(\eps^*)}{\log (1-8^{-m^*})}.
$$
This finishes the proof.\proofendof{\ref{5.2}}

So, how do we get a proof of \ref{2.11} from \ref{5.2}? 
We have to show that our members of $\mathcal K$
as defined in \ref{2.2} Part 2) behave like members of ${\mathcal K}^3$
at sufficiently many points in the domain of the iteration, 
that is we have to define suitable $\name{\Xi}_\alpha^t$ and $\eta$.

\smallskip

Now we shall look at several iteration lengths $\chi$ at the same time.
Recall the definitions of $g_\chi$, $E_\xi^\chi$, $A^\chi_{\chi +\xi}$
from the beginning of the proof of \ref{2.1}.

For $\bar{Q} = \bar{Q}^\chi$ as in \ref{2.2} Part 2)  we set
$\bar{Q}^\chi = P^\chi = P_\chi$ (of length $\chi + \mu$!); 
for $A \subseteq \chi+\mu$, we let
 $P'_A = P'_{\chi,A}$. 

Recall our choice of memories from the beginning
of the proof of \ref{2.1}: $g_\chi \colon \chi \to [\mu]^{<\lambda}$
such that $g_\chi \subseteq g_{\chi'} $ for $\chi < \chi'$ and
such that every point has $\chi$ preimages uner $g_\chi$.
From the $g_\chi$'s we defined:
\begin{eqnarray*}
\mbox{ for } \xi \in \mu \;\; E^\chi_\xi &=&
\{ \alpha < \chi \such 
\xi \not\in g_\chi(\alpha)\},\\
A_{\chi+ \xi}^\chi &=& E_\xi^\chi \cup [\chi,\chi + \xi).
\end{eqnarray*}
We have that $A^\chi_{\chi + \xi} \cap \chi=
A^{\chi'}_{\chi' + \xi} \cap \chi$.

First we need the following

\begin{lemma}\label{5.4}

1. If $\xi \leq \gamma < \mu$ then in $\bar{Q}^\chi$
\begin{myrules}
\item[{\bf (a)}] $P'_{(\chi \cap A_{\chi+\gamma}) \cup [\chi,\chi+ \xi)}
= P'_{E^\chi_\gamma \cup [\chi,\chi + \xi)} \lessdot
P'_{\chi + \xi}$.

\item[{\bf (b)}] if $q  \in P'_{\chi +\xi}$
and $q \restriction (E_\gamma \cup [\chi,\chi + \xi))
\leq p \in P'_{E_\gamma  \cup [\chi,\chi + \xi)}$ then
$$ p \cup q \restriction (\lgg{Q}) \setminus 
(E_\gamma \cup [\chi,\chi+\xi)) \in P'_{\chi +\xi}$$
is the least upper bound of $p$ and $q$.
\end{myrules}

\smallskip
2.  If  $\chi \leq \chi'$, then 
$$
P'_{\chi',\chi \cup [\chi',\chi'+\mu)}  \lessdot P'_{\chi',\chi' + \mu},
$$
and  $P'_{\chi,\chi+\mu}$ is isomorphic to 
$P'_{\chi',\chi \cup [\chi',\chi'+\mu)}$ by
$\hat{h}$ where $h = h^{\chi,\chi'}$ is the canonical mapping, 
i.e.\ $h \colon \chi + \mu  \to \chi'+\mu$ be the identity below
$\chi$ and $h(\chi+\alpha) = \chi'+\alpha$ for $\alpha < \mu$.
\end{lemma}

\proof 1) By \ref{2.6} and \ref{2.7}. For 2): Like in \ref{2.7}, it is
easy to see that $P'_{\chi',\chi \cup [\chi', \chi'+\xi)}
\lessdot  P'_{\chi',\chi' \cup [\chi', \chi'+\xi)}$ as enough types
(see Lemma~\ref{2.7}) are realised in $\chi$.
\proofendof{\ref{5.4}}

\begin{theorem}\label{5.5}
For $\bar{Q}^\chi$ as in Definition~\ref{2.2} Part 2) we have that
$(**)_{\bar{Q}^\chi}$ holds.
\end{theorem}

\proof Given $p^*, \Name{T}, m, \Name{E}$ as in \ref{5.2},
we choose $\bar{\eps}$ and $p'_\zeta$ as in \ref{5.2} (at the end of
\ref{5.2}),
$t$ as in \ref{5.3}.
We let $w_\zeta = \dom(p'_\zeta)$, and $w$ be the heart of the
$\Delta$-system. 
Note that we may choose $p'_\zeta$ such that $w_\zeta \setminus \chi
= w \setminus \chi$, which allows us to avoid \ref{4.2}(f). We now
do so.
We even might choose $p'_\zeta$ such that $w_\zeta \setminus \{\betta_\zeta\}
= w$, but this does not lead to a further simplification.

Let 
$$w \setminus \chi = \{\chi + \gamma_n \such n \in n^* \}, \gamma_0 < 
\gamma_1 < \dots < \gamma_{n^*-1}.
$$

We can replace $\chi$ by $\chi^{+k}$ using $\bar{E}^{\chi^{+k}}$
and thus (by \ref{5.4}) get counterexamples to 
$(**)_{\bar{Q}^{\chi^{+k}}}$ with the same
$t$, $\bar{\eps}$, and with $h^{\chi,\chi^{+k}}(p'_\zeta)$,
$$h^{\chi,\chi^{+k}}{''}(w) \setminus \chi^{+k} = \{\chi^{+k} + \gamma_n 
\such n \in n^* \}, \gamma_0 < 
\gamma_1 < \dots < \gamma_{n^*-1},
$$
and with $A^{\chi^{+k+1}}_{\chi^{+k+1} +\gamma} \cap \chi^{+k}
= A^{\chi^{+k}}_{\chi^{+k} +\gamma} \cap \chi^{+k}$ 
for $\gamma < \mu$.

\smallskip

Now, fixing $\langle \gamma_n \such n <n^*\rangle$ and $\bar{\eps}$, we prove
by  induction on $n < n^*$ that for every $k \in \om$ ($k \leq n^*$
would suffice),
for $\bar{Q}^{\chi^{+k}}$ and for $\gamma_0,\dots, \gamma_{n^*-1}$, 
$\alpha$, and $\langle 
p'_\ell\such \ell \in \om \rangle$ as above, we can find a suitable
modifications $P(n)$ of our original forcing $P^\chi$ and 
$P(n)^{\chi^{+k}}_{\chi^{+k}+\gamma_n +1}$-names for a finitely additive 
measures $(\Name{\Xi}^t_{\chi^{+k} + \gamma_n +1})_{t\in {\mathcal T}}$ 
such that 
\begin{itemize}
\item
 demand (e) of Definition~\ref{4.2} holds for
$\langle \alpha_\ell \such \ell \in \omega \rangle =
\langle f_n \circ \cdots \circ f_0 \circ h^{\chi,\chi^{+k}}
( \gamma_i^\ell) \such \ell \in \omega \rangle$,
$i< i^*$ (from \ref{5.2} 1., only the part before $\chi$
is considered). The $f_i$ are the ``shuffling'' maps coming from the 
L\"owenheim Skolem argument below.
 and such that
\item
(f) and (g) of the
Definition~\ref{4.2} hold for every $n < n^*$ for
$\langle \alpha_\ell \such \ell \in \omega\rangle
= \langle \chi^{+k} + \gamma_n \such \ell \in \omega\rangle $
(so $\alpha_\ell$ is constant) 
and thus  to get the next step in
 the iteration according to \ref{4.6}, and
\item 
though \ref{4.2} (b) is not fulfilled
for $\alpha^* = \chi^{+k} +\mu$, $k \geq 1$,
the original $\eta_\beta \in \,^\kappa 2$ are still
strong enough to code the arguments of
$f_n \circ \cdots \circ f_0 \circ h^{\chi,\chi^{+k}} (p'_\zeta)$, 
$\zeta \in \omega$, according to the \eqref{code} in \ref{5.3}.
Look at the $\gamma_i^\zeta$ to be treated there and at
$f_0, \dots f_{n^*-1}$ and at  $ h^{\chi,\chi^{+k}}$,
how they shift the supports of the $p'_\zeta$.
\end{itemize}

Then we can carry out
the proof of \ref{5.2} and of \ref{5.3}. 
In the end we shall 
first show  $(**)_{P(n^*)^{\chi}}$ for some
modified $P(n^*)^{\chi}$ and mapped $p'_\zeta$,
however with ther same $\mu$, same $\gamma_0,\dots \gamma_{n^*-1}$,
and possibly modified $\betta_\zeta$, $\Name{T}$, $t$.
Thereafter we shall read the automorphisms and 
bijections in the reverse direction
in order to get $(**)_{\bar{Q}^{\chi}}$.

\nothing{
\smallskip
A look at the proof of \ref{5.2} and of \ref{5.3} shows that
\ref{4.2}(e) is relevant only to $\langle \betta_\zeta \such
\zeta \in \omega \rangle$, coming from our choice of $f_\zeta(\betta)
=\betta_\zeta$, $\betta$ as in \ref{2.11}.
\ref{4.2}(f) never occurs, and \ref{4.2}(g) needs to hold only for
$\alpha_\zeta \such \zeta \in \omega \rangle =
\langle \gamma_n^\zeta \such \zeta \in \omega \rangle$, 
$n < n^*$ (from \ref{5.2}),
in order to be applied to
$\langle p'_\zeta(\chi^{+k} +\gamma_n) \such \zeta \in \om \rangle$ 
for $n<n^*$. So we need \ref{4.2}(i) to hold for 
$A^{\chi^{+k}}_{\chi^{+k} + \gamma_n}$
for $n <n^*$ in order to get the next step in
 the iteration according to \ref{4.6}.}
 
\smallskip

In order to proof the claim ``for all $k \in \om$,
$\bar{Q}^{\chi^{+k}}$ can be extended by 
\\
$(\Name{\Xi}_\alpha^{t})_{\alpha \in 
\chi^{+k} + \gamma_n, t\in {\mathcal T}}$ respecting 
the whispering conditions 
at $\chi^{+k} +\gamma_0$, \dots , $\chi^{+k} + \gamma_n$ and
such that $\langle \alpha_\ell \such \ell \in \omega \rangle
= \langle \chi^{+k} + \gamma_n \such \ell \in \omega \rangle$ 
satisfies $(t,\n_n)$
(for the same fixed $t \in {\mathcal T}$, $n < \n^\ast$, with $\n_n =
|\Delta \cap \chi| + n$, 
not depending  on $k$) (let us call this: stage $n+1$)'' ,
we shall use ``for all $k \in \om$,
$\bar{Q}^{\chi^{+k+1}}$ can be extended
by $(\Name{\Xi}_\alpha^{t})_{\alpha \in 
\chi^{+k+1} + \gamma_n}$ respecting the whispering conditions 
at $\chi^{+k+1} +\gamma_0$, $\dots, \chi^{+k+1} + \gamma_{n-1}$ and
such that $\langle \alpha_\ell \such \ell \in \omega \rangle
= \langle \chi^{+k+1} + \gamma_n \such \ell \in 
\omega \rangle $ satisfies $(t,\n_n)$ for $n < n^*$
(let us call this stage $n$)'',
a L\"owenheim and Skolem argument and the uniqueness
of $\n$  in  (d) of Definition~\ref{4.2}.

\smallskip 

To carry out the induction: 
For the stage $n=0$, $k \in \omega$ 
($k= n^*$ would suffice, because we 
need to be able to descend $n^*$ steps in the $k$'s) 
we stipulate that $\gamma_{-1} +1 =0$
and just let $\Name{\Xi}^t_{\chi^{+k}}$ be a $P_{\chi^{+k}}$-name 
for a finitely additive measure
on $\omega$ such that the condition (e) of \ref{4.2}
is fulfilled for the blueprint $t$ and the interesting instances of
$\langle \alpha_\zeta \such \zeta \in \omega \rangle$. 
In the step from stage $n$ to stage $n+1$,
for $\chi^{+k}$, we apply the induction hypothesis to 
$ \gamma_{0}< \cdots <
\gamma_{n-1}$ and $\chi^{+k+1}$ and 
$\langle f_{n-1}^{k+2} \circ \cdots \circ f_0^{k+2+n-1} \circ 
h^{\chi,\chi^{+k+1+n}}(p'_\zeta) \such \zeta \in \om \rangle$, 
(the $f_i^j$ are got from the induction hypothesis, see below, where 
we get $f_n^{k+1}$)
and thus
 we get a $P^{\chi^{+k+1}}_{\chi^{+k+1}+\gamma_{n-1} +1}$-names
$(\Name{\Xi}^t_{\chi^{+k+1} + \gamma_{n-1} +1})_{t\in
{\mathcal T}}$ for   finitely additive measures as required,
i.e.\
the whispering conditions hold for $A^{\chi^{+k+1}}_{\chi^{+k+1} +
\gamma_m}$, $m<n$.

Though we only have $2^\kappa \geq \chi$, the injective coding 
of the indices in the iteration length $\chi + \mu$ by $\eta_{index}
\in 2^\kappa$ works not only for the original $\bar{Q}$ but also 
for $f_{n-1}^{k+2}
\circ \cdots \circ f_0^{k+1+n} \circ
h^{\chi,\chi^{+k+1+n}}{''}(\bar{Q})$, which
is isomorphic to a complete suborder of $ \bar{Q}^{\chi^\kappa}$.

\smallskip

There is a  $P^{\chi^{+k+1}}_{\chi^{+k+1}+\gamma_{n}}$-name 
$\Name{\Xi}^t_{\chi^{+k+1} + \gamma_{n}}$ for a
 finitely  additive measure on $\om$ extending $\Name{\Xi}^t_{\chi^{+k+1}
+ \gamma_{n-1} +1}$: 
This is proved as in 
\ref{4.5} and \ref{4.6}, because 
there are no ``whispering tasks'' (i) of \ref{4.2}
about the $A^{\chi^{+k+1}}_{\chi^{+k+1}+\gamma_n}$ in the stretch
between $\chi^{+k+1} + \gamma_{n-1}+1$ and $\chi^{+k+1} +
\gamma_n$ and no new instances
of (g) of \ref{4.2} as well. 

\smallskip
Now we come to the crucial step from
$\chi^{+k+1} +  \gamma_n$ to $\chi^{+k} +  \gamma_n +1$.
Let 
\begin{equation*}
M_0 \prec M_1 \prec (H(\psi), \in,<_\psi^*),
\end{equation*}
where $\psi= \beth_2(\chi^{+\om})^+$ .

For abbreviation, set $f'=f_{n-1}^{k+2} \circ 
\cdots f_0^{k+2+n-1} \circ h^{\chi,\chi^{+k+n+1}}$, 
and we use $f'$ also for the function
which arises by putting hats over all objects on the right hand side.

\begin{myrules}
\item[$(*)_1$] the objects $\langle \gamma_0, \dots
\gamma_{n^*-1} \rangle$, $\langle g_{\chi^{+l}} \such l \in \om 
\rangle$, $\langle h^{\chi^{+k},\chi^{+k+1}} \such k \in \omega 
\rangle$,\\
 $\mu$, $\chi$,
$\langle f'(p'_\zeta) \such \zeta < \om \rangle$, 
$\langle \bar{Q}^{\chi^{+k}} \such k \in \om \rangle$,
$\langle \wittt{P}_{n-1}^{\chi^{+k}} \such k \in \om \rangle$
$(\Name{\Xi}^t_{\chi^{+k+1} + \gamma_n})_{t\in {\mathcal T}}$,
\\$f'(\Name{T}) = {\mathcal B}(\langle \mbox{truth value}(f'(\delta_\ell)
\in \name{\tau}_{f'(\gamma_\ell)}) \such \ell \in \om \rangle )$
 belong to $M_0$.

\item[$(*)_2$]  $\Vert M_0 \Vert = \Vert M_1 \Vert =
\chi^{+k}$, $\chi^{+k} +1 \subseteq M_0$, $M_0 \in M_1$,
$^{\max(\mu,\kappa)}(M_0) \subseteq M_0$,
$^{\max(\mu,\kappa)}(M_1) \subseteq M_1$.
\end{myrules}

Claim: There is an injective
function $f^{k+1}_n$ from $(\chi^{+k+1} + \gamma_n +1) \cap M_1$ to
$\chi^{+k} + \gamma_n+1$ such that
\begin{myrules}
\item[{\bf (a)}] $f^{k+1}_n(\chi^{+k+1} +\gamma) 
= \chi^{+k} + \gamma$ for $\gamma
                 \leq \gamma_n$,

\item[{\bf (b)}] $f^{k+1}_n$ maps $(\chi^{+k+1}+\gamma_n)
 \cap M_0$ onto $A^{\chi^{+k}}_{\chi^{+k} +
\gamma_n}$ and

\item[{\bf (c)}]  $g_{\chi^{+k}}(f^{k+1}_n(\alpha)) \cap \gamma_n 
                 =g_{\chi^{+k+1}}(\alpha) \cap \gamma_n$ for $\alpha 
                   \in \lambda^{+k+1} \cap M_1$, i.e.\
                  for $\gamma \in \gamma_n,  \alpha 
                   \in \lambda^{+k+1} \cap M_1$: $(f^{k+1}_n(\alpha) \not\in 
                  A^{\chi^{+k}}_{\chi^{+k} + \gamma} \leftrightarrow
                  \alpha \not\in 
                  A^{\chi^{+k+1}}_{\chi^{+k+1} + \gamma})$.

%

\end{myrules}


Proof of the claim: 
\nothing{
\begin{eqnarray*}{\mathcal H}(\psi), M_0, M_1 &\models&
\name{\Xi}^1 \mbox{ is a name for a finitely additive measure on } \om\\
&& \mbox{ such that for $n^* \geq n > m$ its restriction to } 
{\mathcal P}(\om) \mbox{ in }
V^{P^{\chi^{+k+1}}_{A^{\chi^{+k+1}}_{\chi^{+k+1} +\gamma_m}}}\\
&&\mbox{ is a }P^{\chi^{+k+1}}_{A^{\chi^{+k+1}}_{\chi^{+k+1} +\gamma_m}}
\mbox{-name} \mbox{ and }\\
&& \mbox{(e) of \ref{4.2} holds, } \\
&&\langle \alpha_\ell \such \ell \in \om \rangle =
\langle \chi^{+k+1} + \gamma_n \such \ell \in \omega \rangle
\mbox{ satisfies } (t,\n_n) \mbox{ for } \bar{Q} \mbox{ for } n < n^*,\\
&&\mbox{(g) of \ref{4.2} holds for every $m < n$ for }\\
&&\langle h^{\chi,\chi^{+k+1}}(p'_\zeta)(\chi^{+k+1} 
+\gamma_m) \such \zeta \in \om \rangle 
 &&\langle p'_\zeta \such \zeta \in \om \rangle, t, \Name{T},
 m \mbox{ form 
 a counterexample to } (**)_{\bar{Q}} \mbox{ in the }\\
 && \mbox{way described
 in the part of the proof of \ref{5.2} before \ref{5.3}}
\\ &&= \phi(\chi^{+k+1}).
\end{eqnarray*}
} 
Since $M_0 \in M_1$ we have that
 $|\chi^{+k+1} \cap (M_1 \setminus M_0)|=
 |\chi^{+k+1} \cap M_1| = |\chi^{+k+1} \cap M_0|$, and considering
types as in the proof of \ref{2.7} we get
for any $c \in ^{n+1} 2$, with $E^{0} = E$, $E^1 = \chi^{+k} \setminus E$,
\begin{eqnarray*}
\left|
M_1 \cap \bigcap_{m<n+1} 
(E^{\chi^{+k+1}}_{\gamma_m})^{c(m)}\right| & = &\chi^{+k},\mbox{ and }\\
\
\left|
\bigcap_{m<n+1} 
(E^{\chi^{+k}}_{\gamma_m})^{c(m)}\right| & = &\chi^{+k},\mbox{ and }\\
\left|
M_0 \cap \bigcap_{m<n+1} 
(E^{\chi^{+k}}_{\gamma_m})^{c(m)}\right| & = &\chi^{+k}, \mbox{ and }\\
\left|
M_0 \cap \chi^{+k+1} \right| & = & \chi^{+k}.
\end{eqnarray*}
Hence we can find an $f_n^{k+1}$ fulfilling the requirements 
{\bf (a)}, {\bf (b)}, and {\bf (c)}.
\nothing{
bijection $f_n$ such that
\begin{eqnarray*}
f_n &\colon & (\chi^{+k+1} + \gamma_n +1) \cap M_1  
\to \chi^{+k} + \gamma_n +1,
\\
f_n(\chi^{+k+1} +\gamma) & = & \chi^{+k} + \gamma \mbox{ for } \gamma
                 \leq \gamma_n\\
f_n{}^{''} ((\chi^{+k+1} + \gamma_n)\cap M_0) & =& 
A^{\chi^{+k}}_{\chi^{+\kappa} + \gamma_n} \\
&(= & \{ \alpha < \chi^{+k} \such \gamma_n \not\in g_{\chi^{+k}}(\alpha)\}
\cup [ \chi^{+k}, \chi^{+k} + \gamma_n), \\
&& \makebox[3cm]{} \mbox{ see } \ref{2.1})\\ 
\\
g_{\chi^{+k}}(f_n(\alpha)) \cap \gamma_n
                 & = & g_{\chi^{+k+1}}(\alpha) \cap \gamma_n
\mbox{ for } \alpha 
                   \in \lambda^{+k+1} \cap M_1,\\
 \mbox{and hence }
f^{''} (A^{\chi^{+k+1}}_{\chi^{+k+1} +\gamma_m}) & =& 
A^{\chi^{+k}}_{\chi^{+\kappa} +  \gamma_m} \mbox{ for } m <n.
\end{eqnarray*}}
Hence the claim is proved.

\medskip 
Now we change the forcing orders accordingly:
We set  $\wittt{P(0)}^{\chi^{+k}} = P^{\chi^{+k}}$.
As in \ref{2.5} we can define a structure 
$\wittt{P(n)}^{\chi^{+k}}$
by 
$$
\widehat{f_n^{k+1}}\colon  (\wittt{P(n-1)}^{\chi^{+k+1}})\cap M_1 \cong 
\wittt{P(n)}^{\chi^{+k}}
$$
and can extend $\widehat{f_n^{k+1}}$ onto the space of 
$(\wittt{P(n-1)}^{\chi^{+k+1}})\cap M_1$-names.

\smallskip

From $f_n^{k+1} \circ f' \circ h^{\chi,\chi^{+k+n+1}}(\chi + \gamma_m)) = 
\chi^{+k} + \gamma_m$
we get that  $\langle \alpha_\ell \such \ell \in \omega \rangle
= \langle \chi^{+k} + \gamma_m \such \ell \in \omega
\rangle$ still satisfies $(t,\n_m)$
(see \ref{4.2}(d))
for $\wittt{P(n)}^{\chi^{+k}}$ for every $m \leq n^*$. 
Moreover, $f_n^{k+1} \circ f' \circ  h^{\chi,\chi^{+k+n+1}}
(\chi + \gamma_n)
= \chi^{+k} + \gamma_n$ is the argument where
$\langle f_n^{k+1} \circ f' \circ
 h^{\chi,\chi^{+k+n+1}}(p'_\zeta) \such \zeta \in \omega
\rangle$ is treated as in \ref{5.2}.

\smallskip

Now we prove that $\wittt{P(n)}^{\chi^{+k}}$ satisfies the conditions at
$\gamma_0, \gamma_1 \dots \gamma_n$: 

First, for $m =n$,
we have that $\Name{\Xi}^t_{\chi^{k+1} +
\gamma_n}$ is in $M_1$ a 
$\wittt{P(n-1)}^{\chi^{+k+1}}\cap M_1$-name, and 
its restriction to $ 
{\mathcal P}(\om)^{V^{\wittt{P(n-1)}^{\chi^{+k+1}}
_{\chi^{+k+1} +\gamma_n}}}\cap M_i$
is a
$\wittt{P(n-1)}^{\chi^{+k+1}}_{\chi^{+k+1} +\gamma_n}\cap M_i$-name.
We get that
$\widehat{f_n^{k+1}}(\Name{\Xi}^t_{\chi^{k+1} +
\gamma_n}) \restriction M_1=:\Name{\Xi}^t_{\chi^{k} +
\gamma_n +1} (= \Name{\Xi}^t$ in the next paragraphs)  is as required:
We write only $f$ for $f^{k+1}_n$ in the proof of this claim
so that the notation be slightly less clumsy.
\smallskip

We show that it is a 
$\wittt{P(n)}_{\chi^{+k}+ \gamma_n+1}^{\chi^{+k}}$-name 
for a  finitely additive measure on $\om$
such that its restriction to ${\mathcal P}(\om)$
in $V^{\wittt{P(n)}^{\chi^{+k}}
_{A^{\chi^{+k}}_{\chi^{+k} +\gamma_n}}}$ is a 
$\wittt{P(n)}^{\chi^{+k}}_{A^{\chi^{+k}}_{\chi^{+k} +\gamma_n}}$-name, 
so condition (i)  of \ref{4.2} is satisfied:
Let $\Name{A}$ be a 
$\wittt{P(n)}^{\chi^{+k}}_{A^{\chi^{+k}}_{\chi^{+k} +\gamma_n}}$-name:

$$\hat{f}(\Name{\Xi}^t)(\Name{A}) =\hat{f}(\Name{\Xi}^t)(\hat{f_n}
(\hat{f}^{-1}(\Name{A}))),
$$
where $\hat{f}^{-1}(\Name{A}) \in M_0$.

\smallskip

Hence
$$\hat{f}(\Name{\Xi}^t)(\hat{f_n}(\hat{f}^{-1}(\Name{A}))) = \hat{f}
(\Name{\Xi}^t(\hat{f}^{-1}(\Name{A})))$$
and
where $\Name{\Xi}^t(\hat{f}^{-1}(\Name{A})) \in M_0$. 
Hence $\hat{f}(\Name{\Xi}^t (\hat{f}^{-1}(\Name{A})))$ is an 
$f''(M_0 \cap (\chi^{+k+1} + \gamma_n)) = 
A^{\chi^{+k}}_{\chi^{+k} + \gamma_n}$-name.

For $m < n$ the claim that $\Name{\Xi}^t_{\chi^{k} +
\gamma_m  +1} := \Name{\Xi}^t_{\chi^{k} +
\gamma_n  +1} \restriction ({\mathcal P}(\om)
\mbox{ in } V^{\wittt{P(n)}^{\chi^{+k}}
_{\chi^{+k} +\gamma_m+1}})
$ is a  
$\wittt{P(n)}_{\chi^{+k}+ \gamma_m +1}^{\chi^{+k}}$-name for a 
finitely additive measure on $\om$
such that its restriction to ${\mathcal P}(\om)$
in $V^{\wittt{P(n)}^{\chi^{+k}}
_{A^{\chi^{+k}}_{\chi^{+k} +\gamma_m}}}$ is a 
$\wittt{P(n)}^{\chi^{+k}}_{A^{\chi^{+k}}_{\chi^{+k} +\gamma_m}}$-name, 
follows from
\\
$f_n^{''} (A^{\chi^{+k+1}}_{\chi^{+k+1} +\gamma_m})  =
A^{\chi^{+k}}_{\chi^{+\kappa} +  \gamma_m} \mbox{ for } m <n$.

\smallskip

Hence we have
 $\Name{\Xi}^t_{\chi^{+k+1} + \gamma_m +1}$, which are
$\wittt{P(n)}^{\chi^{+k}}_{\chi^{+k} +\gamma_m +1}$-names respecting
the whispering conditions \ref{4.2}(i) 
at $\chi^{+k} + \gamma_0, \dots , \chi^{+k} + \gamma_n$
(which where needed in 
the premises of \ref{4.6} 1)), and the
inductive proof is finished.

\medskip

Now we perform the induction with starting point $h^{\chi,\chi^{+n^*}}(P)$
and get $f_0^{n^*}$, $f_1^{n^*-1}$, \dots ,$f_n^{n^*-n}$, \dots,
$ f^1_{n^*-1}$
and $k := f^1_{n^*-1} \circ \cdots \circ f^{n^*}_0$, $f:= k \circ
h^{\chi,\chi^{+n^*}}$.
After $n^*$ induction steps, we have that the mapped forcing 
$\hat{k}''P^{\chi^{+ n^*}} = P(n^*)^\chi$
is expanded by measures $\Xi^t_{\chi + \gamma_n +1}$, $n \leq n^*$.

\smallskip

So the proofs of \ref{5.2} and \ref{5.3} go through for the
modified forcing and the mapped 
objects: $\hat{f}(\Name{T})$,
$\hat{f}(p'_\zeta)$, $\hat{f}(t)$ (blueprints),
 $\langle \hat{f}(\gamma_i^\zeta) \such i < i^*\rangle$ 
(the domain of $\hat{f}(p'_\zeta)$). Hence the proofs of
\ref{5.2} and of \ref{5.3} show that there is no perfect tree in the
intersection of the mapped trees. So $\hat{f}(\Name{T})$ is
not perfect in the generic extension $V^{\wittt{P(n^*)}^\chi}$.

\smallskip
We have that $h^{\chi,\chi^{+n^*}}$ is a complete embedding,
and that in each step $P(n)^{\chi^k}$ is isomoprhic to
$\wittt{P(n-1)}^{\chi^{+k+1}}\cap M_1$, which is
is a complete suborder of\\
  $\wittt{P(n-1)}^{\chi^{+k+1}}$
(because $M_1 \prec H_\psi$ and all antichains are countable and $^\om M_1
\subseteq M_1$.) Being a perfect tree is absolute
for $\zfc$ models and hence $n^* +1$ applications
of \cite[VII Lemma 13]{Kunen} the condition
$$
p \Vdash_{\wittt{P(n^*)}^\chi}
\mbox{``}\hat{f}(\Name{T}) \mbox{ is 
not perfect in the generic extension } V^{\wittt{P(n^*)}^\chi} \mbox{''}
$$
implies  that some condition in $G$ forces that
$\Name{T}$ is not a perfect tree in $V^P$.
Thus \ref{star} is also proved for the original $\bar{Q}$.

\nothing{
We have $(**)_{\bar{Q}''}$ for some $\bar{Q}''$, got after the
$n^*$  induction steps, for some $p''_\zeta$, $\alpha''_\zeta$,
still satisfying $(t_n,\n_n)$ for $\bar{Q''}$ for $n < n^*$,
and some moved $\Name{T}''$ and moved $\Name{E}''$, such that
 $\bar{Q}''$ is expandible to 
 a member of ${\mathcal K}^3$ at least at the $\gamma_m$
 in the heart of the $\Delta$-system for the domains of the
 (modified) $p'_\zeta$ lying above $\chi$. 
\smallskip
Since all the shufflings $f$ used in the induction leed via
$\hat{f}$ to
isomorphic versions of the initial counterxample, 
we can transfer the situation back to $p_\zeta$ and $\bar{Q}$
and get that there is no perfect tree
in the original intersection and thus that $(**)_{\bar{Q}}$ holds.}
\proofendof{\ref{5.5},\ref{2.11},\ref{2.1}}


\section{The Case of $\cf(\mu) = \om$}\label{S6}

In this section we show a version of Theorem~\ref{2.1}
for the case of $\cf(\mu) = \om$. The main technical point is:
The part of the iteration as in \ref{2.2} Part 2)
lying before $\chi$ and the part thereafter
now are going to take shifts $\omm$ often.

This means a slight increase of the complexity of our notation.
We are going to rework the previous three sections and benefit
from the fact that we did some (but not all) work for the class of forcings 
of \ref{2.2} Part 1). We shall often only hint
to the parallels and give an informal description of the modifications
and strengthenings.

\begin{theorem}\label{6.1}
In \ref{2.1}, we can replace ($\cf(\mu) > \aleph_0$ and $\sup(C) = \mu$) by
\begin{eqnarray*}
\cf^{V_1}(\mu) &=& \om , \mbox{ and there is some $\lambda$ such that }\\
\omm \leq |C|^{V_2} & <& \lambda < \mu, \mbox{ and }\\
\cf^{V_1}(\lambda) & \geq &\omm, \mbox{ and }\\
\forall B \in V_1 && (|B|^{V_1} < \lambda \rightarrow C \not\subseteq B).
\end{eqnarray*}
\end{theorem}

\proof 
We first give an outline:
We define a member of $\mathcal K$ (of \ref{2.2}) that we are going to
use. Then (after adapting \ref{2.6} and
\ref{2.7}) we get the items $(\alpha)$ to $(\delta)$ 
of the conclusion of  \ref{2.1} and of \ref{6.1}.
For item $(\eps)$, we begin with the analogon of the end of
Section~\ref{S2}. Then we slightly modify the blueprints.
 Again we can deal with automorphisms of the iteration length.
We take those automorphisms moving only
some element $\alpha$ within one
 of our $\omm$ intervals $[\chi \cdot \gamma, \chi \cdot (\gamma +1))$. 
So we basically do the old proof in some
interval of the longer iteration. We use that we never required that
there are only partial random forcings after $\chi$.

\smallskip

We take $\chi \geq 2^\mu$ and $\kappa$ such that $2^\kappa \geq \chi$.
Then we define
$$
\bar{Q}^\chi = 
\langle P^\chi_\alpha, \name{Q}_\beta, A^\chi_\beta, \mu_\beta,
\name{\tau}_\beta \such \beta < \chi \cdot \omm,
\alpha \leq \chi \cdot \omm \rangle
\in {\mathcal K} $$
as follows: 

\medskip

We take for $\chi \leq \chi'$
\begin{equation*}
\begin{split}
&g_{\chi,\omm} \colon  \chi \cdot \omm \to (\mu\times \omm)^{<\lambda}\\
&g_{\chi,\omm}(\chi \gamma + \xi) = \emptyset \mbox{ for }  \mu \leq 
\xi < \chi,\\
&g_{\chi',\omm}(\chi' \gamma + \xi) = g_{\chi,\omm}(\chi \gamma + \xi) 
\mbox{ for } \xi < \chi, \gamma \in \omm\\
&\forall \gamma \in \omm \;\, \forall  B \in (\mu \times \omm)^{<\lambda}
\;\; 
\exists^{\chi'} \alpha \in [\chi' \cdot \gamma,\chi'\cdot(\gamma +1)) \;\;
 g_{\chi',\omm} (\alpha) = B.
\end{split}
\end{equation*}

For $\alpha = \chi \cdot \gamma + \xi$, $\gamma \in \omega_1$,
$\xi \in \chi$ we set
$$
A^\chi_\alpha = \left\{ \begin{array}{ll}
\emptyset & \mbox{ if } \gamma = 0 \mbox{ or } \xi > \mu,\\
\{\beta <\chi \cdot \gamma \such 
 (\xi,\gamma) \not\in g_{\chi,\omm}(\beta)\} &
\mbox{ else }.\end{array}
\right.
$$ 
$$\name{Q}_\alpha = 
\left\{ \begin{array}{ll}
(^\om 2, \vartriangleleft), & \mbox{ if } A^\chi_\alpha = \emptyset,\\
{\rm Random}^{V[\name{\tau}_\beta \such \beta \in A^\chi_\alpha]}, &
\mbox{ else.}
\end{array}
\right.
$$

We adopt \ref{2.4} as follows

\begin{definition}\label{6.2}
For $\bar{Q} \in {\mathcal K}$ of the special form of
\ref{6.1}, $\alpha < \chi \cdot \omm$, we let
\begin{eqnarray*}
AUT(\bar{Q}^\chi\restriction \alpha) & = &
\bigl\{ f \colon \alpha \to \alpha \such f \mbox{ is bijective, and },
\\
&&
(\forall \beta , \delta \in \alpha)
\\&& \bigl((|\name{Q}_\alpha|<\kappa \leftrightarrow 
|\name{Q}_{f(\alpha)}|
< \kappa) \wedge \\ 
&&
(\beta \in A_\delta  \leftrightarrow f(\beta) \in A_{f(\delta)})\bigr)\bigr\}.
\end{eqnarray*}
\end{definition}

Then we have that $\hat{f}$ is an automorphisms of $P_\alpha$ and of 
$P'_\alpha$ (from Definition~\ref{3.2} (c)), and
Fact \ref{2.5} holds for ${\mathcal K}$.

\smallskip

Now we get the analogues of \ref{2.6} and of \ref{2.7}
(consider types, similarly to there)
and are ready to prove
 
\begin{myrules}
\item[($\delta'$)] $V_2 \models \; \Vdash_{P_{\chi\cdot\omm}}\;
\mbox{``}\{ \tau_{\chi \cdot \gamma +i} \such i \in C, \gamma \in \omm
 \}$ is not null.''
\end{myrules}

\proof
Let $\Name{N}$ be a $P_{\chi \cdot \omm}$-name for a Borel null set. Hence
for some Borel function ${\mathcal B} \in V_1$ and for some countable
$$
X = \{ x_\ell \such \ell \in \omega \}\subseteq \chi \cup
\bigcup_{\gamma
\in \omm \setminus\{0\}} [\chi \cdot \gamma + \mu,
\chi \cdot (\gamma +1)),$$
$$
Y=\{ y_\ell \such \ell \in \omega \} 
\subseteq \bigcup_{\gamma \in \omm \setminus\{0\}}
[\chi \cdot\gamma , \chi\cdot \gamma +\mu),
$$
$\zeta_\ell$, $\ell \in \omega$,
$\zeta'_\ell$, $\ell \in \omega$, we have that
$$\Name{N}=
{\mathcal B}((
\mbox{truth value}(\zeta_\ell \in 
\name{\tau}_{x_\ell}))_{\ell \in \omega},
(\mbox{truth value}(\zeta'_\ell \in 
\name{\tau}_{y_\ell}))_{\ell \in \omega}).$$

Let $i(\ast) < \omm$ be such that $\chi \cdot i(\ast)> \sup(Y)$. 
Since $\cf^{V_1}(\lambda) > \aleph_0$, we have that
$B:= \bigcup_{\xi \in X \cup Y} g_{\chi,\omm}(\xi) \in
 ([\mu \times \omm]^{<\lambda})^{V_1}$.

\smallskip

Since $C \setminus \pi_{\mu}(B) \neq \emptyset$, there is some $i \in \mu$,
$i \in C \setminus \pi_{\mu}(B)$. 
We claim, that $\tau_{\chi \cdot i(\ast)+ i}$ 
is random over 
a universe, in which $\Name{N}[G]$ has a name. (Moreover
regarding $V_1$ and $V_2$, the same remarks as
in the proof of $(\delta')$ of Theorem~\ref{2.1} apply.)
Then the proof will be finished, because then 
$\tau_{\chi \cdot i(\ast) + i} \not\in \Name{N}[G]$ in $V_2[G]$.
By our construction, we have
$$
\tau_{\chi \cdot i(\ast)+i} \mbox{ is the } {\rm Random}^{V[\name{\tau}_\alpha 
\such \alpha \in A^\chi_{\chi \cdot i(\ast) + i}]}
\mbox{-generic over } V^{P_{\chi \cdot i(\ast)+ i}}.
$$
Since $i \in C \setminus \pi_{\mu}(B)$, we have that 
$\forall \xi \in X \cup Y$ $\forall \gamma \in \omm$
that $g_\chi(\xi) \not\ni (i,\gamma)$, hence
$\forall \xi \in X \cup Y \: \xi \in A_{\chi \cdot \gamma + i}$, so $X \cup Y
 \subseteq A^\chi_{\chi\cdot i(\ast) + i}$. 
Since $P_{A_{\chi \cdot i(\ast) + i}} \lessdot P_{\lgg(\bar{Q})}$ 
the name $\Name{N}$ is
evaluated in the right manner in $V^{P_{A_{\chi \cdot i(\ast) + i}}}$.
Thus the claim is proved. \proofendof{(\delta')}

\begin{myrules}
\item[($\delta$)] $V_2[G] \models \unif{\lebesgue} \leq |C|$.

\end{myrules}
This follows from ($\delta'$).

\smallskip

\begin{myrules}
\item[($\eps$)] $V_1[G] \models \unif{\lebesgue} \geq \lambda$.
\end{myrules}

Again the item $(\eps)$ will be the longest part.
However, it is almost the same as our previous work. Put 
all the $\betta_\zeta$ of an analogue of \ref{2.11}
into one $[\chi \cdot \gamma + \mu, \chi \cdot
(\gamma +1))$. Also the extension of $\chi$ to $\chi'$
now can be done either only in the relevant interval where the
$\alpha_\zeta$ lie, or just all over, thus leading to $h^{\chi,\chi'}$.

\smallskip

More explicit, we start as in the corresponding proof in \ref{2.1}:
 Suppose that $(\eps)$ is  not true. In $V_1$
there is $i(\ast) < \lambda$ and $p \in P_{\chi\cdot \omm}$ such that
$$
p \Vdash_{P_{\chi \cdot \omm}} 
\mbox{``}\name{\eta}_i \in \; ^\omega 2 \mbox{ for } i < i(\ast)
\wedge \{ \name{\eta}_i \such i < i(\ast) \} \mbox{ is not null.''}
$$
A name of  a real in $V_1[G]$ is given  by
$$
\name{\eta}_i = 
{\mathcal B}_i(\langle \mbox{truth value}(\zeta_{i,\ell} \in 
\name{r}_{j_{i,\ell}}) \such \ell \in \omega \rangle)$$ for
suitable $\langle \zeta_{i,\ell}, j_{i,\ell} \such \ell \in 
\omega \rangle$,
$\zeta_{i,\ell} \in \omega$, $j_{i,\ell} \in \chi + \mu$.

We set 
\begin{eqnarray*}
X &=& \{ j_{i,\ell} \such i \in i(\ast), 
\ell \in \omega \} \cap (\chi \cup \bigcup \{ \chi \cdot \gamma + \mu,
\chi \cdot(\gamma +1)) \such \gamma \in \omm \setminus \{0\} \},\\
Y &=& \{j_{i,\ell} \such i \in i(\ast), \ell \in \omega \} \cap 
\bigcup \{ \chi \cdot \gamma,
\chi \cdot\gamma +\mu) \such \gamma \in \omm \setminus \{0\} \}
\end{eqnarray*}

We show the main point:

\smallskip

In $V_1[G]$, $(^\omega 2)^{V[\{\name{\tau}_\xi \such \xi 
\in X \cup Y\}]}$ 
is a Lebesgue null set.

\smallskip

Since   $\exists^{\chi} \alpha \; g_{\chi,\omm}(\alpha) = 
\{ (\gamma,y) \such \chi \cdot \gamma + y  \in Y \}$   
we can fix such an $\alpha \in (\chi \cdot \omm) \setminus
X$
\marginparr{$\alpha$} that is 
not in $A^\chi_{\chi \cdot \gamma + y}$ for $\chi \cdot \gamma + y \in Y$.

\begin{lemma}\label{6.3}
(See \ref{2.8}.)
In $V_1^{P_{\alpha^*}}$, the set
$(^\omega 2)^{V_1[\tau_\xi \such \xi \in X \cup Y]}$ has Lebesgue measure $0$,
and a witness for a definition for a measure zero superset
can be found in $V^{P_{\alpha+1}}$  for  
$\alpha \in \chi\setminus X$  that is 
not in $E_\xi$ for every $\xi \in Y-\chi$.
\end{lemma}

Now proceed through the analogues of Sections~\ref{S2} and~\ref{S3}.
In the definition of a blueprint
we allow $\m^t$ and $\n^t$ to indicate in which intervals
$[\chi \cdot \gamma, \chi \cdot (\gamma +1))$ the
heart of the delta system (intersected with 
the Cohen parts for $\m^t$)
lies, hence $\m^t, \n^t \in [\omm]^{<\om}$ and $\m^t \subseteq \n^t$
in general not as an initial segment, but inserted according to
the type of the heart.
(The old $\n^t$ 
would be just the length of our new $\n^t$.)

\smallskip

Then we modify \ref{4.2} as follows:
In (d) 2.\ we say $n < |\n^t|$ and in (d) 4.\ we say

\smallskip

if  $n < \dom(\m^t) \Leftrightarrow
\forall \ell 
(\alpha_\ell \in [\chi \cdot \m^t(n), \chi \cdot \m^t(n) + \mu) )
\Leftrightarrow
\exists \ell 
(\alpha_\ell \in [\chi \cdot \m^t(n), \chi \cdot \m^t(n) + \mu ))$,
and 

if  $n < \dom(\n^t) \setminus \dom(\m^t) \Leftrightarrow
\forall \ell 
(\alpha_\ell \in [\chi \cdot \m^t(n) +\mu, \chi \cdot (\m^t(n) + 1) )
\Leftrightarrow
\exists \ell 
(\alpha_\ell \in [\chi \cdot \m^t(n) + \mu, \chi \cdot (\m^t(n) + 1) )$.


\smallskip

The rest of Section \ref{S4} shows that the
new ${\mathcal K}^3$ has the desired members.
In \ref{5.3}, the choice of the blueprint has to 
be modified accordingly.
Thus we get $(**)_{\bar{Q}}$ for the modified class ${\mathcal K}^3$.

\smallskip

Since the analogue of \ref{2.7} holds, we
also get analogues to \ref{5.4}
and to \ref{5.5} and hence can finish the proof of \ref{6.1}
\proofendof{\ref{6.1}}


\section{Getting the Premises of \ref{1.1} and \ref{2.1}}
\label{S7}

In this section we discuss how to get the bare set-theoretic
premises of Theorems~\ref{1.2} and \ref{2.1}.

\smallskip

If we do not insist on $(V_1, V_2)$ having the same cardinals but
just require $(^\om V_1)^{V_2} \subseteq V_1$, then we can get the
situation in the premise of \ref{1.2} for example as follows:

Take for $V_1$ any model of $\zfc$ and 
let $\aleph_1 \leq \nu < \nu'$ be regular cardinals in $V_1$.
We extend $V_1$ by forcing with $P=
(\{ f \such  f \colon \nu \to \nu', |\dom(f)|^{V_1}
 \leq \aleph_0\}, \subseteq)$.
Since $P$ is $\om$-closed we have that $(^\om V_1)^{V_2} \subseteq V_1$.
We set
$$
N= \{(\mu,\mu') \in V_1 \such \exists f \in V_2 \: f\colon
\mu \stackrel{\mbox{\footnotesize cofinal}}{\to} \mu', 
\mu, \mu' \mbox{ regular in } V_1,
\mu < \mu'\}.
$$
Let $\lambda = \min (\pi_0(N))$, where $\pi_0$ denotes the projection onto
the first coordinate. Then we have that $\cf^{V_2}(\lambda)$ is uncountable.
Let $\mu' = \mu'(\lambda)$ a minimal witness that $\lambda \in 
\pi_0(N)$ and let $f \in V_2$, $f \colon \lambda 
 \stackrel{\mbox{\footnotesize cofinal}}{\to}
\mu'$. Let $C = \rge(f) \in V_2$. Then $|C|^{V_2} = |\lambda|^{V_2}
=\lambda$. Let ${\mathcal I} \in V_2$ be the set of all bounded subsets of $C$.
For any $B \in V_1$ such that $|B|^{V_1} < \mu'$ we have that
$B \cap C$ is not cofinal in $\mu'$.

\smallskip

If we allow cofinalities to be changed, there is the following 
constellation with consistency strength $\exists \kappa \; o(\kappa) =
\omega_1$:
Gitik \cite{Gitik} shows that assuming 
 $\exists \kappa \; o(\kappa) = \omega_1$ there is some $V$ 
(got with a preparatory forcing) such that in $V$, 
there is a regular cardinal $\kappa > \omega_1$ and 
a notion of forcing $P$ that adds a cofinal sequence of length
$\omega_1$ to $\kappa$ and does not add any countable sequences and
does not add any bounded subsets of $\kappa$.
Now we have $V_1=V$, $V_2 = V^P$, $C =$ the range of the new
cofinal sequence, $\mu= \kappa$, $\lambda = \aleph_1$,
${\mathcal I} = \{ C' \subseteq \kappa \such C' \in V_2, |C'| < \aleph_1 \}$.

\medskip

In order to get  $(V_1,V_2)$ with the same cofinality function, we
take a model announced in the ``Added in proof'' in
 Gitik \cite{Gitik89}:

\begin{theorem}\label{7.1}[Gitik]
Assume that there is a measurable $\kappa$ of Mitchell order $\kappa^{++} + 
\theta$, $\theta$ regular and $\theta \geq \omega_1$.
Then the singular cardinal hypothesis can be violated in the following manner:
There is some model $V$ such that $2^\kappa = \kappa^+$ in $V$ and 
such that there is a notion
of forcing $P$ such that $P$ does not change cofinalities above 
$\kappa$ and such that in $V^P$, 
$\kappa$ is a singular strong limit, $\aleph_0 <
\cf(\kappa)= \theta$, $2^\kappa = \kappa^{++}$ and such that
$\forall x (x \in V^P \wedge x \subseteq {\rm Ord} \wedge
  |x|^{V^P} < \kappa^+ 
\rightarrow \exists y \in V (y \in {\rm Ord} \wedge
|y|^{V^P} < \kappa^+ \wedge
x \subseteq y))$. \proofend
\end{theorem}

\begin{remark}
 By \cite{GitikMitchell96} the lower bound for the consistency strength
is of such a failure of $\sch$ is between 
$\exists \kappa \; o(\kappa) = \kappa^{++}$ and 
$\exists \kappa \; o(\kappa) = \kappa^{++}+\theta$, and if $\theta 
> \aleph_1$ then the strength is $o(\kappa) = \kappa^{++}+\theta$.
\end{remark}

\begin{theorem}\label{7.2} 
Suppose that we have $V$, $P$, $\kappa$, $\theta$ 
as in Theorem~\ref{7.1}.

Then there are $V_1$, $V_2$ such that
 \begin{enumerate}
\item $V \subseteq V_1 \subseteq V_2 \subseteq V[G]$,
\item $(H(\kappa))^{V_1} =(H(\kappa))^{V_2} = (H(\kappa))^{V[G]}$, 
\item $(^{<\theta} V_1)^{V_2} \subseteq V_1$,
\item $V_1$ and $V_2$ have the same cofinality function,
\item in $V_2$ there is a subset $C$ of $\kappa$ of size $\theta$ such that
$C$ is not covered by any set in $V_1$ of size less than $\kappa$.
\end{enumerate}
\end{theorem}
\proof
Let $A = H(\kappa)^{V[G]}$.

\smallskip

By the ``cov versus pp (= pseudo power)  theorem''
\cite[II, 5.4]{Sh:g}  we get that ${\rm pp}(\kappa) = 2^\kappa
= \kappa^{++}$ in $V_2$, and hence by
the definition of pp there is a $\langle \kappa_i \such i < \theta \rangle
\in V[G]$ be a sequence of regular cardinals cofinal in $\kappa$
and an ideal $I $ on $\theta$ containing all the bounded sets in $\theta$ such 
that $\tcf(\prodt \kappa_i/ I) = \kappa^{++}$.
That means:
There is a $<_I$-cofinal
scale $\langle f_\alpha \such \alpha \in \kappa^{++} \rangle$ in $V_2$, i.e.\
for $\alpha < \beta \in \kappa^{++}$ we have
\begin{equation*}\begin{split}
&f_\alpha \colon  \theta \to \kappa,\\
&f_\alpha(\gamma) \in \kappa_\gamma\mbox{ for }\gamma \in \theta,\\
& f_\alpha <_I f_\beta \mbox{ for } \alpha < \beta \in \kappa^{++}\\
&\forall g \in \prodt_{i \in \theta} \kappa_i \; \exists \alpha
\in \kappa^{++} \:\: g <_I f_\alpha,
\end{split}\end{equation*}
where $f <_I g$ iff $\{ i < \theta \such f(i) \geq g(i) \}\in I$.
(By \cite[VIII, \S1]{Sh:g} that there is even a scale with respect to the
ideal $J_\theta^{bd}$ of the bounded subsets of $\theta$.)

We set
$$
V_1  = V[A, \langle \kappa_i \such i < \theta \rangle].
$$
Then we have that there is some $f_\alpha \in V^P$ that $<_I$-dominates
$V_1$:

Proof: In $V$,
 in the subalgebra $P'$ of the Gitik algebra $P$ that is generated by
$H(\kappa)^{V[G]} \cup \{ \langle\kappa_i \such i < \theta \rangle \}$
there are only $\leq \kappa^+$ elements (since
the Gitik algebra $P$ hat the $\kappa^+$-c.c.)
 and it has the $\kappa^+$ c.c.
Hence there are only $\kappa^+$ many $P'$-names for
subsets of $\kappa$ in $V$, so
we have that in $V_1 = V^{P'}$,
 $2^\kappa = \kappa^+$.

Since $C_\alpha =
\{ f \in \, ^\theta \!\kappa \cap V_1 \such f \not\leq_I f_\alpha \}$
is decreasing, of length $\kappa^{++}$ and has empty intersection, there
is some $\alpha < \kappa^{++}$  such that $C_\alpha =
\emptyset$ 
and hence $f_\alpha$ that $<_I$-dominates 
$^\theta \kappa \cap V_1$.

\smallskip

We fix such an $f_\alpha$ and set
$$
V_2 = V_1[f_\alpha].
$$
For $C$ we take $\rge(f_\alpha)$.
Now all the items claimed in \ref{7.2} are true:

We give a proof of item 5, the others are easier.
We show that $\rge(f_\alpha) = C$ is a set in 
$V_2$ that is not covered
by any set $B$ in $V_1$ of size less than $\kappa$.

Suppose the contrary: $B \supseteq C$, $B \in V_1$ and $|B| < \kappa$.
We show that these premises imply $f_\alpha \in V_1$.
We have that $\langle \sup(B \cap \kappa_{i}) \such 
i < \theta \rangle \in V_1$. Since $|B|<\kappa$, 
there is some 
$\theta_0 < \theta$ such that for $i >  \theta_0$
 we have that $\sup(B\cap \kappa_{i}) <\kappa_{i}$. 

We set 
$$g(i) = \left\{ \begin{array}{ll}
\sup(B\cap \kappa_{i}) +1, & \mbox{ if } i > \theta_0, \\
                                0, & \mbox{ else.}
\end{array}
\right.
$$
But we have that $f_\alpha(\gamma) < g(\gamma)$ for $\gamma > \theta_0$. 
Since that latter is in $V_1$ and since $I$ contains all the bounded subsets
of $\theta$ and is proper, 
this is a contradiction to $f_\alpha $ 
being $<_I$-unbounded
 and hence to being $<_I$-dominating over $V_1$.

Remark: Unboundedness with respect to $<_I$  
instead of being dominating w.r.t.\
$<_I$
would suffice for the proof of item 5 and all other items. 
\proofendof{\ref{7.2}}

\end{document}